\documentclass[11pt,twoside]{article}

\usepackage{latexsym,amssymb,amsfonts,amsmath,amsthm}

\usepackage{mathrsfs}

\allowdisplaybreaks[3]

\newcommand{\listinit}{
  \setlength{\labelwidth}{0.9cm}
  \setlength{\labelsep}{0.2cm}
  \setlength{\rightmargin}{0.5cm}
  \setlength{\leftmargin}{1.5cm}
  \setlength{\itemindent}{0cm}}
\newcommand{\hf}{\hfill}

\theoremstyle{plain}
\newtheorem{thm}{Theorem}[section]
\newtheorem{cor}[thm]{Corollary}
\newtheorem{lem}[thm]{Lemma}
\newtheorem{prop}[thm]{Proposition}

\theoremstyle{definition}
\newtheorem{defn}[thm]{Definition}

\theoremstyle{remark}
\newtheorem{rem}[thm]{Remark}

\newtheorem{ex}[thm]{Example}

\newenvironment{claim}[1][]{\smallskip \noindent\emph{#1}\ }{\smallskip}

\newcommand{\vtsp}{\hspace{0.1em}}
\newcommand{\negvtsp}{\hspace{-0.1em}}

\newlength{\phantomheight}

\DeclareMathOperator{\Ker}{Ker}

\DeclareMathOperator{\Ran}{Ran}
\DeclareMathOperator{\Index}{Index}
\DeclareMathOperator{\re}{Re}
\DeclareMathOperator{\im}{Im}
\DeclareMathOperator{\codim}{codim}

\DeclareMathOperator{\Span}{Sp}

\DeclareMathOperator{\diag}{diag}

\DeclareMathOperator{\Closure}{Cl}

\newcommand{\loc}{\mathrm{loc}}

\newcommand{\supp}{\mathrm{supp}}

\newcommand{\R}{\mathbb{R}}
\newcommand{\N}{\mathbb{N}}
\newcommand{\NZ}{\mathbb{N}_0}
\newcommand{\C}{\mathbb{C}}
\newcommand{\Z}{\mathbb{Z}}
\renewcommand{\epsilon}{\varepsilon}

\newcommand{\ipd}[2]{\langle{#1},{#2}\rangle}
\newcommand{\bigipd}[2]{\bigl\langle{#1},\vtsp{#2}\bigr\rangle}

\newcommand{\abs}[1]{\lvert{#1}\rvert}

\newcommand{\bigabs}[1]{\bigl\lvert{#1}\bigr\rvert}
\newcommand{\norm}[1]{\lVert{#1}\rVert}

\newcommand{\bignorm}[1]{\bigl\lVert{#1}\bigr\rVert}

\renewcommand{\d}{\,d}

\usepackage{ifthen}

\newcommand{\cly}{{\Pi^n}}
\newcommand{\don}{{\R^n_*}}

\newcommand{\dfo}{D}

\DeclareMathOperator{\order}{ord}
\newcommand{\sumord}[2]{[#1-#2]}
\newcommand{\detsym}[1]{\mathrm{det}_{#1}}

\newcommand{\Sob}[2]{\ifthenelse{\equal{0}{#2}}{L^{#1}}%
  {\ifthenelse{\equal{2}{#1}}{H^{#2}}{H^{#1,#2}}}}
\newcommand{\Hol}[2]{C^{#1+#2}}
\newcommand{\Zyg}[1]{\mathcal C^{#1}}
\newcommand{\Cnt}[1]{C^{#1}}

\newcommand{\negvvtsp}{\hspace{-0.05em}}
\newcommand{\clysp}[2]{\mathsf{Z}_{#1}{\negvvtsp #2}}
\newcommand{\consp}[2]{\mathsf{X}_{#1}{\negvvtsp #2}}
\newcommand{\donsp}[2]{\mathsf{Y}_{\negvtsp #1}{\negvvtsp #2}}

\newcommand{\clyvsp}[3]{\mathsf{Z}_{#1}^{#2\!}{#3}}
\newcommand{\convsp}[3]{\mathsf{X}_{#1}^{#2\!}{#3}}
\newcommand{\donvsp}[3]{\mathsf{Y}_{\! #1}^{#2\!}{#3}}

\newcommand{\ssym}[1]{\mathrm{Sc}^{#1}}
\newcommand{\psym}[1]{\mathrm{Sym}^{#1}}
\newcommand{\genXYfns}{B\negvtsp S}
\newcommand{\XYfns}{B\negvtsp S_{01}}
\newcommand{\pdo}[1]{\Psi\mathrm{Op}^{#1}}

\newcommand{\fspec}[1]{\sigma(#1)}
\newcommand{\Res}[1]{\Gamma(#1)}

\title{Fredholm Properties of Elliptic Operators on $\R^n$}
\author{Daniel M.~Elton}

\begin{document}

\maketitle

\begin{abstract}
We study the Fredholm properties of a general class of elliptic differential operators on $\R^n$. These results are expressed in terms of a class of weighted function spaces, which can be locally modeled on a wide variety of standard function spaces, and a related spectral pencil problem on the sphere, which is defined in terms of the asymptotic behaviour of the coefficients of the original operator. 
\end{abstract}

\section{Introduction}

\label{intro}

In this paper we consider systems of differential operators on $\R^n$ whose coefficients have certain asymptotic properties as $\abs x\to\infty$. These elliptic operators define continuous maps between weighted function spaces on $\R^n$ which are locally modeled on Sobolev, H\"older or other types of spaces, but which contain a derivative dependent weight that controls behaviour as $\abs x\to\infty$. In this setting we obtain results of the following type which parallel those of the standard theory for elliptic operators on compact manifolds;
\begin{itemize}
\item A priori estimates for solutions.
\item Regularity results relating solutions in different weighted spaces.
\item The Fredholm property for operators acting between certain weighted spaces.
\item Dependence of the Fredholm index on the weighted spaces.
\end{itemize}

In order to avoid a large number of lengthy definitions in the introduction we presently restrict our attention to the formulation of our results for the case of scalar operators acting on weighted function spaces modeled on Sobolev spaces of integral order. At the end of the introduction we indicate where the corresponding results for the general case can be found. 

\bigskip

For the purpose of defining our class of elliptic operators we introduce the following symbol classes (see Section \ref{notation} for basic notation).

\begin{defn}
\label{defofpsym}
For any $\beta\in\R$ define $\ssym\beta$ to be the set of those functions $p\in C^\infty_\loc$ for which there exists $a\in C^\infty(S^{n-1})$ and $q\in C^\infty_\loc$ satisfying the following conditions.
\begin{list}{}{\listinit}
\item[(i)\hf]
Writing $x\in\don$ in polar coordinates as $x=(r,\omega)$ we have 
\[
p(x)\;=\;a(\omega)r^{-\beta}+q(x)
\]
whenever $\abs x=r\ge1$.
\item[(ii)\hf]
For any multi-index $\alpha$ we have an estimate of the form
\[
\dfo_x^\alpha q(x)
\,=\,o(\abs x^{-\beta-\abs\alpha})
\quad\text{as $\abs x\to\infty$}.
\]
\end{list}
The function $a(\omega)r^{-\beta}$ defined on $\don$ will be called the \emph{principal part} of $p$.
\end{defn}

Let $A(x,\dfo_x)$ be a differential operator on $\R^n$ of order $m$. Thus we can write
\begin{equation}
\label{ftwhadop}
A(x,\dfo_x)
\;=\;\sum_{\abs\alpha\le m}p^\alpha(x)\,\dfo_x^\alpha
\end{equation}
where $p^\alpha\in C^\infty_\loc$ for each multi-index $\alpha$ with $\abs\alpha\le m$. We say that $A$ is an \emph{admissible elliptic operator} provided $p^\alpha\in\ssym{m-\abs\alpha}$ for each $\abs\alpha\le m$ and $A$ is uniformly elliptic on $\R^n$ in the sense that 
\begin{equation}
\label{introunifell}
\biggl\lvert\sum_{\abs\alpha=m}p^\alpha(x)\xi^\alpha\biggr\rvert
\;\ge\;C\abs\xi^m
\end{equation}
for all $x,\xi\in\R^n$ (where $C$ is some positive constant). Examples of admissible elliptic operators include the Laplacian $-\Delta$ and the Schr\"odinger operator $-\Delta+V$ whenever $V\in\ssym2$; in particular $V$ must decay at least as quickly as $\abs x^{-2}$.

In order to define the function spaces on which $A$ shall act we need to introduce the weight function $\Lambda$ defined on $\R^n$ by $\Lambda(x)=(1+\abs x^2)^{1/2}$.

\begin{defn}
\label{defnofHpkbeta1}
For $p\in[1,\infty)$, $k\in\NZ$ and $\beta\in\R$ define a norm $\norm{\cdot}_{H^{p,k}_\beta}$ on $C^\infty_0$ by
\[
\norm{u}_{H^{p,k}_\beta}^{\,p}
\;=\;\sum_{\abs\alpha\le k}\int\Lambda^{p(\beta+\abs{\alpha})}(x)\,\abs{\dfo_x^\alpha u(x)}^p\d^nx,
\]
and let $H^{p,k}_\beta$ denote the Banach space obtained by taking the completion of $C^\infty_0$ with respect to this norm. 

For $p\in(1,\infty)$, $k\in\Z\!\setminus\!\NZ$ and $\beta\in\R$ let $q\in(1,\infty)$ be given by $1/p+1/q=1$ and define $H^{p,k}_\beta$ to be the Banach space obtained by taking the dual of $H^{q,-k}_{-\beta}$ with respect to the $L^2$ pairing on $\R^n$.
\end{defn}

If $A$ is an admissible elliptic operator of order $m$ then $A$ defines a continuous map $H^{p,k+m}_{\beta-m}\to H^{p,k}_\beta$ for any $p\in(1,\infty)$, $k\in\NZ$ and $\beta\in\R$. We obtain the following regularity result relating solutions of the equation $Au=f$ for some different values of $p$, $k$ and $\beta$.

\begin{thm}
\label{introthm1}
Let $p,q\in(1,\infty)$, $k,l\in\Z$, $\beta,\gamma\in\R$ and suppose we either have $\beta+n/p<\gamma+n/q$ or $\beta+n/p\le\gamma+n/q$ and $p\ge q$. If $Au\in H^{p,k}_\beta\cap H^{q,l}_\gamma$ for some $u\in H^{q,l+m}_{\gamma-m}$ then we also have $u\in H^{p,k+m}_{\beta-m}$. Furthermore
\[
\norm{u}_{H^{p,k+m}_{\beta-m}}
\;\le\;C\bigl(\norm{Au}_{H^{p,k}_\beta}+\norm{u}_{H^{q,l+m}_{\gamma-m}}\bigr)
\]
for all such $u$.
\end{thm}

In order to proceed with further regularity results and Fredholm properties for the map $A:H^{p,k+m}_{\beta-m}\to H^{p,k}_\beta$ we must eliminate a countable set of values of $\beta$. These values are related to the eigenvalues of an associated spectral problem which we now introduce.

Suppose $A$ is an admissible elliptic operator of order $m$ given by \eqref{ftwhadop}. We define the \emph{principal part} of $A$ to be the operator $A_0$ on $\don$ given by
\[
A_0(x,\dfo_x)
\;=\;\sum_{\abs\alpha\le m} a^\alpha(\omega)r^{\abs\alpha-m}\,\dfo_x^\alpha,
\]
where, for each $\abs\alpha\le m$, $a^\alpha(\omega)r^{\abs\alpha-m}$ is the principal part of $p^\alpha$. It is easy to see that the ellipticity  estimate \eqref{introunifell} for $A$ implies $A_0$ is elliptic on $\don$.

The principal part of $A$ can be rewritten in the form 
\[
A_0(x,\dfo_x)
\;=\;\sum_{j=0}^m A^{m-j}(\omega,\dfo_\omega)\vtsp(r\dfo_r)^j\bigl(r^{-m}\,\cdot\,\bigr),
\]
where, for $j=0,\dots,m$, $A^j(\omega,\dfo_\omega)$ is a differential operator on $S^{n-1}$ of order at most~$j$. Associated to $A$ we now define an operator pencil $\:\mathfrak B_A\negvtsp:\C\to\mathscr L\bigl(\Sob 2m(S^{n-1}),\vtsp\Sob 20(S^{n-1})\bigr)$ by
\[
\mathfrak B_A(\lambda)
\;=\;\sum_{j=0}^m A^{m-j}(\omega,\dfo_\omega)\vtsp\lambda^j
\]
for each $\lambda\in\C$. The spectrum of this operator pencil is the set 
\[
\fspec{\mathfrak B_A}
\;=\;\bigl\{\lambda\in\C\;\big|\;
\mathfrak B_A(\lambda):\text{$\Sob 2m(S^{n-1})\to\Sob 20(S^{n-1})$ is not invertible}\bigr\}.
\]
The geometric multiplicity of $\lambda_0\in\fspec{\mathfrak B_A}$ is just $\dim\Ker\mathfrak B_A(\lambda_0)$, whilst the algebraic multiplicity can be defined as the sum of the lengths of a set of maximal Jordan chains corresponding to $\lambda_0$ (see Section \ref{oppensec} or \cite{GGK} for more details).

Using the ellipticity of $A_0$ it can be shown that $\fspec{\mathfrak B_A}$ consists of isolated points of finite algebraic multiplicity and that any strip of finite width parallel to the real axis contains at most finitely many points of $\fspec{\mathfrak B_A}$ (see Theorem 5.2.1 in \cite{KMR} or Theorem 1.2.1 in \cite{NP2} for example). 

The projection of $\fspec{\mathfrak B_A}$ onto the imaginary axis is of particular importance and will be denoted by $\Res{A}$; that is 
\[
\Res{A}
\;=\;\bigl\{\im\lambda\,\big|\,\lambda\in\fspec{\mathfrak B_A}\bigr\}
\;\subset\;\R.
\]
In particular, the above discussion implies $\Res{A}$ consists of isolated points and, given $\gamma\in\Res{A}$, the total algebraic multiplicity of all those $\lambda\in\fspec{\mathfrak B_A}$ with $\im\lambda=\gamma$ is finite. 

The Fredholm property for $A$ is related to the spectrum of the associated operator pencil through the set $\Res{A}$ as follows.

\begin{thm}
\label{introthm2}
Let $p\in(1,\infty)$, $k\in\Z$ and $\beta\in\R$. If $\beta+n/p\notin\Res{A}$ then the map $A:H^{p,k+m}_{\beta-m}\to H^{p,k}_\beta$ is Fredholm.
\end{thm}

We also obtain a $\Res{A}$ dependent regularity result complementing Theorem \ref{introthm1}.

\begin{thm}
\label{introthm3}
Let $p,q\in(1,\infty)$, $k,l\in\Z$, $\beta,\gamma\in\R$ and suppose $\beta+n/p$ and $\gamma+n/q$ belong to the same component of $\R\!\setminus\!\Res{A}$. If $Au\in H^{p,k}_\beta\cap H^{q,l}_\gamma$ for some $u\in H^{q,l+m}_{\gamma-m}$ then we also have $u\in H^{p,k+m}_{\beta-m}$. Furthermore
\[
\norm{u}_{H^{p,k+m}_{\beta-m}}
\;\le\;C\bigl(\norm{Au}_{H^{p,k}_\beta}+\norm{u}_{H^{q,l+m}_{\gamma-m}}\bigr)
\]
for all such $u$.
\end{thm}

As a consequence of Theorems \ref{introthm2} and \ref{introthm3} we also obtain a stability result for the Fredholm index of $A$.

\begin{thm}
\label{introthm4}
Let $p,q\in(1,\infty)$, $k,l\in\Z$, $\beta,\gamma\in\R$ and suppose $\beta+n/p$ and $\gamma+n/q$ belong to the same component of $\R\!\setminus\!\Res{A}$. Then the Fredholm maps $A:H^{p,k+m}_{\beta-m}\to H^{p,k}_\beta$ and $A:H^{q,l+m}_{\gamma-m}\to H^{q,l}_\gamma$ have the same index.
\end{thm}

If the parameter $\beta$ is varied so that $\beta+n/p$ moves between components of $\R\!\setminus\!\Res{A}$ then the index of the corresponding map will change. This change is related to more detailed information about the spectrum of the operator pencil $\mathfrak B_A$.

\begin{thm}
\label{introthm5}
Let $p\in(1,\infty)$, $k\in\Z$ and $\beta_1,\beta_2\in\R$ with $\beta_1\le\beta_2$ and $\beta_i+n/p\notin\Res{A}$ for $i=1,2$. Set $\Sigma=\{\lambda\in\fspec{\mathfrak B}\,|\,\im\lambda\in[\beta_1,\beta_2]\}$ and, for each $\lambda\in\Sigma$, let $m_\lambda$ denote the algebraic multiplicity of $\lambda$. Then we have 
\[
\Index A^{(\beta_1)}
\;=\;\Index A^{(\beta_2)}+\sum_{\lambda\in\Sigma}m_\lambda,
\]
where $A^{(\beta_i)}$ denotes the map $A:H^{p,k+m}_{\beta_i-m}\to H^{p,k}_{\beta_i}$ for $i=1,2$.
\end{thm}

\medskip

Let $A^*(x,\dfo)$ denote the differential operator obtained by taking the formal adjoint of $A(x,\dfo_x)$ (with respect to the standard Lebesgue measure on $\R^n)$. Using the definition of the symbol classes $\ssym\beta$ it is straightforward to check that $A^*$ is also an admissible elliptic operator. Furthermore we have 
\begin{equation}
\label{flistofrelofspecAAs}
\Res{A^*}
\;=\;(n+m)-\Res{A}.
\end{equation}

In the case that $A$ is formally self-adjoint we can determine the Fredholm index entirely from knowledge of the spectrum of the operator pencil $\mathfrak B_A$.

\begin{thm}
\label{introthm6}
Let $p\in(1,\infty)$, $k\in\Z$ and, for any $\gamma\in\R$, let $A^{(\gamma)}$ denote the map $A:H^{p,k+m}_{\gamma-m}\to H^{p,k}_\gamma$. Now suppose $A$ is formally self-adjoint. Then $\Res{A}$ is symmetric about $(n+m)/2$ and, for any $\beta\in\R$ with $\beta+n/p\notin\Res{A}$, we have 
\[
\Index A^{(n+m-\beta-2n/p)}
\;=\;-\Index A^{(\beta)}.
\]
In particular we either have $(n+m)/2\notin\Res{A}$, in which case $\Index A^{(m/2+n/2-n/p)}=0$, or we have $(n+m)/2\in\Res{A}$, in which case the sum of the algebraic multiplicities of those $\lambda\in\fspec{\mathfrak B_A}$ with $\im\lambda=(n+m)/2$ is even (say $2d$ for some $d\in\N$) and 
\[
\Index A^{(m/2+n/2-n/p-\epsilon)}
\;=\;d
\;=\;-\Index A^{(m/2+n/2-n/p+\epsilon)}
\]
for all sufficiently small $\epsilon>0$.
\end{thm}

\bigskip

The paper is arranged as follows. In Section \ref{funspace} we define the general class of weighted function spaces on which our elliptic operators act, as well as related weighted function spaces for the associated ``model operators''. The majority of this Section is devoted to establishing the basic properties of these spaces that are necessary in order to work with them. In particular Section \ref{equivnormssssec} gives details of how some previously defined weighted function spaces (including the weighted Sobolev spaces of Definition \ref{defnofHpkbeta1}) arise in this general setting. 

The full class of elliptic operators on $\R^n$ to which our results apply is introduced at the beginning to Section \ref{mainressec} (see Definition \ref{defnofadellop}). This class is basically a generalisation of the class of (scalar) admissible elliptic operators introduced above to cover the case of systems with Douglis-Nirenberg type ellipticity. The generalisations of Theorems \ref{introthm1} to \ref{introthm6} are given in Theorems \ref{regnorest}, \ref{sfomres1}, \ref{changebetareggen}, \ref{stabofindex}, \ref{mressec5} and \ref{selfadjthm} respectively (Remarks \ref{locandbndincegs} and \ref{actidentifofequivnorms} provide the details needed to derive the results given above from their counterparts in Section \ref{mainressec}). Additionally it is shown that the finite dimensionality of the kernel implied by Theorem \ref{introthm2} (or Theorem \ref{sfomres1}) remains valid without restriction on the parameter $\beta$ (see Theorem \ref{genfdkerresult}). Finally, at the end of Section \ref{mainressec}, we give some index formulae for elliptic operators whose principal part is homogeneous with constant coefficients (see Theorems \ref{thethmihaflofse} and \ref{corofthethmihaflofse}).

The main results are established from results for ``model operators'' on $\cly$ and $\don$. These operators provide the necessary generalisation of the operators $A_0$ and $\mathfrak B_A$ introduced above and are dealt with in Section \ref{modprobsec}. The isomorphism results contained in Sections \ref{subsecisooncly} and \ref{seconmodprobondon} (for $\cly$ and $\don$ respectively) are also of interest in their own right.

\medskip

There appear to be at least two independent areas of research related to the results presented here. One of these areas has involved the detailed study of homogeneous constant coefficient elliptic operators on $\R^n$. Theorem \ref{thethmihaflofse} was established for the spaces $E=\Sob pk$, $p\in(1,\infty)$, $k\in\NZ$ in a series of papers culminating with \cite{LM}. In \cite{Be} the same Theorem was given for scalar operators and the spaces $E=\Hol l\sigma$, $l\in\NZ$, $\sigma\in(0,1)$. The approach used in these papers is based on establishing mapping properties for an explicit class of related convolution operators. This approach seems to be well suited to problems involving constant coefficient operators and allows for the computation of the index and the characterisation of kernels and cokernels; on the other hand, this approach does not appear to generalise easily to cover operators with variable coefficients. 

Operators with variable coefficients have been considered in another area of research. In \cite{BK} Theorem \ref{sfomres1} was established for scalar operators and the spaces $E=\Sob 2k$, $k\in\NZ$, whilst Theorem \ref{genfdkerresult} was given for scalar operators and the spaces $E=\Sob p0$, $p\in(1,\infty)$ when $\beta=n/p$. The possibility of generalisation to systems of operators was also observed in \cite{BK}. Our method for obtaining Theorem \ref{sfomres1} from results for the ``model operators'' is essentially similar to the method used in \cite{BK} with the exception that our consideration of more general classes of weighted function spaces allows us to use a duality argument in going from the semi-Fredholm property to the Fredholm property (i.e.~from Theorem \ref{wfomres1} to Theorem \ref{sfomres1}; a regulariser is used for the corresponding step in \cite{BK}).

Since the paper \cite{BK} work on operators with variable coefficients appears to have been concentrated on the closely related problem of elliptic operators on domains whose boundaries contain conical singularities. For this related problem many results similar to those of Section \ref{mainressec} have been established in the case of weighted $L^p$ Sobolev spaces of integer order and weighted H\"older spaces (see \cite{KMR} or \cite{NP2} for an exposition of this theory). The presumed existence of parallel results for elliptic operators on $\R^n$ has been remarked upon by several authors (see in particular Remark 4.1.5 in \cite{NP2}) --- however a detailed study of this problem does not appear to have been carried out as of yet.

The problem of elliptic operators on $\R^n$ and that of elliptic operators on domains with boundaries containing conical singularities both give rise to the same type of ``model operators'' on cylindrical and conical domains (which, in our case, are $\cly$ and $\don$ respectively). Furthermore the arguments used to go from results for the ``model operators'' to results for the actual problems are similar in many aspects. Important developments in the study of the relevant ``model problems'' include \cite{Ko} where the case of scalar operators on weighted $L^2$ Sobolev spaces of non-negative integral order was considered, and \cite{NP1} where generalisations were made to systems of operators on weighted $L^p$ Sobolev spaces of integer order and weighted H\"older spaces. Weighted $L^2$ Sobolev spaces of fractional order were considered in \cite{D} (see \cite{KMR} for more details about the historical development of this theory and references to further results). These results are obtained in Section \ref{modprobsec} using a different approach based on a characterisation of pseudo-differential operators given in \cite{B}. The advantage to this approach is that it enables the systematic consideration of a much more general class of weighted function spaces; perhaps the most important of the new types of spaces to be considered include weighted $L^p$ Sobolev spaces of arbitrary real order and non-separable weighted H\"older spaces (rather than the separable subspaces of these spaces which were considered previously). It should be possible to extend these methods to give more general results for the case of elliptic operators on domains with boundaries containing conical singularities, although no attempt to do this is made here.

One application of the results presented here is to the study of zero modes (or zero energy bound states) of the Dirac-Weyl operator $\sigma.(\dfo-A)$ on $\R^3$ (here $\sigma$ is the vector of Pauli matrices and $A$ is a real vector potential); this will appear in \cite{E}.

\subsection{Notation}

\label{notation}

In this section we introduce some (not necessarily standard) notation and conventions that will be used throughout the paper.

We define $\NZ=\N\cup\{0\}$ to be the set of non-negative integers; thus a multi-index $\alpha$ is simply an element of $\NZ^n$ with $\abs\alpha:=\alpha_1+\dots+\alpha_n$. 

The sets $\R^n\!\setminus\!0$ and $\R\!\times\!S^{n-1}$ are denoted by $\don$ and $\cly$ respectively. The letter $\omega$ is used to denote a point on $S^{n-1}$, whilst $(r,\omega)$ and $(t,\omega)$ denote polar coordinates on $\don$ and cylindrical coordinates on $\cly$ respectively. When necessary $S^{n-1}$, $\don$ and $\cly$ will be considered as Riemannian manifolds with the obvious choice of Riemannian metrics (i.e.~that given by the standard embedding into $\R^n$ as the unit sphere for $S^{n-1}$, the restriction of the Euclidean metric for $\don$ and the product metric for $\cly$).

For $i\in\{1,\dots,n\}$ we use $\dfo_i$ to denote the differential operator $-i\partial/\partial x_i$ on $\R^n$. By $\dfo_x$ we mean the vector differential operator $(\dfo_1,\dots,\dfo_n)$ whilst $\dfo_r$ and $\dfo_t$ are used to denote the differential operators $-i\partial/\partial r$ and $-i\partial/\partial t$ on $\R^+$ and $\R$ respectively. Finally, the notation $A(\omega,\dfo_\omega)$ is used to mean that $A$ is a differential operator on $S^{n-1}$. 

For any manifold $M$ with volume measure $dM$ let $(\cdot,\cdot)_M$ denote the $L^2$ pairing on $M$; that is $(u,v)_M=\int_M u v\d M$ for all appropriate $u$ and $v$. The associated sesquilinear pairing is then $\ipd uv_M:=(u,\overline{v})_M$. We use these pairings with $M$ given by $\R^n$, $S^{n-1}$, $\don$ or $\cly$. In all cases $dM$ is the volume measure induced by our choice of the Riemannian metric; in particular the volume measures on $\R^n$, $\don$ and $\cly$ are $d^nx$ (the Lebesgue measure), $d^nx=r^{n-1}\d r\d S^{n-1}$ and $dt\d S^{n-1}$ respectively, where $\d S^{n-1}$ is the standard volume measure on $S^{n-1}$ inherited from its embedding into $\R^n$ as the unit sphere.

Let $\mathscr D'$ denote the set of distributions on $\R^n$ and $\mathscr D'(M)$ the set of distributions on an arbitrary manifold $M$ (see Section 6.3 of \cite{Hor1} for further details). The pairing $(u,v)_M$ can be defined for all $u\in C^\infty_0(M)$ and $v\in\mathscr D'(M)$; in particular this pairing can be viewed as a way of identifying elements of $\mathscr D'(M)$ with distributional densities on $M$.

If $E\subset\mathscr D'(M)$ for some manifold $M$ we use $E_\loc$ to denote the set of $u\in\mathscr D'(M)$ with $\phi u\in E$ for all $\phi\in C^\infty_0(M)$. On the other hand, for open $U\subseteq M$, $E(U)$ denotes the set of (equivalence classes of) restrictions of elements $u\in E$ to $U$. If $\phi\in C^\infty_0(U)$ and $u\in E(U)$ we can extend $\phi u$ by 0 outside $U$ to enable us to consider it as an element of $E$.

We use $\Lambda$ to denote the weight function defined on $\R^n$ by $\Lambda(x)=(1+\abs x^2)^{1/2}$. Let $\mathscr S$ and $\mathscr S'$ denote respectively the locally convex spaces of Schwartz class functions on $\R^n$ and its dual, the set of tempered distributions on $\R^n$. The topology of the former is provided by the semi-norms
\[
p_l(u)\;:=\;\sum_{\abs\alpha\le l}\,\sup_{x\in\R^n}\,\Lambda^l(x)\abs{\dfo_x^\alpha u(x)}
\]
for any $l\in\NZ$, whilst we choose the weak dual topology for the latter.

For any $l\in\NZ$ we use $\Cnt l$ to denote the set of $l$-times bounded continuously differentiable functions on $\R^n$, provided with the norm
\[
\norm{u}_{\Cnt l}\;=\;\sum_{\abs\alpha\le l}\,\sup_{x\in\R^n}\,\abs{\dfo_x^\alpha u(x)}.
\]
We put $C^\infty=\bigcap_{l\in\NZ}\Cnt l$ and provide this set with the locally convex topology induced by the collection of semi-norms $\{\norm\cdot_{\Cnt l}\,|\,l\in\NZ\}$. The set of smooth functions on $\R^n$ (without restrictions on growth at infinity) is then denoted by $C^\infty_\loc$. We also use $C^\infty_0$ for the set of smooth functions on $\R^n$ with compact support.

For $p\in[1,\infty]$ and $s\in\R$ we use $\Sob ps$ to denote the Sobolev space on $\R^n$ of ``functions with $s$ $p$-integrable derivatives''. This notation will be simplified to $\Sob 2s$ in the case $p=2$ and $\Sob p0$ in the case $s=0$. Other spaces appearing as examples include the H\"older spaces $\Hol l\sigma$ for $l\in\NZ$ and $\sigma\in[0,1)$ (n.b.~we set $\Hol l0=\Cnt l$) and the Zygmund spaces $\Zyg s$ for $s\in\R^+$. A detailed account of all these spaces can be found in \cite{Tr1}.

Let $\genXYfns$ denote the set of $\R$-valued functions $\zeta\in C^\infty$ which are constant in a neighbourhood of 0 and $\infty$. We also use $\XYfns$ to denote the subset of $\genXYfns$ containing those $\zeta$ with $\zeta=0$ in a neighbourhood of 0 and $\zeta=1$ in a neighbourhood of $\infty$. 

If $\chi_1$ and $\chi_2$ are $\R$-valued functions we write $\chi_1\prec\chi_2$ (or, alternatively, $\chi_2\succ\chi_1$) provided $\chi_2=1$ on $\supp(\chi_1)$. If $\chi_1\prec\chi_2$ it clearly follows that $\chi_1\chi_2=\chi_1$. 

We use $C$ to denote any positive real constant whose exact value is not important but which may depend only on the things it is allowed to in a given problem (i.e.~parameters defining function spaces but not the actual element of the space under consideration etc.). Constants depending on something extra are indicated with appropriate function type notation whilst subscripts are added if we need to keep track of the value of a particular constant (e.g.~$C_1(u)$ etc.).

We use the notation $\kappa(u)\asymp\rho(u)$ to indicate the quantities $\kappa(u)$ and $\rho(u)$ satisfy the inequalities $C\kappa(u)\le\rho(u)\le C\kappa(u)$ for all relevant $u$ (possibly including parameter values). Generally $\kappa$ and $\rho$ will be norms of some description.

\section{Function spaces}

\label{funspace}

In this section we introduce classes of weighted function spaces for our elliptic operators and associated `model' operators. In order to define these spaces and establish the basic results necessary to work with them, we use general constructions and arguments applied to specific `model spaces'. 

\begin{defn}
A \emph{model space} $E$ is a Banach space of functions satisfying the following properties.
\begin{list}{}{\listinit}
\item[(A1)\hf]
We have continuous inclusions $\mathscr S\hookrightarrow E\hookrightarrow\mathscr S'$.
\item[(A2)\hf]
Multiplication defines a continuous bilinear map $C^\infty\times E\to E$. 
\item[(A3)\hf]
The norm $\norm{\cdot}_E$ is translationally invariant.
\item[(A4)\hf]
If $\psi$ is a diffeomorphism on $\R^n$ which is linear outside some compact set then the pull back $\psi^*:\mathscr S'\to\mathscr S'$ restricts to give an isomorphism on $E$.
\end{list}
Throughout this paper the letters $E$, $F$ and $G$ are used to denote model spaces.
\end{defn}

\begin{rem}
Examples of model spaces include the Sobolev spaces $\Sob ps$ for $p\in[1,\infty]$ and $s\in\R$, the H\"older spaces $\Hol k\sigma$ for $k\in\NZ$ and $\sigma\in[0,1)$ and the Zygmund spaces $\Zyg s$ for $s\in\R^+$. Clearly if $E$ is a model space and $E_0$ is the separable subspace of $E$ formed by taking the closure of $\mathscr S$ (or, equivalently, of $C^\infty_0$), then $E_0$ is also a model space. 
\end{rem}

\begin{lem}
\label{useobstobeshted}
Suppose $\psi:U\to V$ is a diffeomorphism between open subsets of $\R^n$ and $\chi\in C^\infty_0(V)$. Then $\norm{\psi^*(\chi u)}_E\asymp\norm{\chi u}_E$ for all $u\in E(V)$.
\end{lem}

\begin{proof}
Any point $x\in U$ has a neighbourhood in $U$ outside of which $\psi$ can be extended linearly to give a diffeomorphism of $\R^n$. Choose an finite collection of such neighbourhoods $\{U_i\}_{i\in I}$ which cover $\psi^{-1}(\supp(\chi))$ and, for each $i\in I$, let $\psi_i:\R^n\to\R^n$ denote a diffeomorphism which is linear outside a compact region and satisfies $\psi_i(x)=\psi(x)$ for all $x\in U_i$. Thus $\{\psi_i(U_i)\}_{i\in I}\cup\{\R^n\!\setminus\!\supp(\chi)\}$ is an open cover of $\R^n$. Choosing any partition of unity $\{\phi_i\}_{i\in I}\cup\{\phi_\infty\}$ subordinate to this cover we clearly have $\sum_{i\in I}\phi_i=1$ on $\supp(\chi)$. Conditions (A2) and (A4) now give 
\begin{multline*}
\norm{\psi^*(\chi u)}_E
\ \le\ \sum_{i\in I}\,\norm{\psi^*(\phi_i\chi u)}_E
\ =\ \sum_{i\in I}\,\norm{\psi_i^*(\phi_i\chi u)}_E\\
\le\ C\sum_{i\in I}\,\norm{\phi_i\chi u}_E
\ \le\ C\vtsp\sup_{i\in I}\vtsp\norm{\phi_i\chi u}_E
\ \le\ C\norm{\chi u}_E
\end{multline*}
for all $u\in E(V)$. Symmetry completes the result.
\end{proof}

\begin{defn}
\label{defnofElver1}
Suppose $E$ is a model space and $l\in\NZ$. We define the space $E^l$ to be the set of all $u\in\mathscr S'$ satisfying $\dfo_x^\alpha u\in E$ for each multi-index $\alpha$ with $\abs\alpha\le l$, equipped with the norm
\[
\norm{u}_{E^l}
\;=\;\sum_{\abs\alpha\le l}\,\norm{\dfo_x^\alpha u}_E.
\]
\end{defn}

\begin{rem}
\label{remnofElver1}
It is straightforward to check that $E^l$ is again a model space. Furthermore the differential operator $\dfo_x^\alpha$ clearly defines a continuous map $E^l\to E$ whenever $\abs\alpha\le l$.
\end{rem}

\begin{defn}
\label{defofsss}
For any model space $E$ we define $E_0$ to be the separable subspace of $E$ obtained by taking the closure of $\mathscr S$ in $E$. 
\end{defn}

\begin{rem}
\label{defofsssrem}
It is straightforward to check that $E_0$ is again a model space. Furthermore $C^\infty_0$ is a dense subset of $\mathscr S$ (see Proposition VI.1.3 in \cite{Y} for example) so $E_0=\Closure(\mathscr S)=\Closure(C^\infty_0)$; in particular $\mathscr S$ is dense in $E$ (i.e.~$E=E_0$) iff $C^\infty_0$ is dense in $E$.
\end{rem}

\subsection{Weighted function spaces on $\boldsymbol{\cly}$}

\label{subsecweightfunspcly}

For any chart $(\psi,U)$ on $S^{n-1}$ we can define a corresponding chart $(\Psi,\R\times U)$ on $\cly$ by setting $\Psi(t,\omega)=(t,\psi(\omega))\in\R\times\psi(U)\subseteq\R^n$ for all $(t,\omega)\in\R\times U$. Now suppose $\{(\psi_i,U_i)\}_{i\in I}$ is a finite atlas for $S^{n-1}$ and let $\{(\Psi_i,\R\times U_i)\}_{i\in I}$ be the corresponding atlas for $\cly$. Choose a partition of unity $\{\chi_i\}_{i\in I}$ which is subordinate to the cover $\{U_i\}_{i\in I}$ of $S^{n-1}$. We also consider $\{\chi_i\}_{i\in I}$ to be a partition of unity subordinate to the cover $\{\R\times U_i\}_{i\in I}$ of $\cly$ by regarding each $\chi_i$ as a function on $\cly$ which is independent of $t$.

\smallskip

Let $E$ and $F$ be model spaces on $\R^n$ and $\R^{n-1}$ respectively. We define $E(\cly)$ to be the set of $u\in\mathscr D'(\cly)$ with $(\Psi_i^{-1})^*(\chi_i u)\in E$ for each $i\in I$. On this set we define a norm 
\begin{equation}
\label{defofEclynobeta}
\norm u_{E(\cly)}
\;=\;\sum_{i\in I}\,\bignorm{(\Psi_i^{-1})^*(\chi_i u)}_E\,.
\end{equation}
We define the normed space $F(S^{n-1})$ in a similar fashion with $\Psi_i$ replaced by $\psi_i$. Standard calculations using Conditions (A2) and (A4) show that $E(\cly)$ and $F(S^{n-1})$ are Banach spaces which are independent of the choice of atlas $\{(\psi_i,U_i)\}_{i\in I}$ and partition of unity $\{\chi_i\}_{i\in I}$ (up to equivalent norms); the next result is a somewhat more general statement of this fact for the space $E(\cly)$.  

\begin{lem}
\label{verydetverofinvarlem}
Suppose $\{(\phi_j,V_j)\}_{j\in J}$ is a finite atlas for $S^{n-1}$ and $\{\zeta_j\}_{j\in J}$ is a collection of functions in $C^\infty(S^{n-1})$ with $\supp(\zeta_j)\subset V_j$ for each $j\in J$ and $\bigabs{\sum_{j\in J}\zeta_j}\ge C>0$ on $S^{n-1}$. Let $\{(\Phi_j,\R\times V_j)\}_{j\in J}$ be the corresponding atlas for $\cly$ and consider $\zeta_j$ as a function on $\cly$ which is independent of $t$. Then $u\in E(\cly)$ iff $u\in\mathscr D'(\cly)$ and $(\Phi_j^{-1})^*(\zeta_j u)\in E$ for each $j\in J$. Furthermore
\[
\norm{u}_{E(\cly)}
\;\asymp\;\sum_{j\in J}\,\bignorm{(\Phi_j^{-1})^*(\zeta_j u)}_E\,.
\]
\end{lem}

\begin{proof}
For each $i\in I$ choose $\chi_i'\in C^\infty_0(U_i)$ with $\chi_i'\succ\chi_i$. Then, for each $j\in J$, set $\widetilde\zeta_{ij}=(\zeta_j\chi_i')\circ\Psi_i^{-1}\in C^\infty$ (here we are considering $\zeta_j\chi_i'$ to be a function on $\cly$ which is independent $t$ and extending $\widetilde\zeta_{ij}$ by 0 outside $\R\times\psi_i(U_i)\,$). For $u\in E(\cly)$ and $j\in J$
\[
(\Phi_j^{-1})^*(\zeta_j u)
\ =\ \sum_{i\in I}\,(\Phi_j^{-1})^*(\zeta_j\chi_i'\vtsp\chi_i u)
\ =\ \sum_{i\in I}\,(\Psi_i\circ\Phi_j^{-1})^*
\bigl(\vtsp\widetilde\zeta_{ij}\vtsp(\Psi_i^{-1})^*(\chi_i u)\bigr).
\]
Lemma \ref{useobstobeshted} (with the obvious minor modification) and Condition (A2) now give $(\Phi_j^{-1})^*(\zeta_j u)\in E$ and 
\[
\bignorm{(\Phi_j^{-1})^*(\zeta_j u)}_E
\ \le\ C\sum_{i\in I}\,\bignorm{(\Psi_i^{-1})^*(\chi_i u)}_E
\ \le\ C\norm{u}_{E(\cly)}.
\]

Our assumption on the $\zeta_j$'s implies $\zeta:=\bigl(\sum_{j\in J}\zeta_j\bigr)^{-1}\in C^\infty(S^{n-1})$. Using the partition of unity on $S^{n-1}$ given by $\{\zeta \zeta_j\}_{j\in J}$, the remainder of the result can be completed with an argument similar to that above.
\end{proof}

\begin{rem}
\label{enlargechirem}
We note the following technically useful consequence of Lemma \ref{verydetverofinvarlem}. For each $i\in I$ choose $\chi'_i\in C^\infty_0(U_i)$ with $\chi'_i\succ\chi_i$. Then
\[
\norm{u}_{E(\cly)}
\;=\;\sum_{i\in I}\,\bignorm{(\Psi_i^{-1})^*(\chi_iu)}_E
\;\asymp\;\sum_{i\in I}\,\bignorm{(\Psi_i^{-1})^*(\chi'_iu)}_E.
\]
for all $u\in E(\cly)$.
\end{rem}

\begin{lem}
\label{cordequivnorm}
Suppose $(\psi,U)$ is a chart for $S^{n-1}$ and let $(\Psi,\R\times U)$ be the corresponding chart for $\cly$. Also suppose $\chi\in C^\infty_0(U)$ (considered as a function on $\cly$ which is independent of $t$) and define $\widetilde\chi=\chi\circ\Psi^{-1}\in C^\infty$. Then $\chi u\in E(\cly)$ iff $(\Psi^{-1})^*(\chi u)\in E$ whilst $\norm{\chi u}_{E(\cly)}\asymp\norm{(\Psi^{-1})^*(\chi u)}_E$. On the other hand $\widetilde\chi v\in E$ iff $\Psi^*(\widetilde\chi v)\in E(\cly)$ whilst $\norm{\widetilde\chi u}_E\asymp\norm{\Psi^*(\widetilde\chi u)}_{E(\cly)}$.
\end{lem}

\begin{proof}
The first part of this result follows easily from Lemma \ref{verydetverofinvarlem} applied to the atlas $\{(\psi,U)\}\cup\{(\psi_i,U_i)\}$ for $S^{n-1}$ and any partition of unity $\{\zeta\}\cup\{\zeta_i\}_{i\in I}$ which is subordinate to the covering $\{U\}\cup\{U_i\!\setminus\!\supp(\chi)\}_{i\in I}$ of $S^{n-1}$ (n.b.~$\zeta\chi=\chi$ and $\zeta_i\chi=0$ for each $i\in I$ in this case). 

Choose $\chi'\in C^\infty_0(U)$ with $\chi'\succ\chi$ and set $\widetilde\chi'=\chi'\circ\Psi^{-1}\in C^\infty$. The second part of the result follows from the first part by taking $u=\Psi^*(\widetilde\chi'v)$.
\end{proof}

\begin{rem}
\label{clyforCinfty}
Suppose $K$ is a locally convex space which satisfies Conditions (A1), (A2) and (A4). By applying \eqref{defofEclynobeta} to individual semi-norms we can obviously define a new locally convex space $K(\cly)\subset\mathscr D'(\cly)$. In particular we shall need the locally convex spaces $\mathscr S(\cly)$, $\mathscr S'(\cly)$ and $C^\infty(\cly)$. It is straightforward to check that Condition (A1) gives us continuous inclusions $\mathscr S(\cly)\hookrightarrow E(\cly)\hookrightarrow\mathscr S'(\cly)$ for any model space $E$. It is also clear that $C^\infty(\cly)=\bigcap_{l\in\NZ} C^l(\cly)$ whilst the topology on $C^\infty(\cly)$ is that induced by the collection of semi-norms $\{\norm{\cdot}_{C^l(\cly)}\vtsp|\vtsp l\in\NZ\}$.
\end{rem}

\begin{lem}
\label{PropA2clyset}
For any model space $E$ multiplication defines a continuous bilinear map $C^\infty(\cly)\times E(\cly)\to E(\cly)$. 
\end{lem}

\begin{proof}
Condition (A2) for $E$ means we can find $l\in\NZ$ such that 
\[
\norm{\phi u}_E
\;\le\;C\norm{\phi}_{\Cnt l}\norm{u}_E
\]
for all $\phi\in C^\infty$ and $u\in E$. With the notation of Remark \ref{enlargechirem} it follows that 
\[
\bignorm{(\Psi_i^{-1})^*(\chi_i\phi u)}_E
\ =\ \bignorm{(\Psi_i^{-1})^*(\chi_i\phi)\,(\Psi_i^{-1})^*(\chi'_i u)}_E
\ \le\ C\bignorm{(\Psi_i^{-1})^*(\chi_i\phi)}_{\Cnt l}\,
\bignorm{(\Psi_i^{-1})^*(\chi'_i u)}_E
\]
for any $\phi\in C^\infty(\cly)$ and $u\in E(\cly)$. Combining this with Remark \ref{enlargechirem} and the fact that $I$ is finite, we then get
\begin{multline*}
\norm{\phi u}_{E(\cly)}
\ \le\ C\sum_{i\in I}\,\bignorm{(\Psi_i^{-1})^*(\chi_i\phi u)}_E\\
\le\ C\sum_{i\in I}\,\bignorm{(\Psi_i^{-1})^*(\chi_i\phi)}_{\Cnt l}\,
\bignorm{(\Psi_i^{-1})^*(\chi'_i\phi)}_E
\ \le\ C\norm{\phi}_{\Cnt l(\cly)}\;\norm{u}_{E(\cly)}
\end{multline*}
for all $\phi\in C^\infty(\cly)$ and $u\in E(\cly)$. This completes the proof.
\end{proof}

\bigskip

Given a model space $E$ and $\beta\in\R$ we define $\clysp\beta E$ to be the set of $u\in\mathscr D'(\cly)$ with $e^{\beta t}u\in E(\cly)$. On this set we define a norm
\begin{equation}
\label{defofnormonclysp}
\norm u_{\clysp\beta E}
\,=\,\norm{e^{\beta t}u}_{E(\cly)}.
\end{equation}
Clearly $\clysp\beta E$ is a Banach space and $\clysp\beta E\subset E(\cly)_\loc$.

\medskip

\begin{rem}
\label{PropA2clysp}
As an easy consequence of Lemma \ref{PropA2clyset} we have that multiplication defines a continuous bilinear map $C^\infty(\cly)\times\clysp\beta E\to\clysp\beta E$ for any model space $E$ and $\beta\in\R$.
\end{rem}

\begin{rem}
\label{clyspforCinfty}
With $K$ as in Remark \ref{clyforCinfty} it is clear that we can define a locally convex space $\clysp\beta K\subset\mathscr D'(\cly)$ for any $\beta\in\R$. It is straightforward to check that Condition (A1) gives us continuous inclusions $\clysp\beta{\mathscr S}\hookrightarrow\clysp\beta E\hookrightarrow\clysp\beta{\mathscr S'}$ for any model space $E$ and $\beta\in\R$. 
\end{rem}

\smallskip

In order to work with the spaces $\clysp\beta E$ (and other spaces to be defined below) we need to introduce a set of auxiliary functions. Choose $\phi_0\in C^\infty(\R)$ with $\phi_0=0$ on $(-\infty,-2/3]$, $\phi_0=1$ on $[-1/3,\infty)$ and $\Ran\phi_0=[0,1]$. For any $i,j\in\Z\cup\{\pm\infty\}$ with $j\ge i$ define $\phi_{ij}\in C^\infty(\cly)$ by 
\[
\phi_{ij}(t,\omega)
\,=\,\phi_0(t\negvtsp -\negvtsp i)-\phi_0(t\negvtsp -\negvtsp j\negvtsp -\negvtsp 1).
\]
Therefore $\phi_{ij}$ is non-negative, $\supp(\phi_{ij})\subseteq(i\negvtsp-\negvtsp1,j\negvtsp+\negvtsp1)\!\times\!S^{n-1}$ and $\phi_{ij}=1$ on $[i,j]\!\times\!S^{n-1}$. We also set $\phi_i=\phi_{i\vtsp \infty}$ so $\phi_{ij}=\phi_{j+1}-\phi_i$.

\begin{rem}
\label{multibyphiijrem}
Suppose $\phi\in C^\infty(\cly)$ satisfies $\supp(\phi)\subseteq[i,j]\times S^{n-1}$ for some $i,j\in\Z\cup\{\pm\infty\}$ with $i\le j$. Then we have 
\[
\norm{\phi u}_{\clysp\beta E}
\;=\;\biggl\lVert\,\sum_{k=i}^j\phi_{kk}\phi u\,\biggr\rVert_{\clysp\beta E}
\;\le\;\sum_{k=i}^j\,\norm{\phi_{kk}\phi u}_{\clysp\beta E}
\]
for all $u\in E(\cly)_\loc$. On the other hand $\{\phi_{kk}\}_{k\in\Z}$ is a bounded subset of $C^\infty(\cly)$ so Remark \ref{PropA2clysp} implies
\[
\sup_{i\le k\le j}\norm{\phi_{kk}\phi u}_{\clysp\beta E}
\;\le\;C\norm{\phi u}_{\clysp\beta E}
\]
for all $u\in E(\cly)_\loc$.
\end{rem}

For any $k\in\Z$ we can write $e^{\beta t}=e^{\beta k}e^{\beta(t-k)}$ where $e^{\beta(t-k)}$ and its derivatives can be bounded independently of $k\in\Z$ on $\supp(\phi_{kk})\subseteq(k-1,k+1)$. Lemma~\ref{PropA2clyset} then gives
\begin{equation}
\label{rathimpteqntoreplem}
\norm{\phi_{kk} u}_{\clysp\beta E}
\,\asymp\,e^{\beta k}\norm{\phi_{kk} u}_{E(\cly)}\,,
\end{equation}
where the equivalence constants are independent of $k$.

\subsection{Weighted function spaces on $\pmb{\R}^{\boldsymbol n}_{\boldsymbol *}$ and $\pmb{\R}^{\boldsymbol n}$}

\label{wfnspondon}

Let $\Theta:\cly\to\don$ denote the diffeomorphism defined by $\Theta(t,\omega)=(r,\omega)$ where $r=e^t$. Under the pull back $\Theta^*$ we clearly have 
\begin{equation}
\label{simpleconjuncov}
r\to e^t,\quad
r\dfo_r\to\dfo_t
\quad\text{and}\quad
r^{-1}\d r\to dt. 
\end{equation}
Now $\Theta^*$ defines an isomorphism $\mathscr D'(\don)\to\mathscr D'(\cly)$. For any model space $E$ and $\beta\in\R$ we define $\donsp\beta E\subset\mathscr D'(\don)$ to be the preimage of $\clysp\beta E$ under $\Theta^*$ with the induced norm. Therefore the restricted map
\begin{equation}
\label{Thetaisaniso}
\Theta^*:\donsp\beta E\to\clysp\beta E
\quad\text{is an isomorphism}
\end{equation}
and 
\begin{equation}
\label{defofdonspnorm}
\norm{u}_{\donsp\beta E}
\,=\,\norm{\Theta^* u}_{\clysp\beta E}
\quad\text{for any $u\in\donsp\beta E$.}
\end{equation}

\bigskip

Choose $\theta\in\XYfns$ with $\theta=1$ on $\{\abs x\ge2\}$, $\theta=0$ on $\{\abs x\le1\}$ and $\Ran\theta=[0,1]$. For any model space $E$ and $\beta\in\R$ we define $\consp\beta E$ to be the set of $u\in\mathscr D'$ with $(1-\theta)u\in E$ and $\theta u\in\donsp\beta E$. On this set we define a norm
\[
\norm u_{\consp\beta E}
\,=\,\norm{(1\negvtsp-\negvtsp\theta)u}_E+\norm{\theta u}_{\donsp\beta E}.
\]
Straightforward calculations show that $\consp\beta E\subset E_\loc$, $\consp\beta E$ is a Banach space and the definition is independent of the choice of $\theta$ (up to equivalence of norms). These and other observations are summarised in the following result.

\begin{lem}
\label{equivofnorms1}
Suppose $\beta\in\R$, $\zeta\in\XYfns$ and $\eta\in C^\infty_0$ is a non-negative function for which $\eta+\zeta$ is bounded away from 0. Then $u\in\consp\beta E$ iff $\eta u\in E$ and $\zeta u\in\donsp\beta E$. Furthermore, all such $u$ satisfy an estimate of the from 
\[
\norm{u}_{\consp\beta E}
\,\asymp\,\norm{\eta u}_E+\norm{\zeta u}_{\donsp\beta E}.
\]
Also, for $u\in E(\don)_\loc$ we have $\zeta u\in\consp\beta E$ iff $\zeta u\in\donsp\beta E$ whilst all such $u$ satisfy an estimate of the form $\norm{\zeta u}_{\consp\beta E}\asymp\norm{\zeta u}_{\donsp\beta E}$. 
\end{lem}

\begin{rem}
\label{newequivofnormsrem1}
Suppose $f\in C^\infty$ is constant on a neighbourhood of 0 and satisfies $f\circ\Theta\in C^\infty(\cly)$ (we could take $f\in\genXYfns$ for example). For any $\beta\in\R$, Remark~\ref{PropA2clysp} and \eqref{Thetaisaniso} then imply that multiplication by $f$ defines a continuous map $\donsp\beta E\to\donsp\beta E$. Coupling this observation with Lemma \ref{equivofnorms1} and the fact that $f$ is constant on a neighbourhood of 0, it follows easily that multiplication by $f$ also defines a continuous map $\consp\beta E\to\consp\beta E$.
\end{rem}

\begin{rem}
\label{newequivofnormsrem2}
Suppose $\beta,\gamma\in\R$ and $\zeta\in\XYfns$. Making straightforward applications of Lemma \ref{equivofnorms1} (with $\eta=1-\zeta$) we can obtain the following. 
\begin{list}{}{\listinit}
\item[\emph{(i)}\hf]
If $\zeta u\in\consp\gamma E$ for some $u\in\consp\beta E$ then we also have $u\in\consp\gamma E$. 
\item[\emph{(ii)}\hf]
If $\zeta u\in\donsp\gamma E$ and $u=0$ in a neighbourhood of 0 for some $u\in\donsp\beta E$ then we also have $u\in\donsp\gamma E$. 
\end{list}
\end{rem}

\medskip

For each $i,j\in\Z\cup\{\pm\infty\}$ with $j\ge i$ define $\zeta_i,\zeta_{ij}\in C^\infty$ by $\zeta_i=\phi_i\circ\Theta^{-1}$ and $\zeta_{ij}=\phi_{ij}\circ\Theta^{-1}$ (n.b.~there are no problems with smoothness at 0 since $\phi_i$ and $\phi_{ij}$ are constant in a neighbourhood of $-\infty\times S^{n-1}$). Also set $\eta_i=\zeta_{-\infty\vtsp i-1}$. Therefore $\zeta_i\in\XYfns$, $\eta_i\in C^\infty_0$ and $\eta_i+\zeta_i=1$. The next result follows from Lemma \ref{equivofnorms1} with $\zeta=\zeta_0$, $\eta=\eta_0$ and the observation that $\zeta_0\zeta_i=\zeta_i$, $\eta_0\zeta_i=0$ for any $i\in\N$.

\begin{lem}
\label{equivofnorms2}
Given $i\in\N$ we have $\norm{\zeta_i u}_{\consp\beta E}\asymp\norm{\zeta_i u}_{\donsp\beta E}$ for any $u\in\consp\beta E$, where the equivalence constants are independent of $i$.
\end{lem}

\subsection{Basic properties}

\label{subsecBasicFunSpProp}

\subsubsection{Isomorphisms}

\begin{prop}
\label{basicchangbetaiso}
Suppose $E$ is a model space and $\beta,\gamma\in\R$. Then multiplication by $e^{\gamma t}$, $r^\gamma$ and $\Lambda^\gamma$ defines isomorphisms $\clysp{\beta+\gamma} E\to\clysp\beta E$, $\donsp{\beta+\gamma} E\to\donsp\beta E$ and $\consp{\beta+\gamma} E\to\consp\beta E$ respectively.
\end{prop}

\begin{proof}
The first isomorphism is an immediate consequence of \eqref{defofnormonclysp} whilst the second then follows from \eqref{simpleconjuncov} and \eqref{Thetaisaniso}. Now define a smooth function $f$ on $\R^n$ by $f=\zeta_0\Lambda^\gamma r^{-\gamma}$. It follows that $f\circ\Theta$ is a smooth function on $\cly$ with $f(\Theta(t,\omega))=0$ for $t\le-1$ and $f(\Theta(t,\omega))=(1+e^{-2t})^{\gamma/2}$ for $t\ge 1$. Hence $f\circ\Theta\in C^\infty(\cly)$, so multiplication by $f$ defines a continuous map $\donsp\beta E\to\donsp\beta E$ by Remark \ref{newequivofnormsrem1}. Combining this observation with Lemma \ref{equivofnorms1}, the identity $\zeta_1\Lambda^\gamma=f\zeta_1 r^\gamma$ and the second part of the present result, we now have
\[
\norm{\Lambda^\gamma u}_{\consp\beta E}
\ \le\ C\bigl(\norm{\eta_1\Lambda^\gamma u}_E+\norm{f\zeta_1r^\gamma u}_{\donsp\beta E}\bigr)
\ \le\ C\bigl(\norm{\eta_1\Lambda^\gamma u}_E+\norm{\zeta_1 u}_{\donsp{\beta+\gamma}E}\bigr)
\ \le\ C\norm{u}_{\consp{\beta+\gamma}E}
\]
for all $u\in\consp{\beta+\gamma}E$ (n.b.~$\eta_1\Lambda^\gamma\in C^\infty_0$ is non-negative whilst $\eta_1\Lambda^\gamma+\zeta_1$ is bounded away from 0). A similar argument for $\Lambda^{-\gamma}$ now completes the result. 
\end{proof}

\subsubsection{Inclusions}

\begin{rem}
\label{easyincrem}
It is easy to see that any continuous inclusion $E\hookrightarrow F$ between model spaces $E$ and $F$ induces continuous inclusions $\clysp\beta E\hookrightarrow\clysp\beta F$, $\donsp\beta E\hookrightarrow\donsp\beta F$ and $\consp\beta E\hookrightarrow\consp\beta F$ for any $\beta\in\R$. 
\end{rem}

By a \emph{local inclusion $E_\loc\hookrightarrow F_\loc$} we mean that $E_\loc\subseteq F_\loc$ and, given $\phi\in C^\infty_0$, $\norm{\phi u}_F\le C\norm{\phi u}_E$ for all $u\in E_\loc$ (where $C$ may depend on $\phi$). 

\begin{rem}
\label{locandbndincegs}
We have a local inclusion $\Sob{p_1}{s_1}_\loc\hookrightarrow\Sob{p_2}{s_2}_\loc$ whenever $p_1,p_2\in[1,\infty)$ and $s_1,s_2\in\R$ satisfy $s_1\ge s_2$ and $s_1-n/p_1\ge s_2-n/p_2$. The additional condition $p_1\le p_2$ is needed in order to get a continuous inclusion $\Sob{p_1}{s_1}\hookrightarrow\Sob{p_2}{s_2}$. We also have continuous inclusions $\Sob{p_1}{s_1}\hookrightarrow\Zyg{s_2}$ whenever $s_1-n/p_1\ge s_2>0$ and $\Sob{p_1}{s_1}\hookrightarrow\Cnt k$ whenever $k\in\NZ$ satisfies $s_1-n/p_1>k$. Further details can be found in Sections 2.3.2, 2.7.1 and 3.3.1 of \cite{Tr1}.
\end{rem}

Obviously a continuous inclusion $E\hookrightarrow F$ leads to a local inclusion $E_\loc\hookrightarrow F_\loc$. Although the converse does not hold in general we can obtain the following related results, the first of which is an easy consequence of the compactness of $S^{n-1}$.

\begin{lem}
\label{incforfunsponSnm1}
A local inclusion $E_\loc\hookrightarrow F_\loc$ for model spaces $E$ and $F$ on $\R^{n-1}$ induces a continuous inclusion $E(S^{n-1})\hookrightarrow F(S^{n-1})$.
\end{lem}

\begin{lem}
\label{oftios}
Suppose we have a local inclusion $E_\loc\hookrightarrow F_\loc$ for some model spaces $E$ and $F$. If $\beta\in\R$, $\epsilon>0$, $l\in\Z$ and $\phi_\pm\in C^\infty(\cly)$ with $\supp(\phi_\pm)\subseteq\pm[l,+\infty)\times S^{n-1}$, then 
\[
\norm{\phi_\pm u}_{\clysp\beta F}
\;\le\;C(\epsilon,\negvtsp l)\vtsp\norm{\phi_\pm u}_{\clysp{\beta\pm\epsilon} E}
\]
for any $u\in E(\cly)_\loc$. Furthermore $C(\epsilon,\negvtsp l)$ can be chosen independently of $\phi_\pm$.
\end{lem}

\begin{proof}
Clearly Condition (A3) implies the norm $\norm\cdot_{F(\cly)}$ is invariant under translations with respect to the first variable of $\cly$. With the help of \eqref{rathimpteqntoreplem} it follows that the local inclusion $E_\loc\hookrightarrow F_\loc$ leads to an estimate
\begin{equation}
\label{locincforclysps}
\norm{\phi_{kk}v}_{\clysp\beta F}
\;\le\;Ce^{-\epsilon k}\norm{\phi_{kk}v}_{\clysp{\beta+\epsilon}E}
\end{equation}
for any $k\in\Z$ and $v\in E(\cly)_\loc$, where $C$ is independent of $k$. On the other hand $\sum_{k\ge l}\phi_{kk}=1$ on $\supp(\phi_+)$. Combined with Remark \ref{multibyphiijrem} and \eqref{locincforclysps} we now get 
\begin{multline*}
\norm{\phi_+ u}_{\clysp\beta F}
\;\le\;\sum_{k\ge l}\norm{\phi_{kk}\phi_+ u}_{\clysp\beta F}
\;\le\;C\sum_{k\ge l}e^{-\epsilon k}\norm{\phi_{kk}\phi_+ u}_{\clysp{\beta+\epsilon}E}\\
\le\;C\sum_{k\ge l}e^{-\epsilon k}\;
\Bigl(\vtsp\sup_{k\ge l}\norm{\phi_{kk}\phi_+ u}_{\clysp{\beta+\epsilon}E}\Bigr)
\;\le\;C\norm{\phi_+ u}_{\clysp{\beta+\epsilon}E}
\end{multline*}
for any $u\in E(\cly)_\loc$. Clearly a similar argument can be used for $\phi_-$. 
\end{proof}

\begin{rem}
\label{previouslemforZS}
If $\phi_\pm$ are as in the previous lemma it is easy to check that multiplication by $\phi_\pm$ defines  continuous maps $\clysp\beta{\mathscr S}\to\clysp{\beta\mp\epsilon}{\mathscr S}$ and $\clysp\beta{\mathscr S'}\to\clysp{\beta\mp\epsilon}{\mathscr S'}$ for any $\beta\in\R$ and $\epsilon\ge0$. 
\end{rem}

\begin{prop}
\label{genformincforconsp}
Suppose we have a local inclusion $E_\loc\hookrightarrow F_\loc$ for some model spaces $E$ and $F$. If $\beta,\gamma\in\R$ with $\beta<\gamma$ then we have a continuous inclusion $\consp\gamma E\hookrightarrow\consp\beta F$.
\end{prop}

This result obviously implies that we have a continuous inclusion 
\begin{equation}
\label{cntsincforXE}
\consp\gamma E\hookrightarrow\consp\beta E
\quad\text{for $\beta\le\gamma$.}
\end{equation}

\begin{proof}
Let $u\in\consp\gamma E$. Therefore $u\in E_\loc\subseteq F_\loc$. Now Lemma \ref{equivofnorms1} gives
\begin{equation}
\label{logapsinthepfatmom}
\norm{u}_{\consp\beta F}
\;\le\;C\bigl(\norm{\eta_1 u}_F+\norm{\zeta_1 u}_{\donsp\beta F}\bigr).
\end{equation}
However \eqref{defofdonspnorm} and Lemma \ref{oftios} (with $\phi_+=\phi_1$) combine to give
\[
\norm{\zeta_1u}_{\donsp\beta F}
\;=\;\norm{\phi_1\vtsp\Theta^* u}_{\clysp\beta F}
\;\le\;C\norm{\phi_1\vtsp\Theta^* u}_{\clysp\gamma E}
\;=\norm{\zeta_1u}_{\donsp\gamma E}.
\]
Together with the definition of local inclusion and \eqref{logapsinthepfatmom} we now have 
\[
\norm{u}_{\consp\beta F}
\;\le\;C\bigl(\norm{\eta_1 u}_E+\norm{\zeta_1 u}_{\donsp\gamma E}\bigr)
\;\le\;C\norm{u}_{\consp\gamma E},
\]
the last inequality following from a further application of Lemma \ref{equivofnorms1}.
\end{proof}

\bigskip

We finish this section with some results which are direct consequences of Conditions (A1) to (A3) for model spaces. We will make use of the collection of semi-norms $\{p_l\,|\,l\in\NZ\}$ for $\mathscr S$ introduced in Section \ref{notation}; in particular we observe that $p_{l'}(u)\le p_l(u)$ whenever $0\le l'\le l$.

\begin{lem}
\label{CntlintoEloc}
We have $\Cnt l_\loc\hookrightarrow E_\loc$ for all sufficiently large $l\in\NZ$. 
\end{lem}

\begin{proof}
The continuous inclusion $\mathscr S\hookrightarrow E$ given by Condition (A1) simply means that we can find $l'\in\NZ$ such that $\norm{u}_E\le C p_{l'}(u)$ for all $u\in\mathscr S$. Now let $l\in\NZ$ with $l\ge l'$. Also let $\phi\in C^\infty_0$ and choose $\phi_1\in C^\infty_0$ with $\phi_1\succ\phi$. Therefore $\dfo_x^\alpha(\phi u)=\phi_1\dfo_x^\alpha(\phi u)$ for any multi-index $\alpha$ whilst $\Lambda^l\phi_1\in C^\infty_0$. Hence
\[
\norm{\phi u}_E
\;\le\;Cp_l(\phi u)
\;=\;C\sum_{\abs\alpha\le l}\,\norm{\Lambda^l\phi_1\dfo_x^\alpha(\phi u)}_{L^\infty}
\;\le\;C\sum_{\abs\alpha\le l}\,\norm{\dfo_x^\alpha(\phi u)}_{L^\infty}
\;=\;C\norm{\phi u}_{\Cnt l}
\]
for all $u\in\mathscr S$. The result now follows from the fact that if $u\in\Cnt l$ then we can approximate $\phi u$ arbitrarily closely (in the $\Cnt l$ norm) by $\phi u'$ for some $u'\in\mathscr S$.
\end{proof}

\begin{lem}
\label{techlemmfromfana}
Suppose $L$ is a locally convex space whose topology is given by a countable collection of semi-norms $\{q_l\,|\,l\in\N\}$. Also suppose $X$ and $Y$ are Banach spaces and $T:L\times X\to Y$ is a bilinear mapping which is continuous in each variable. Then we can find $l\in\N$ and a constant $C$ such that 
\begin{equation}
\label{bilinearest}
\norm{T(a,x)}_Y\;\le\;C\norm{x}_X\sum_{i=1}^lq_i(a)
\end{equation}
for all $(a,x)\in L\times X$.
\end{lem}

\begin{proof}
Without loss of generality we may assume that the semi-norms are defined so that $q_i\le q_j$ for all $i,j\in\N$ with $i\le j$. Therefore \eqref{bilinearest} can be rewritten as 
\[
\norm{T(a,x)}_Y\;\le\;C\norm{x}_X\, q_l(a)
\]
for all $(a,x)\in L\times X$. Suppose an estimate of this form is not valid. It follows that we can find a sequence $\{a_i\}_{i\in\N}$ in $L$ such that $q_i(a_i)\le 1$ and $\norm{T(a_i,\cdot)}_{\mathscr L(X,Y)}\to\infty$ as $i\to\infty$. Now the continuity of the map $T(\cdot,x):L\to Y$ gives us $j\in\N$ and a constant $C$ such that $\norm{T(a_i,x)}_Y\le C q_j(a_i)$ for all $i\in\N$. It follows that the set $\{T(a_i,x)\vtsp|\vtsp i\in\N\}$ is bounded in $Y$ (by the maximum of $Cq_j(a_1),\dots,Cq_j(a_{j-1})$ and $C$). The Uniform Boundedness Theorem (see Section II.1 in \cite{Y} for example) then implies $\{\norm{T(a_i,\cdot)}_{\mathscr L(X,Y)}\vtsp|\vtsp i\in\N\}$ must also be bounded. The result now follows by contradiction.
\end{proof}

\begin{lem}
\label{ElintoCntglob}
We have $E^l\hookrightarrow\Cnt0$ for all sufficiently large $l\in\NZ$.
\end{lem}

\begin{proof}
Using the inclusion $E\hookrightarrow\mathscr S'$ and the dual pairing $\mathscr S\times\mathscr S'\to\C$ we can define a continuous bilinear map $\mathscr S\times E\to\C$ by $(\phi,u)\mapsto(\phi,u)_{\R^n}$. By Lemma \ref{techlemmfromfana} we can thus find $j\in\N$ and a constant $C$ such that 
\[
\abs{(\phi,u)_{\R^n}}
\;\le\;C\norm{u}_E\,p_j(\phi)
\]
for all $\phi\in\mathscr S$ and $u\in E$. Now suppose $\chi\in C^\infty_0$ and choose $\chi_1,\chi_2\in C^\infty_0$ with $\chi_2\succ\chi_1\succ\chi$. Thus $\dfo_x^\alpha(\chi_1\phi)=\chi_2\dfo_x^\alpha(\chi_1\phi)$ and $\Lambda^j\chi_2\in C^\infty_0$ so
\[
p_j(\chi_1\phi)
\;=\;\sum_{\abs\alpha\le j}\,\norm{\Lambda^j\chi_2\dfo_x^\alpha(\chi_1\phi)}_{L^\infty}
\;\le\;C\norm{\chi_1\phi}_{\Cnt j}
\;\le\;C\norm{\phi}_{\Sob 2{j+n}}
\]
for all $\phi\in\mathscr S$, where we have used the continuous inclusion $\Sob 2{j+n}\hookrightarrow\Cnt j$ in the last inequality. Therefore 
\[
\abs{(\phi,\chi u)_{\R^n}}
\;=\;\abs{(\chi_1\phi,\chi u)_{\R^n}}
\;\le\;C\norm{\chi u}_E\,\norm{\phi}_{\Sob 2{j+n}}
\]
for all $u\in E$ and $\phi\in\mathscr S$. Now let $l\in\NZ$ with $l>2n+j$. Choose $l'\in2\NZ$ so that $2n+j\le l'\le l$. Now the Fourier multiplier $\Lambda^{-l'}(\dfo)$ defines an isomorphism on $\mathscr S$ whose inverse is simply the constant coefficient differential operator $\Lambda^{l'}(\dfo)$ of order $l'$ (recall that $l'$ is even). Therefore 
\begin{eqnarray}
\abs{(\phi,\chi u)_{\R^n}}
&=&\bigabs{\bigl(\Lambda^{l'}(\dfo)\Lambda^{-l'}(\dfo)\phi,\vtsp\chi u\bigr)_{\R^n}}
\nonumber\\
&=&\bigabs{\bigl(\Lambda^{-l'}(\dfo)\phi,\vtsp\Lambda^{l'}(\dfo)(\chi u)\bigr)_{\R^n}}
\nonumber\\
&\le&C\norm{\Lambda^{l'}(\dfo)(\chi u)}_E\vtsp\norm{\Lambda^{-l'}(\dfo)\phi}_{\Sob 2{j+n}}
\nonumber\\
\label{dualspres}
&\le&C\norm{\chi u}_{E^l}\norm{\phi}_{\Sob 2{-n}}
\end{eqnarray}
for all $u\in E^l$ and $\phi\in\mathscr S$, where the last inequality follows from Remark \ref{remnofElver1}, the fact that $\Lambda^{-l'}(\dfo)$ defines an isomorphism $\Sob 2{j+n-l'}\to\Sob 2{j+n}$ and the inequality $j+n-l'\le-n$. Since $\mathscr S$ is dense in $\Sob 2{-n}$, \eqref{dualspres} implies that for any $u\in E^l$ we have $\chi u\in\Sob 2n$ (the dual space of $\Sob 2{-n}$) with a corresponding norm estimate. Since we have a continuous inclusion $\Sob 2n\hookrightarrow\Cnt 0$ whilst $\chi\in C^\infty_0$ was arbitrary, we finally arrive at a local inclusion $E^l_\loc\hookrightarrow\Cnt0_\loc$.

\smallskip

Let $\{\chi_I\}_{I\in\Z^n}$ be a partition of unity of $\R^n$ where $\chi_0\in C^\infty_0$ and, for each $I\in\Z^n$, $\chi_I(x)=\chi_0(x-I)$. Using Conditions (A2) and (A3) we thus have $\norm{\chi_I u}_{E^l}\le C\norm{u}_{E^l}$ for all $u\in E^l$, where $C$ is independent of $I\in\Z^n$. On the other hand it is clear that $\norm{u}_{\Cnt 0}\le C\sup_{I\in\Z^n}\norm{\chi_I u}_{\Cnt 0}$ (where for $C$ we can take $\#\{I\in\Z^n\vtsp|\vtsp \supp(\chi_0)\cap\supp(\chi_I)\ne\emptyset\}$). The fact that we have a continuous inclusion $E^l\hookrightarrow\Cnt 0$ now follows from the existence of a local inclusion.
\end{proof}

Clearly a continuous inclusion $E\hookrightarrow F$ for model spaces $E$ and $F$ leads to a continuous inclusion $E^l\hookrightarrow F^l$ for any $l\in\NZ$. Coupling this observation with Lemmas \ref{CntlintoEloc} and \ref{ElintoCntglob} we immediately get the following.

\begin{cor}
\label{EFlocincpositivever}
For any model spaces $E$ and $F$ we have a local inclusion $E^l_\loc\hookrightarrow F_\loc$ for all sufficiently large $l\in\NZ$. 
\end{cor}

\subsubsection{Derivatives}

For each $i\in\{1,\dots,n\}$ we can write 
\begin{equation}
\label{defofpiPi}
\dfo_i\,=\,r^{-1}\bigl(b_i(\omega)(r\dfo_r)+P_i(\omega,\dfo_\omega)\bigr)
\end{equation}
where $b_i(\omega)=x_i/r\in C^\infty(S^{n-1})$ and $P_i(\omega,\dfo_\omega)$ is a first order differential operator on $S^{n-1}$. Now any first order differential operator $A$ on $S^{n-1}$ we can be written as a (not necessarily unique) linear combination 
\[
A\,=\,a_0+\sum_{i=1}^n\, a_i P_i
\]
where $a_0,\dots,a_n\in C^\infty(S^{n-1})$. On the other hand $b_1^2+\dots+b_n^2=1$ while $b_1P_1+\dots+b_nP_n=0$. Defining differential operators $B_i$ on $\cly$ by $B_i=b_i(\omega)\dfo_t+P_i(\omega,\dfo_\omega)$ for $i=1,\dots,n$ we thus arrive at the following result.

\begin{lem}
\label{locrepofscamodop1}
If $A(\omega,\dfo_\omega,\dfo_t)$ is a first order differential operator on $\cly$ whose coefficients do not depend upon $t$ then we can write
\[
A\,=\,a_0+\sum_{i=1}^n\, a_i B_i
\]
for some $a_0,\dots,a_n\in C^\infty(S^{n-1})$.
\end{lem}

\begin{lem}
\label{diffoponclysplem}
Suppose $E$ be a model space and $\beta\in\R$. Then 
\[
\norm{u}_{\clysp\beta{E^1}}
\;\asymp\;\norm{u}_{\clysp\beta E}
+\sum_{i=1}^n\,\norm{B_i u}_{\clysp\beta E}
\]
for all $u\in\mathscr D'(\cly)$ (where we define the norm of a function not belonging to the relevant space to be $+\infty$). In particular $B_i$ defines a continuous map $\clysp\beta{E^1}\to\clysp\beta E$ for $i=1,\dots,n$.
\end{lem}

\begin{proof}
Consider the notation of Remark \ref{enlargechirem} and let $i\in I$ and $j\in\{1,\dots,n\}$. Now $(\Psi_i^{-1})^*\chi_iB_j\chi_i'\Psi_i^*$ is a first order differential operator on $\R^n$ whose coefficients are contained in $C^\infty$ (in fact the coefficients do not depend upon the first variable and are compactly supported with respect to the rest). Thus we have 
\[
\norm{(\Psi_i^{-1})^*\chi_i(B_j+i\beta b_j)\chi_i'\Psi_i^*v}_E
\;\le\;C\norm{v}_{E^1}
\]
for any $v\in\mathscr D'$. On the other hand $e^{\beta t}B_ju=(B_j+i\beta b_j)e^{\beta t}u$ so
\[
(\Psi_i^{-1})^*\chi_i\vtsp e^{\beta t}B_ju
\;=\;(\Psi_i^{-1})^*\chi_i(B_j+i\beta b_j)\chi_i'\Psi_i^*
\vtsp(\Psi_i^{-1})^*\chi_i'\vtsp e^{\beta t}u
\]
for any $i\in I$, $j\in\{1,\dots,n\}$ and $u\in\mathscr D'(\cly)$. Since we clearly have $\norm{u}_{\consp\beta E}\le C\norm{u}_{\consp\beta{E^1}}$, the above results combine to give 
\[
\norm{u}_{\clysp\beta E}
+\sum_{i=1}^n\,\norm{B_i u}_{\clysp\beta E}
\;\le\;C\norm{u}_{\clysp\beta{E^1}}
\]
for all $u\in\mathscr D'(\cly)$.

\medskip

Let $i\in I$. Using Definition \ref{defnofElver1} we get 
\[
\norm{(\Psi_i^{-1})^*\chi_i e^{\beta t}u}_{E^1}
\;\le\;\norm{(\Psi_i^{-1})^*\chi_i e^{\beta t}u}_E
+\sum_{j=1}^n\,\norm{\dfo_j(\Psi_i^{-1})^*\chi_i e^{\beta t}u}_E
\]
for all $u\in\mathscr D'(\cly)$. On the other hand we can write
\[
\dfo_j(\Psi_i^{-1})^*\chi_i e^{\beta t}u
\;=\;(\Psi_i^{-1})^*\chi'_i e^{\beta t}A_{ij} u
\]
where $A_{ij}=\chi_i'\Psi^*(D_j-i\beta\delta_{j1})(\Psi^{-1})^*\chi_i$ is a first order differential operator on $\cly$ whose coefficients are independent of $t$. By Lemma \ref{cordequivnorm}, Remark \ref{enlargechirem} and \eqref{defofnormonclysp} we thus have
\[
\norm{u}_{\clysp\beta{E^1}}
\;\le\;C\Bigl(\norm{u}_{\clysp\beta E}
+\sum_{i\in I}\sum_{j=1}^n\,\norm{A_{ij} u}_{\clysp\beta E}\Bigr)
\]
for all $u\in\mathscr D'(\cly)$. On the other hand, using Lemma \ref{locrepofscamodop1} and the fact that multiplication by an element of $C^\infty(S^{n-1})$ defines a continuous map on $\clysp\beta E$ (see Remark \ref{PropA2clysp}), we get 
\[
\norm{A_{ij} u}_{\clysp\beta E}
\;\le\;C\Bigl(\norm{u}_{\clysp\beta E}+\sum_{j'=1}^n\,\norm{B_{j'}u}_{\clysp\beta E}\Bigr)
\]
for all $u\in\mathscr D'(\cly)$. Clearly the last two estimates complete the result.
\end{proof}

For any $i\in\{1,\dots,n\}$, \eqref{simpleconjuncov} and \eqref{defofpiPi} give us $\dfo_i=r^{-1}(\Theta^{-1})^*B_i\Theta^*$. With the help of \eqref{Thetaisaniso} and Proposition \ref{basicchangbetaiso}, Lemma \ref{diffoponclysplem} now implies that 
\begin{equation}
\label{estfordonspE1}
\norm{u}_{\donsp\beta{E^1}}
\;\asymp\;\norm{u}_{\donsp\beta E}
+\sum_{i=1}^n\,\norm{\dfo_i u}_{\donsp{\beta+1}E}
\end{equation}
for all $u\in\mathscr D'(\don)$.

\begin{prop}
\label{lemoncontofpderiv}
Let $E$ be a model space, $l\in\NZ$ and $\beta\in\R$. Given $u\in\mathscr D'$ we have $u\in\consp\beta{E^l}$ iff $\dfo_x^\alpha u\in\consp{\beta+\abs\alpha}E$ for all multi-indices $\alpha$ with $\abs\alpha\le l$. Furthermore
\[
\norm{u}_{\consp\beta{E^l}}
\;\asymp\;\sum_{\abs\alpha\le l}\,\norm{\dfo_x^\alpha u}_{\consp{\beta+\abs\alpha}E}
\]
for all $u\in\consp\beta{E^l}$. In particular the differential operator $D_x^\alpha$ defines a continuous map $\consp\beta{E^l}\to\consp{\beta+\abs\alpha}{E}$ whenever $\abs\alpha\le l$.
\end{prop}

\begin{proof}
Induction clearly reduces the proof to the case $l=1$. Now let $i\in\{1,\dots,n\}$ and $u\in\mathscr D'$.
Since $\dfo_i\eta_1=-\dfo_i\zeta_1\in C^\infty_0(\don)$ and multiplication by such functions defines continuous maps on $E$ and $\donsp\beta E$ (see Remark \ref{newequivofnormsrem1}), Lemma \ref{equivofnorms1} gives us
\begin{eqnarray*}
\lefteqn{\norm{(\dfo_i\eta_1)u}_E+\norm{(\dfo_i\zeta_1)u}_{\donsp\beta E}}\qquad\\
&&=\ \norm{(\dfo_i\eta_1)\eta_2u}_E+\norm{(\dfo_i\zeta_1)\zeta_0u}_{\donsp\beta E}
\ \le\ C\bigl(\norm{\eta_2u}_E+\norm{\zeta_0u}_{\donsp\beta E}\bigr)
\ \le\ C\norm{u}_{\consp\beta E}\,.
\end{eqnarray*}
Combining this with Definition \ref{defnofElver1}, Lemma \ref{equivofnorms1} and \eqref{estfordonspE1} we thus have
\begin{eqnarray*}
\norm{\dfo_iu}_{\consp{\beta+1}E}
&\le&C\bigl(\norm{\eta_1\dfo_iu}_E+\norm{\zeta_1\dfo_iu}_{\donsp{\beta+1}E}\bigr)\\
&\le&C\bigl(\norm{\dfo_i(\eta_1u)}_E+\norm{\dfo_i(\zeta_1u)}_{\donsp{\beta+1}E}
+\norm{u}_{\consp\beta E}\bigr)\\
&\le&C\bigl(\norm{\eta_1u}_{E^1}+\norm{\zeta_1u}_{\donsp\beta{E^1}}
+\norm{u}_{\consp\beta E}\bigr)\\
&\le&C\norm{u}_{\consp\beta{E^1}}
\end{eqnarray*}
and 
\begin{eqnarray*}
\norm{u}_{\consp\beta{E^1}}
&\le&C\bigl(\norm{\eta_1u}_{E^1}+\norm{\zeta_1u}_{\donsp\beta{E^1}}\bigr)\\
&\le&C\Bigl(\norm{\eta_1u}_E+\norm{\zeta_1u}_{\donsp\beta E}
+\sum_{i=1}^n\bigl(\norm{\dfo_i(\eta_1u)}_E
+\norm{\dfo_i(\zeta_1u)}_{\donsp{\beta+1}E}\bigr)\Bigr)\\
&\le&C\Bigl(\norm{u}_{\consp\beta E}
+\sum_{i=1}^n\bigl(\norm{\eta_1\dfo_iu}_E
+\norm{\zeta_1\dfo_iu)}_{\donsp{\beta+1}E}\bigr)\Bigr)\\
&\le&C\Bigl(\norm{u}_{\consp\beta E}
+\sum_{i=1}^n\norm{\dfo_iu}_{\consp{\beta+1}E}\Bigr)
\end{eqnarray*}
for all $u\in\mathscr D'$. The result follows.
\end{proof}

\subsubsection{Separable Subspaces}

Suppose $E$ is a model space and let $E_0$ denote the separable subspace obtained by taking the closure of $C^\infty_0$ (see Definition \ref{defofsss} and Remark \ref{defofsssrem}). 

\begin{lem}
\label{denofCinft0incertsp}
For any $\beta\in\R$, $\clysp\beta{E_0}$ is the closure of $C^\infty_0(\cly)$ in $\clysp\beta E$, $\donsp\beta{E_0}$ is the closure of $C^\infty_0(\don)$ in $\donsp\beta E$ and $\consp\beta{E_0}$ is the closure of $C^\infty_0$ in $\consp\beta E$.
\end{lem}

\begin{proof}
For any $\beta\in\R$ it is easy to check that $C^\infty_0(\cly)\subset\clysp\beta E$ whilst the isometric inclusion $E_0\hookrightarrow E$ leads to an isometric inclusion $\clysp\beta{E_0}\hookrightarrow\clysp\beta E$. Since $\clysp\beta{E_0}$ is complete it thus remains to show that $C^\infty_0(\cly)$ is dense in $\clysp\beta{E_0}$. In turn, using the isomorphism $e^{\beta t}:\clysp\beta{E_0}\to E_0(\cly)$ (see \eqref{defofnormonclysp}) and the definition of the norm on $E_0(\cly)$, it is clear that the first part of the result is completed by the following.

\begin{claim}[Claim: Suppose $(\Psi,\R\times U)$ and $\chi$ are as in Lemma \ref{cordequivnorm}. If $u\in E_0(\cly)$ and $\epsilon>0$ then there exists $u_\epsilon\in E_0(\cly)$ with $\norm{\chi u-u_\epsilon}_{E_0(\cly)}<\epsilon$.]
Choose $\chi_1\in C^\infty_0(U)$ with $\chi_1\succ\chi$ and set $\widetilde\chi_1=\chi_1\circ\Psi^{-1}\in C^\infty$. Now Lemma \ref{cordequivnorm} gives us $(\Psi^{-1})^*(\chi u)\in E_0$ so, given any $\delta>0$, we can find $v_\delta\in C^\infty_0$ with $\norm{(\Psi^{-1})^*(\chi u)-v_\delta}_E<\delta$. On the other hand $\widetilde\chi_1\in C^\infty$ so Condition (A2) implies 
\[
\norm{(\Psi^{-1})^*(\chi u)-\widetilde\chi_1 v_\delta}_E
\;=\;\bignorm{\widetilde\chi_1\bigl((\Psi^{-1})^*(\chi u)-v_\delta\bigr)}_E
\;<\;C\delta
\]
where $C$ is independent of $\delta$. Setting $u_\delta=\Psi^*(\widetilde\chi_1 v_\delta)\in C^\infty_0(\cly)$ we thus get $\norm{\chi u-u_\delta}_{E(\cly)}< C\delta$ by Lemma \ref{cordequivnorm}, where $C$ is again independent of $\delta$. This completes the claim.
\end{claim}

\medskip

The second part of the result follows from \eqref{Thetaisaniso} and the fact that $\Theta^*(C^\infty_0(\don))=C^\infty_0(\cly)$.

\medskip

Arguing as above it is clear that $\consp\beta{E_0}$ is a closed subspace of $\consp\beta E$ which contains $C^\infty_0$. Now let $u\in\consp\beta{E_0}$. By Lemma \ref{equivofnorms1} $\eta_1u\in E_0$ and $\zeta_1u\in\donsp\beta{E_0}$ so, given $\delta>0$, we can find $v\in C^\infty_0$ and $w\in C^\infty_0(\don)$ with $\norm{\eta_1 u-v}_E,\norm{\zeta_1 u-w}_{\donsp\beta E}<\delta$. Thus $\eta_2 v+\zeta_0 w\in C^\infty_0$. On the other hand Condition (A2) gives
\begin{equation}
\label{newlabonrewr1}
\norm{\eta_1 u-\eta_2 v}_E
\;=\;\norm{\eta_2(\eta_1 u-v)}_E
\;<\;C\delta,
\end{equation}
whilst from Remark \ref{newequivofnormsrem1} we have 
\begin{equation}
\label{newlabonrewr2}
\norm{\zeta_1 u-\zeta_0w}_{\donsp\beta E}
\;=\;\norm{\zeta_0(\zeta_1 u-w)}_{\donsp\beta E}
\;<\;C\delta.
\end{equation}
Combining \eqref{newlabonrewr1}, \eqref{newlabonrewr2} and Lemma \ref{equivofnorms1} we get 
\begin{multline*}
\norm{u-(\eta_2 v+\zeta_0 w)}_{\consp\beta E}
\;\le\;\norm{\eta_1 u-\eta_2 v}_{\consp\beta E}+\norm{\zeta_1 u-\zeta_0 w}_{\consp\beta E}\\
\le\;C\bigl(\norm{\eta_1 u-\eta_2 v}_E+\norm{\zeta_1 u-\zeta_0 w}_{\donsp\beta E}\bigr)
\;<\;C\delta.
\end{multline*}
The fact that $C$ is independent of $\delta$ completes the result.
\end{proof}

\subsubsection{Equivalent Norms for Some Model Spaces}

\label{equivnormssssec}

In this section we consider some equivalent norms on the spaces $\consp\beta E$ when $E$ is either a Sobolev space of positive integral order or a H\"older space. In particular this will allow us to identify the spaces $H^{p,k}_\beta$ given in Definition \ref{defnofHpkbeta1} in the Introduction with $\consp\beta E$ for appropriate $\beta$ and $E$.

\begin{prop}
\label{equivnormmodspSob}
Suppose $\beta\in\R$, $p\in[1,\infty)$ and $k\in\NZ$. Given $u\in\mathscr D'$ we have $u\in\consp\beta{\Sob pk}$ iff $\dfo_x^\alpha u$ is a measurable function for all multi-indices $\alpha$ with $\abs\alpha\le k$ and 
\begin{equation}
\label{sumgivealtnorm1}
\sum_{\abs\alpha\le k}\int\Lambda^{p(\beta+\abs\alpha)-n}(x)
\,\abs{\dfo_x^\alpha u(x)}^p\d^nx
\;<\;+\infty.
\end{equation}
Furthermore the quantity on the left hand side of \eqref{sumgivealtnorm1} is equivalent to $\norm{u}_{\consp\beta{\Sob pk}}^p$.
\end{prop}

\begin{rem}
\label{actidentifofequivnorms}
Let $\beta\in\R$, $p\in[1,\infty)$ and $k\in\NZ$. Since the Sobolev space $\Sob pk$ contains $C^\infty_0$ as a dense subset, Lemma \ref{denofCinft0incertsp} and Proposition \ref{equivnormmodspSob} immediately imply that the space $H^{p,k}_\beta$ given in Definition \ref{defnofHpkbeta1} in the Introduction is simply $\consp{\beta+n/p}{\Sob pk}$ (up to equivalent norms).
\end{rem}

\begin{proof}[Proof of Proposition \ref{equivnormmodspSob}]
Propositions \ref{basicchangbetaiso} and \ref{lemoncontofpderiv} reduce our task to proving the result under the assumptions that $\beta=0$ and $k=0$.

Let $I$, $\psi_i$, $\Psi_i$, $U_i$ and $\chi_i$ be as given in the introduction to Section \ref{subsecweightfunspcly}. Now the pull back of the density $d^nx$ on $\R\times\psi_i(U_i)\subseteq\R^n$ is a density on $\R\times U_i\subset\cly$ so we can write $(\Psi_i)^*(d^nx)=J_i\d t\d S^{n-1}$ for some positive function $J_i$ defined on $\R\times U_i$. In particular
\[
\int(\Psi_i^{-1})^*v\d^nx
\;=\;\int v\vtsp J_i\d t\d S^{n-1}
\]
for any measurable function $v$ which is supported on $\R\times U_i$. Now $J_i$ is independent of $t$ whilst $\supp(\chi_i)\subset U_i$ is compact. Therefore $J_i$ is bounded and bounded away from 0 on $\R\times\supp(\chi_i)$. Hence we have 
\[
\bignorm{(\Psi_i^{-1})^*(\chi_i u)}_{L^p}^p
\ =\ \int\bigabs{(\Psi_i^{-1})^*(\chi_i u)}^p\d^nx
\ =\ \int\bigabs{\chi_iu}^p\vtsp J_i\d t\d S^{n-1}
\ \asymp\ \int\bigabs{\chi_iu}^p \d t\d S^{n-1}
\]
for any measurable function $u$ on $\cly$. From the definition of $\clysp0{L^p}$ (see \eqref{defofEclynobeta} and \eqref{defofnormonclysp}) and the finiteness of $I$ it follows that
\[
\norm{u}_{\clysp0{L^p}}^p
\;\asymp\;\int\Bigl(\sum_{i\in I}\,\abs{\chi_i u}\Bigr)^p\d t\d S^{n-1}
\;=\;\int\abs{u}^p\d t\d S^{n-1},
\]
where we have used the identity
\begin{equation}
\label{modulusidforpartuni}
\sum_{i\in I}\,\abs{\chi_i u}
\;=\;\abs u
\end{equation}
(recall that $\{\chi_i\}_{i\in I}$ is a partition of unity). From \eqref{simpleconjuncov} we have 
\[
(\Theta^{-1})^*\d t\d S^{n-1}
\;=\;r^{-1}\d r\d S^{n-1}
\;=\;r^{-n}\vtsp r^{n-1}\d r\d S^{n-1}
\;=\;\abs{x}^{-n}\d^nx
\]
so \eqref{defofdonspnorm} now gives 
\[
\norm{u}_{\donsp0{L^p}}^p
\;\asymp\;\int\abs{\Theta^* u}^p\d t\d S^{n-1}
\;=\;\int\abs{u}^p\vtsp(\Theta^{-1})^*(dt\d S^{n-1})
\;=\;\int\abs{x}^{-n}\abs{u}^p\d^n x
\]
for all measurable functions $u$ on $\don$. Finally Lemma \ref{equivofnorms1} gives
\begin{multline*}
\norm{u}_{\consp0{L^p}}^p
\;\asymp\;\norm{\eta_1u}_{L^p}^p+\norm{\zeta_1u}_{\donsp0{L^p}}^p\\
\asymp\;\int\bigl(\abs{\eta_0(x)}^p+\abs{x}^{-n}\abs{\zeta_0(x)}^p\bigr)
\abs{u(x)}^p\d^nx
\;\asymp\;\int\Lambda^{-n}(x)\vtsp\abs{u(x)}^p\d^n x
\end{multline*}
for all measurable functions $u$ on $\R^n$, where the last line follows from the existence of constants $C_1,C_2>0$ such that 
\[
C_1\Lambda^{-n}(x)
\;\le\;\abs{\eta_0(x)}^p+\abs{x}^{-n}\abs{\zeta_0(x)}^p
\;\le\;C_2\Lambda^{-n}(x)
\]
for all $x\in\R^n$.
\end{proof}

\medskip

A simpler version of the previous argument can be used to find an equivalent norm for the spaces $\consp\beta{\Cnt l}$ when $\beta\in\R$ and $l\in\NZ$.

\begin{prop}
\label{equivnormmodspCnt}
Suppose $\beta\in\R$ and $l\in\NZ$. Given $u\in\mathscr D'$ we have $u\in\consp\beta{\Cnt l}$ iff $u$ is $l$ times continuously differentiable and 
\begin{equation}
\label{sumgivealtnorm2}
\sum_{\abs\alpha\le l}\,\sup_{x\in\R^n}\, \Lambda^{\beta+\abs\alpha}(x)
\,\abs{\dfo_x^\alpha u(x)}
\;<\;+\infty.
\end{equation}
Furthermore the quantity on the left hand side of \eqref{sumgivealtnorm2} is equivalent to $\norm{u}_{\consp\beta{\Cnt l}}$.
\end{prop}

\begin{proof}
Once again Propositions \ref{basicchangbetaiso} and \ref{lemoncontofpderiv} reduce our task to proving the result in the case $\beta=0$ and $l=0$. 

In $u$ is a continuous function on $\cly$ then \eqref{defofEclynobeta} and \eqref{defofnormonclysp} give
\[
\norm{u}_{\clysp0{\Cnt 0}}
\ =\ \sum_{i\in I}\,\sup_{x\in\R\times\psi_i(U_i)}\bigabs{(\chi_iu)(\Psi_i^{-1}(x))}
\ \asymp\ \sup_{(t,\omega)\in\cly}\sum_{i\in I}\,
\bigabs{(\chi_i u)(t,\omega)}
\ =\ \sup_{(t,\omega)\in\cly}\,\abs{u(t,\omega)},
\]
where we have used the finiteness of $I$ and \eqref{modulusidforpartuni} in the second last and last steps respectively. From \eqref{defofdonspnorm} we now get
\[
\norm{u}_{\donsp0{\Cnt 0}}
\;\asymp\;\sup_{(t,\omega)\in\cly}\bigabs{u(\Theta(t,\omega))}
\;=\;\sup_{x\in\don}\abs{u(x)},
\]
for all continuous functions $u$ on $\don$, whilst Lemma \ref{equivofnorms1} finally gives
\begin{multline*}
\norm{u}_{\consp0{\Cnt 0}}
\;\asymp\;\sup_{x\in\R^n}\,\abs{\eta_0(x)u(x)}+\sup_{x\in\R^n}\,\abs{\zeta_0(x)u(x)}\\
\asymp\;\sup_{x\in\R^n}\,\bigl(\abs{\eta_0(x)u(x)}+\abs{\zeta_0(x)u(x)}\bigr)
\;=\;\sup_{x\in\R^n}\,\abs{u(x)}
\end{multline*}
for all continuous functions $u$ on $\R^n$.
\end{proof}

It is clear from the proof of Proposition \ref{equivnormmodspCnt} that we have 
\begin{equation}
\label{equivnormsdonspCnt}
\norm{u}_{\donsp\beta{\Cnt l}}
\;\asymp\;\sum_{\abs\alpha\le l}\,\sup_{x\in\don}
\,\abs x^{\beta+\abs\alpha}\,\abs{\dfo_x^\alpha u(x)}
\end{equation}
for any $\beta\in\R$ and $l\in\NZ$. This observation will be useful below.

\bigskip

For any $l\in\NZ$ and $\sigma\in(0,1)$ the Holder space $\Hol l\sigma$ can be defined as the collection of all those functions $u\in\Cnt l$ for which 
\[
\sum_{\abs\alpha\le l}\,\sup_{x\in\R^n}\,\abs{\dfo_x^\alpha u(x)}
\;+\;\sum_{\abs\alpha=l}\,\sup_{\substack{x,y\in\R^n\\ x\ne y}}\,
\bigabs{\dfo_x^\alpha u(x)\negvtsp-\negvtsp\dfo_x^\alpha u(y)}
\vtsp\abs{x\negvtsp-\negvtsp y}^{-\sigma}
\ <\ +\infty.
\]
We can use this sum to define the norm $\norm\cdot_{\Hol l\sigma}$ on $\Hol l\sigma$.

\begin{rem}
\label{altdefofHolnorm}
Suppose $V\subseteq\R^n$ is open and $\{B_x\}_{x\in V}$ is a collection of open subsets of $V$ for which there exists a constant $\kappa>0$ such that $V\cap\{y\,|\,\abs{x\negvtsp-\negvtsp y}\negvtsp<\negvtsp\kappa\}\subseteq B_x$ for all $x\in V$. It is straightforward to check that
\[
\norm{u}_{\Hol 0\sigma}
\;\asymp\;\sup_{x\in V}\,\abs{u(x)}
\;+\negthickspace\sup_{\substack{x\in V\\ y\in B_x\setminus\{x\}}}\!
\abs{u(x)\negvtsp-\negvtsp u(y)}\vtsp\abs{x\negvtsp-\negvtsp y}^{-\sigma}
\]
for all continuous functions $u$ with $\supp(u)\subseteq V$. 
\end{rem}

\medskip

As was the case for the model spaces $\Sob pk$ and $\Cnt l$ we can obtain an explicit description of the space $\consp\beta{\Hol l\sigma}$. 

\begin{prop}
\label{equivnormforHolcase}
Suppose $\beta\in\R$, $l\in\NZ$ and $\sigma\in(0,1)$. Given $u\in\mathscr D'$ we have $u\in\consp\beta{\Hol l\sigma}$ iff $u$ is $l$ times continuously differentiable and 
\begin{eqnarray}
\lefteqn{\sum_{\abs\alpha\le l}\,\sup_{x\in\R^n}\,\abs{(\Lambda^{\beta+\abs\alpha}\dfo_x^\alpha u)(x)}
}\qquad\quad\nonumber\\
\label{letsseehowthisgoes}
&+&\negthickspace\sum_{\abs\alpha=l}\,\sup_{\substack{x,y\in\R^n\\ x\ne y}}\,
\bigabs{(\Lambda^{\beta+l+\sigma}\dfo_x^\alpha u)(x)-(\Lambda^{\beta+l+\sigma}\dfo_x^\alpha u)(y)}
\vtsp\abs{x\negvtsp-\negvtsp y}^{-\sigma}
\ <\ +\infty.\ \quad
\end{eqnarray}
Furthermore the quantity on the left hand side of \eqref{letsseehowthisgoes} is equivalent to $\norm\cdot_{\consp\beta{\Hol l\sigma}}$.
\end{prop}

The proof of Proposition \ref{equivnormforHolcase} will be proceeded by some technical results for which we introduce the following notation; if $u$ is a continuous function on $\R^n$ let $\norm u_\sigma$ denote the left hand side of \eqref{letsseehowthisgoes} when $\beta=0$ and $l=0$.

\begin{lem}
\label{localtdeflemforHol}
For continuous functions $u$ on $\R^n$
\[
\norm{u}_\sigma
\;\asymp\;\sup_{x\in\R^n}\,\abs{u(x)}
\,\;+\!\negthickspace\sup_{\substack{x,y\in\R^n\\ 0<\abs{x-y}<\Lambda(x)}}
\negthickspace\Lambda^\sigma(x)\vtsp\abs{u(x)\negvtsp-\negvtsp u(y)}\vtsp
\abs{x\negvtsp-\negvtsp y}^{-\sigma}\,.
\]
\end{lem}

\begin{proof}
Since $\sigma\in(0,1)$ we have $1\negvtsp-\negvtsp t^\sigma\le 1\negvtsp-\negvtsp t\le(1\negvtsp-\negvtsp t)^\sigma$ for all $t\in[0,1]$. It follows that
\begin{equation}
\abs{\Lambda^\sigma(x)\negvtsp-\negvtsp \Lambda^\sigma(y)}
\;\le\;\abs{\Lambda(x)\negvtsp-\negvtsp \Lambda(y)}^\sigma
\;\le\;\abs{x\negvtsp-\negvtsp y}^\sigma
\end{equation}
for all $x,y\in\R^n$. If $u$ is any continuous function on $\R^n$ we thus have
\[
\Bigl\lvert{\vtsp\bigabs{(\Lambda^\sigma u)(x)\negvtsp-\negvtsp(\Lambda^\sigma u)(y)}
-\Lambda^\sigma(x)\vtsp\abs{u(x)\negvtsp-\negvtsp u(y)}\vtsp}\Bigr\rvert
\ \le\ \abs{\Lambda^\sigma(x)\negvtsp-\negvtsp \Lambda^\sigma(y)}\vtsp\abs{u(y)}
\ \le\ \abs{x-y}^\sigma\vtsp\abs{u(y)}
\]
for all $x,y\in\R^n$. Hence 
\[
\norm{u}_\sigma
\;\asymp\;\sup_{x\in\R^n}\,\abs{u(x)}
\,+\sup_{\substack{x,y\in\R^n\\ x\ne y}}
\Lambda^\sigma(x)\vtsp\abs{u(x)\negvtsp-\negvtsp u(y)}\vtsp\abs{x\negvtsp-\negvtsp y}^{-\sigma}.
\]
On the other hand $\abs{x\negvtsp-\negvtsp y}^{-\sigma}\le\Lambda^{-\sigma}(x)$ whenever $\abs{x\negvtsp-\negvtsp y}\ge\Lambda(x)$. Therefore
\[
\sup_{\substack{x,y\in\R^n\\ \abs{x-y}\ge\Lambda(x)}}
\Lambda^\sigma(x)\vtsp\abs{u(x)\negvtsp-\negvtsp u(y)}\vtsp\abs{x\negvtsp-\negvtsp y}^{-\sigma}
\ \le\ 2\sup_{x\in\R^n}\,\abs{u(x)}\,,
\]
completing the result.
\end{proof}

\begin{lem}
\label{locmultilemforHol}
We have $\norm{uv}_\sigma\le C\norm{u}_{\consp0{\Cnt 1}}\norm{v}_\sigma$ for all $u\in\consp0{\Cnt 1}$ and continuous functions $v$ on $\R^n$.
\end{lem}

\begin{proof}
For any $u\in\consp0{\Cnt 1}$ Proposition \ref{equivnormmodspCnt} gives us
\begin{equation}
\label{fisrest1}
\abs{u(x)}\;\le\;C\norm{u}_{\consp0{\Cnt 1}}
\end{equation}
and
\begin{equation}
\label{fisrest15}
\abs{\dfo_i u(x)}\;\le\;C\norm{u}_{\consp0{\Cnt 1}}\vtsp\Lambda^{-1}(x),
\quad i=1,\dots,n,
\end{equation}
for all $x\in\R^n$. Now suppose $x,y\in\R^n$ with $0<\abs{x\negvtsp-\negvtsp y}\le\Lambda(y)$. Thus $x$ and $y$ both belong to the ball of radius $2\Lambda(y)$ centred at the origin. Combining this observation with \eqref{fisrest15} and the fact that $\Lambda^{\sigma-1}(y)\abs{x\negvtsp-\negvtsp y}^{1-\sigma}\le 1$ we then get
\[
\Lambda^{\sigma}(y)\vtsp\abs{u(x)\negvtsp-\negvtsp u(y)}\vtsp\abs{x\negvtsp-\negvtsp y}^{-\sigma}
\;\le\;C\Lambda^\sigma(y)\abs{x\negvtsp-\negvtsp y}^{1-\sigma}\,
\norm{u}_{\consp0{\Cnt 1}}\vtsp\Lambda^{-1}(y)
\;\le\;C\norm{u}_{\consp0{\Cnt 1}}.
\]
On the other hand, if $x,y\in\R^n$ with $\abs{x\negvtsp-\negvtsp y}\ge\Lambda(y)$ then $\Lambda^{\sigma}(y)\vtsp\abs{x\negvtsp-\negvtsp y}^{-\sigma}\le 1$ so
\[
\Lambda^{\sigma}(y)\vtsp\abs{u(x)\negvtsp-\negvtsp u(y)}\vtsp\abs{x\negvtsp-\negvtsp y}^{-\sigma}
\;\le\;C\norm{u}_{\consp0{\Cnt 1}}
\]
by \eqref{fisrest1}. Combining the above estimates we thus get
\begin{equation}
\label{fisrest2}
\Lambda^{\sigma}(y)\vtsp\abs{u(x)\negvtsp-\negvtsp u(y)}\vtsp\abs{x\negvtsp-\negvtsp y}^{-\sigma}
\;\le\;C\norm{u}_{\consp0{\Cnt 1}}
\end{equation}
for all $x,y\in\R^n$ with $x\ne y$. 

Suppose $v$ is a continuous function and $x,y\in\R^n$ with $x\ne y$. Using \eqref{fisrest1} and \eqref{fisrest2} we thus have
\begin{eqnarray*}
\lefteqn{
\bigabs{(\Lambda^\sigma uv)(x)\negvtsp-\negvtsp (\Lambda^\sigma uv)(y)}\vtsp
\abs{x\negvtsp-\negvtsp y}^{-\sigma}}\quad\\
&\le&\abs{u(x)}\vtsp\bigabs{(\Lambda^\sigma v)(x)\negvtsp-\negvtsp (\Lambda^\sigma v)(y)}
\vtsp\abs{x\negvtsp-\negvtsp y}^{-\sigma}
+\abs{v(y)}\,\Lambda^{\sigma}(y)\vtsp\abs{u(x)\negvtsp-\negvtsp u(y)}\vtsp
\abs{x\negvtsp-\negvtsp y}^{-\sigma}\\
&\le&C\norm{u}_{\consp0{\Cnt 1}}\bigl(\bigabs{(\Lambda^\sigma v)(x)
\negvtsp-\negvtsp (\Lambda^\sigma v)(y)}\vtsp\abs{x\negvtsp-\negvtsp y}^{-\sigma}
+\abs{v(y)}\bigr).
\end{eqnarray*}
The result now follows from the definition of $\norm\cdot_\sigma$.
\end{proof}

\begin{proof}[Proof of Proposition \ref{equivnormforHolcase}]
Propositions \ref{basicchangbetaiso} and \ref{lemoncontofpderiv} reduce our task to proving the result in the case $\beta=0$ and $l=0$. 

Consider the notation introduced in the first paragraph of Section \ref{subsecweightfunspcly} and let $\Theta$ be as given in Section \ref{wfnspondon}. Now, for each $i\in I$, define open sets $V_i,W_i\subset\R^n$ by $V_i=(0,\infty)\times\psi_i(U_i)$ and $W_i=\Theta(\Psi_i^{-1}(V_i))$. Thus the map $\Phi_i:=\Theta\circ\Psi_i^{-1}:V_i\to W_i$ is a diffeomorphism. Also define a function $\rho_i$ by $\rho_i(x)=\zeta_1(x)\chi_i(\Phi_i^{-1}(x))$ for $x\in W_i$ and $\rho_i(x)=0$ for $x\notin W_i$. It is easy to see that $\rho_i$ is smooth and independent of $\abs{x}$ for sufficiently large $\abs{x}$; it follows from Proposition \ref{equivnormmodspCnt} that $\rho_i\in\consp0{\Cnt 1}$. Finally set $V=\{\abs{x}<3\}\subset\R^n$, so $\eta_1\in C^\infty_0(V)\subset\consp0{\Cnt 1}$. 

Since $\{\chi_i\}_{i\in I}$ is a partition of unity on $\cly$ we have $\eta_1+\sum_{i\in I}\rho_i=\eta_1+\zeta_1=1$. Together with Lemma \ref{locmultilemforHol} and the fact that $I$ is finite we then get 
\[
\norm{u}_\sigma\;\asymp\;\norm{\eta_1u}_\sigma
+\sum_{i\in I}\,\norm{\rho_iu}_\sigma
\]
for all continuous functions $u$ on $\R^n$. On the other hand, the definition of $\consp0{\Hol 0\sigma}$ means we have 
\[
\norm{u}_{\consp0{\Hol 0\sigma}}\;\asymp\;\norm{\eta_1u}_{\Hol 0\sigma}
+\sum_{i\in I}\,\norm{\Phi^*_i(\rho_i u)}_{\Hol 0\sigma}
\]
for all $u\in\consp0{\Hol 0\sigma}$. The following claims thus complete the result.

\begin{claim}[Claim(i): We have $\norm{v}_{\Hol 0\sigma}\asymp\norm{v}_\sigma$ for all continuous functions $v$ with $\supp(v)\subseteq V$.]
This is a straightforward consequence of Remark \ref{altdefofHolnorm}, Lemma \ref{localtdeflemforHol} and the fact that $\Lambda(x)\asymp 1$ for $x\in V$ (n.b.~$V$ is bounded).
\end{claim}

\begin{claim}[Claim(ii): If $i\in I$ then $\norm{\Phi^*_i v}_{\Hol 0\sigma}\asymp\norm{v}_\sigma$ for all continuous functions $v$ with $\supp(v)\subseteq V_i$.]
If we write $y\in\R^n$ in the form $y=(t,w)$ with $t\in\R$ and $w\in\R^{n-1}$ then $\Phi_i(y)=e^t\psi_i^{-1}(w)$ (where we are considering $S^{n-1}$ to be the unit sphere in $\R^n$). Thus, for all $y_1,y_2\in V_i$ with $\abs{y_1-y_2}<1$,
\begin{eqnarray}
\abs{\Phi_i(y_1)\negvtsp-\negvtsp\Phi_i(y_2)}
&\asymp&e^{t_1}\bigl(\abs{1\negvtsp-\negvtsp e^{t_2-t_1}}\,\abs{\psi_i^{-1}(w_1)}
+\bigabs{\psi_i^{-1}(w_1)\negvtsp-\negvtsp\psi_i^{-1}(w_2)}\bigr)
\nonumber\\
&\asymp&e^{t_1}\bigl(\abs{t_1\negvtsp-\negvtsp t_2}+\abs{w_1\negvtsp-\negvtsp w_2}\bigr)
\nonumber\\
\label{aranarwego}
&\asymp&\Lambda(\Phi_i(y_1))\vtsp\abs{y_1\negvtsp-\negvtsp y_2},
\end{eqnarray}
where the inequalities $\abs{t_1\negvtsp-\negvtsp t_2}\le\abs{y_1\negvtsp-\negvtsp y_2}<1$ and $\abs{\Phi_i(y_1)}=e^{t_1}\ge1$ have been used in the second last and last lines respectively. For each $y_1\in V_i$ set
\[
B_{y_1}\;=\;\bigl\{y_2\in V_i\;\big|\;\abs{\Phi_i(y_1)\negvtsp-\negvtsp\Phi_i(y_2)}
<\Lambda(\Phi(y_1))\bigr\}
\]
Estimate \eqref{aranarwego} implies there exists $\kappa>0$ such that $V_i\cap\{y_2\,|\,\abs{y_1\negvtsp-\negvtsp y_2}\negvtsp<\negvtsp\kappa\}\subseteq B_{y_1}$ for all $y_1\in V_i$. Further use of \eqref{aranarwego} together with Remark \ref{altdefofHolnorm} and Lemma \ref{localtdeflemforHol} then gives
\begin{eqnarray*}
\norm{\Phi_i^*v}_{\Hol 0\sigma}
&\asymp&\sup_{y\in V_i}\,\abs{v(\Phi_i(y))}
\;+\negthickspace\sup_{\substack{y_1\in V_i\\ y_2\in B_{y_1}\!\setminus\{y_1\}}}\!
\bigabs{v(\Phi_i(y_1))\negvtsp-\negvtsp v(\Phi_i(y_2))}\vtsp\abs{y_1\negvtsp-\negvtsp y_2}^{-\sigma}\\
&\asymp&\sup_{x\in W_i}\,\abs{v(x)}
\;+\negthickspace
\sup_{\substack{x_1,x_2\in W_i\\ 0<\abs{x_1\negvtsp-\negvtsp x_2}<\Lambda(x_1)}}\!
\Lambda^\sigma(x_1)\,\abs{v(x_1)\negvtsp-\negvtsp v(x_2)}\,\abs{x_1\negvtsp-\negvtsp x_2}^{-\sigma}\\
&\asymp&\norm{v}_\sigma,
\end{eqnarray*}
for all continuous functions $v$ with $\supp(v)\subset V_i$.
\end{claim}
\end{proof}

\begin{rem}
Let $\beta\in\R$, $l\in\NZ$ and $\sigma\in(0,1)$. Using Proposition \ref{equivnormforHolcase} it can be seen that $\consp\beta{\Hol l\sigma}$ coincides with the space $C^{\sigma+l}_\beta(\R^n)$ defined in \cite{Be}. Also, using Proposition \ref{basicchangbetaiso}, Lemma \ref{denofCinft0incertsp} and the obvious modification of Proposition \ref{equivnormforHolcase} for the spaces $\donsp\beta{\Hol l\sigma}$, it can be seen that $\donsp\beta{\Hol l\sigma_0}$ coincides with the space $\Lambda^{l,\sigma}_{\beta+l+\sigma}(\don)$ defined in Section 3.6.4 of \cite{NP2} (here $\Hol l\sigma_0$ denotes the separable subspace of $\Hol l\sigma$ obtained by taking the completion of $C^\infty_0$). The non-separable space $\donsp\beta{\Hol l\sigma}$ contains elements which behave as $O(\abs x^{-\beta})$ for $\abs x\to0,\infty$ and is strictly larger than $\donsp\beta{\Hol l\sigma_0}$ (see Section~\ref{someresrelsymsssec} for a related discussion). 
\end{rem}

\subsubsection{Dual spaces}

Suppose $E$ is a model space. By definition $\mathscr S$ is dense in $E_0$ so any element in the dual of $E_0$ is uniquely determined by its action on $\mathscr S$. We can thus uniquely identify elements of the dual of $E_0$ with tempered distributions on $\R^n$; i.e.~elements of $\mathscr S'$. Furthermore a norm can be defined on $E_0^*$ by the expression
\begin{equation}
\label{defofdualnorm}
\norm{v}_{E_0^*}
\;=\;\sup_{0\ne u\in\mathscr S}\frac{\abs{(u,v)_{\R^n}}}{\norm{u}_E}
\;=\;\sup_{0\ne u\in C^\infty_0}\frac{\abs{(u,v)_{\R^n}}}{\norm{u}_E}\,,
\end{equation}
where the second equality follows from the density of $C^\infty_0$ in $E_0$ (see Remark \ref{defofsssrem}). Applying standard duality arguments to Conditions (A1) to (A4) we obtain the following.

\begin{lem}
\label{EmodimpliesE0stmod}
If $E$ is a model space then so is $E_0^*$.
\end{lem}

If $E$ and $F$ are model spaces we shall write $E_0^*=F$ provided these spaces agree as subsets of $\mathscr S'$ and $\norm{\cdot}_F$ is equivalent to the norm given by \eqref{defofdualnorm}. In other words $E_0^*=F$ iff there are constants $C_1,C_2>0$ such that the following hold.
\begin{list}{}{\listinit}
\item[(D1)\hf]
For each $u\in C^\infty_0$ and $v\in F$ we have $\abs{(u,v)_{\R^n}}\le C_1\norm{u}_E\norm{v}_F$.
\item[(D2)\hf]
For each $v\in F$ there exists $0\ne u\in C^\infty_0$ with $\abs{(u,v)_{\R^n}}\ge C_2\norm{u}_E\norm{v}_F$.
\end{list}

\bigskip

Suppose $E$ is a model space and $\beta\in\R$. Using Lemma \ref{denofCinft0incertsp} and an argument similar to that above, we can identify the dual spaces $(\clysp\beta{E_0})^*$, $(\donsp\beta{E_0})^*$ and $(\consp\beta{E_0})^*$ with subspaces of $\mathscr D'(\cly)$, $\mathscr D'(\don)$ and $\mathscr D'$ respectively. The pairings $(\cdot,\cdot)_\cly$, $(\cdot,\cdot)_\don$ and $(\cdot,\cdot)_{\R^n}$ then allow us to define norms on these dual spaces (as in the second part of \eqref{defofdualnorm}) and compare them with existing spaces.

\begin{rem}
\label{propsA1A4innewform}
For any $\beta\in\R$ the pairing $(\cdot,\cdot)_{\cly}$ on $\cly$ extends to a dual pairing $\clysp\beta{\mathscr S'}\times\clysp{-\beta}{\mathscr S}\to\C$.
\end{rem}

\begin{lem}
\label{dualforclysps}
If $E$ is a model space and $\beta\in\R$ then $(\clysp\beta{E_0})^*=\clysp{-\beta}{(E_0^*)}$.
\end{lem}

Here `equality' is understood in the sense of equivalent norms; that is, we have expressions similar to (D1) and (D2) above.

\begin{proof}
Let $F=E_0^*$ and consider the notation of Remark \ref{enlargechirem}. Now, for $i\in I$, the pull-back under $\Psi_i^{-1}$ of the density $\chi_i'\d t\d S^{n-1}$ is a smooth density on $\R^n$. Thus we can write $(\Psi_i^{-1})^*(\chi_i'\d t\d S^{n-1})=J_i\d^nx$ where $J_i\in\Cnt\infty$ (in fact $J_i$ is independent of $t$ and compactly supported in the remaining variables). It follows that we have
\[
(u,\chi_iv)_\cly
\;=\;\bigl(\chi_i'\vtsp\chi_i'e^{\beta t}u,\,\chi_i e^{-\beta t}v\bigr)_{\R^n}
\;=\;\bigl(J_i(\Psi_i^{-1})^*(\chi_i'e^{\beta t}u),
\,(\Psi_i^{-1})^*(\chi_i e^{-\beta t}v)\bigr)_{\R^n}
\]
for all $u\in C^\infty_0(\cly)$ and $v\in\mathscr D'(\cly)$. Using (D1), (A2) (to deal with $J_i$), Remark \ref{enlargechirem} and \eqref{defofnormonclysp}, we therefore have
\begin{eqnarray*}
\abs{(u,v)_\cly}
&\le&\sum_{i\in I}\,\abs{(u,\chi_iv)_\cly}\\
&\le&C\sum_{i\in I}\,\bignorm{J_i(\Psi_i^{-1})^*(\chi_i'e^{\beta t}u)}_E\,
\bignorm{(\Psi_i^{-1})^*(\chi_i e^{-\beta t}v)}_F\\
&\le&C\sum_{i\in I}\,\bignorm{(\Psi_i^{-1})^*(\chi_i'e^{\beta t}u)}_E\;
\sum_{i\in I}\,\bignorm{(\Psi_i^{-1})^*(\chi_i e^{-\beta t}v)}_F\\
&\le&C\norm{u}_{\clysp\beta E}\vtsp\norm{v}_{\clysp{-\beta}F}
\end{eqnarray*}
for all $u\in C^\infty_0(\cly)$ and $v\in\clysp{-\beta}F$.

\medskip

Now let $v\in\clysp{-\beta}F$. By definition
\[
\norm{v}_{\clysp{-\beta}F}
\;=\;\sum_{i\in I}\,\bignorm{(\Psi_i^{-1})^*(\chi_i e^{-\beta t}v)}_F\,.
\]
Hence we can find $i\in I$ such that
\[
\norm{v}_{\clysp{-\beta}F}
\;\le\;C_1\bignorm{(\Psi_i^{-1})^*(\chi_i e^{-\beta t}v)}_F\,;
\]
n.b.~$i$ may depend on $v$ but $C_1$ does not (we can define $C_1$ to be the number of elements in $I$). Using (D2) we can now choose $0\ne\phi\in C^\infty_0$ such that
\begin{equation}
\label{titfbev1}
\bigabs{\bigl(\phi,\,(\Psi_i^{-1})^*(\chi_i e^{-\beta t}v)\bigr)_{\R^n}}
\ \ge\ C_2\norm{\phi}_E\vtsp\bignorm{(\Psi_i^{-1})^*(\chi_i e^{-\beta t}v)}_F
\ \ge\ C_1^{-1}C_2\norm{\phi}_E\vtsp\norm{v}_{\clysp{-\beta}F}.
\end{equation}
Define $J'_i\in C^\infty(\cly)$ by $\chi_i'\Psi_i^*(d^nx)=J'_i\vtsp d t\d S^{n-1}$ and set $u=\chi_i J_i'e^{-\beta t}\Psi_i^*\phi\in C^\infty_0(\cly)$. Using \eqref{defofnormonclysp}, Lemma \ref{cordequivnorm} and Condition (A2) (n.b.~$(\Psi_i^{-1})^*(\chi_iJ_i')\in\Cnt\infty$), it follows that
\begin{equation}
\label{titfbev2}
\norm{u}_{\clysp\beta E}
\;=\;\norm{\chi_iJ_i'\Psi_i^*\phi}_{E(\cly)}
\;\le\;C\norm{(\Psi_i^{-1})^*(\chi_iJ_i')\,\phi}_E
\;\le\;C_3\norm{\phi}_E.
\end{equation}
On the other hand the definition of $J_i'$ gives us
\begin{equation}
(\label{titfbev3}
u,v)_\cly
\ =\ \bigl(J_i' \chi_ie^{-\beta t}\Psi_i^*\phi,\,v\bigr)_\cly
\ =\ \bigl(J_i' \Psi_i^*\phi,\,\chi_ie^{-\beta t}v\bigr)_\cly
\ =\ \bigl(\phi,\,(\Psi_i^{-1})^*(\chi_ie^{-\beta t}v)\bigr)_{\R^n}.
\end{equation}
Combining \eqref{titfbev1}, \eqref{titfbev2} and \eqref{titfbev3}, we then get
\[
\abs{(u,v)_\cly}
\;\ge\;C_1^{-1}C_2\norm{\phi}_E\vtsp\norm{v}_{\clysp{-\beta}F}
\;\ge\;C_1^{-1}C_2C_3^{-1}\norm{u}_{\clysp\beta E}\vtsp\norm{v}_{\clysp{-\beta}F}.
\]
This completes the result.
\end{proof}

Under the pull back induced by the diffeomorphism $\Theta:\cly\to\don$ (from Section \ref{wfnspondon}) we have $\Theta^*\d^nx=e^{nt}\d t\d S^{n-1}$. It follows that $(u,v)_\don=(\Theta^*u,\vtsp e^{nt}\Theta^*v)_\cly$ for any $u\in C^\infty_0(\don)$ and $v\in\mathscr D'(\don)$. Now suppose $E$ is a model space and $\beta\in\R$. With the help of \eqref{defofnormonclysp} and \eqref{Thetaisaniso}, Lemma \ref{dualforclysps} then gives
\begin{multline*}
\norm{v}_{\donsp{n-\beta}{(E_0^*)}}
\ =\ \norm{e^{nt}\Theta^*v}_{\clysp{-\beta}{(E_0^*)}}
\ \asymp\ \sup_{0\ne w\in C^\infty_0(\cly)}
\frac{\bigabs{(w,\vtsp e^{nt}\Theta^*v)_\cly}}{\norm{w}_{\clysp\beta E}}\\
=\ \sup_{0\ne u\in C^\infty_0(\don)}
\frac{\bigabs{\bigl(\Theta^*u,\vtsp e^{nt}\Theta^*v\bigr)_\cly}}{\norm{\Theta^* u}_{\clysp\beta E}}
\ =\ \sup_{0\ne u\in C^\infty_0(\don)}
\frac{\abs{(u,v)_\don}}{\norm{u}_{\donsp\beta E}}
\ =\ \norm{v}_{(\donsp\beta{E_0})^*}
\end{multline*}
for all $v\in\donsp{n-\beta}{(E_0^*)}$; that is, we have  
\begin{equation}
\label{dualeqnfordonsp}
(\donsp\beta{E_0})^*\;=\;\donsp{n-\beta}{(E_0^*)}.
\end{equation}

\begin{prop}
\label{dualpropforconsp}
If $E$ is a model space and $\beta\in\R$ then $(\consp\beta{E_0})^*=\consp{n-\beta}{(E_0^*)}$. 
\end{prop}

\begin{proof}
Let $F=E_0^*$ and choose $0\ne v\in\consp{n-\beta}F$. Now $\eta_2\eta_1=\eta_1$, $\zeta_0\zeta_1=\zeta_1$ and $\eta_2+\zeta_0\ge\eta_1+\zeta_1=1$ so, for any $u\in C^\infty_0$,
\[
(u,v)_{\R^n}
\,=\,(\eta_2 u,\eta_1 v)_{\R^n}+(\zeta_0 u,\zeta_1 v)_{\don}
\]
while
\begin{equation}
\label{theOMGpf5}
\norm{u}_{\consp\beta E}
\,\asymp\,\norm{\eta_2u}_E+\norm{\zeta_0u}_{\donsp\beta E}
\qquad\text{and}\qquad
\norm{v}_{\consp{n-\beta}F}
\,\asymp\,\norm{\eta_1v}_F+\norm{\zeta_1v}_{\donsp{n-\beta}F},
\end{equation}
by Lemma \ref{equivofnorms1}. Together with \eqref{dualeqnfordonsp} we then get
\begin{eqnarray*}
\abs{(u,v)_{\R^n}}
&\le&\abs{(\eta_2 u,\eta_1 v)_{\R^n}}+\abs{(\zeta_0 u,\zeta_1 v)_{\don}}\\
&\le&C\bigl(\norm{\eta_2u}_E\norm{\eta_1u}_F
+\norm{\zeta_0u}_{\donsp\beta E}\norm{\zeta_1v}_{\donsp{n-\beta}F}\bigr)\\
&\le&C\bigl(\norm{\eta_2u}_E+\norm{\zeta_0u}_{\donsp\beta E}\bigr)
\bigl(\norm{\eta_1v}_F+\norm{\zeta_1v}_{\donsp{n-\beta}F}\bigr)\\
&\le&C\norm{u}_{\consp\beta E}\norm{v}_{\consp{n-\beta}F}
\end{eqnarray*}
for any $u\in C^\infty_0$.

\smallskip

By the second part of \eqref{theOMGpf5} we have 
\[
\norm{v}_{\consp{n-\beta}F}
\;\le\;C_1\bigl(\norm{\eta_1v}_F+\norm{\zeta_1v}_{\donsp{n-\beta}F}\bigr). 
\]
Consider the following cases. 

\begin{claim}[Case (i): $\norm{v}_{\consp{n-\beta}F}\le2C_1\norm{\eta_1v}_F$.]
Using (D2) choose $w\in C^\infty_0$ with $\norm{w}_{E}=1$ and $\norm{v}_{\consp{n-\beta}F}\le 4C_1\abs{(w,\eta_1v)_{\R^n}}$. Set $u=\eta_1w$ so $u\in C^\infty_0$, $\norm{v}_{\consp{n-\beta}F}\le 4C_1\abs{(u,v)_{\R^n}}$ and 
\[
\norm{u}_{\consp\beta E}
\;\le\;C\bigl(\norm{\eta_2u}_E+\norm{\zeta_2u}_{\donsp\beta E}\bigr)
\;=\;C\norm{\eta_1 w}_E
\;\le\;C_2,
\]
by Lemma \ref{equivofnorms1}, where $C_2$ is independent of $u$ and $v$. 
\end{claim}

\begin{claim}[Case (ii): $\norm{v}_{\consp{n-\beta}F}\le2C_1\norm{\zeta_1v}_{\donsp{n-\beta}F}$.]
Using \eqref{dualeqnfordonsp} choose $w\in C^\infty_0(\don)$ with $\norm{w}_{\donsp\beta E}=1$ and $\norm{v}_{\consp{n-\beta}F}\le 4C_1\abs{(w,\zeta_1v)_{\R^n}}$. Set $u=\zeta_1w$ so $u\in C^\infty_0$, $\norm{v}_{\consp{n-\beta}F}\le 4C_1\abs{(u,v)_{\R^n}}$ and 
\[
\norm{u}_{\consp\beta E}
\;\le\;C\bigl(\norm{\eta_0u}_E+\norm{\zeta_0u}_{\donsp\beta E}\bigr)
\;=\;C\norm{\zeta_1 w}_{\donsp\beta E}
\;\le\;C_2,
\]
by Lemma \ref{equivofnorms1}, where $C_2$ is independent of $u$ and $v$. 
\end{claim}

\smallskip

By combining the 2 cases it follows that we can find $0\ne u\in C^\infty_0$ with 
\[
\norm{u}_{\consp\beta E}\,\norm{v}_{\consp{n-\beta}F}
\;\le\;4C_1C_2\abs{(u,v)_{\R^n}}.
\]
This completes the result.
\end{proof}

\subsubsection{Multiplication}

\begin{prop}
\label{multionallsp}
Suppose multiplication defines a continuous bilinear map $E\times F\to G$ for some model spaces $E$, $F$ and $G$. Then multiplication also defines continuous bilinear maps $\clysp\beta E\times\clysp\gamma F\to\clysp{\beta+\gamma} G$, $\donsp\beta E\times\donsp\gamma F\to\donsp{\beta+\gamma} G$ and $\consp\beta E\times\consp\gamma F\to\consp{\beta+\gamma} G$ for any $\beta,\gamma\in\R$.
\end{prop}

\begin{proof}
We can prove multiplication defines a continuous bilinear map $E(\cly)\times F(\cly)\to G(\cly)$ by an argument identical to that given for Lemma \ref{PropA2clyset}. The first two parts of the result now follow from \eqref{defofnormonclysp} and \eqref{Thetaisaniso} respectively.

\smallskip

Now $\eta_2\eta_1=\eta_1$ and $\zeta_1\zeta_0=\zeta_1$ so $\zeta_1u\vtsp\zeta_0v=\zeta_1uv$ and $\eta_2u\vtsp\eta_1v=\zeta_1uv$. However $\eta_2+\zeta_1,\eta_1+\zeta_0\ge 1$, so Lemma \ref{equivofnorms1} and the continuity of multiplication as a map $E\times F\to G$ and as a map $\donsp\beta E\times\donsp\gamma F\to\donsp{\beta+\gamma}G$ give
\begin{eqnarray*}
\norm{uv}_{\consp{\beta+\gamma}G}
&\le&C\bigl(\norm{\eta_1 uv}_{G}+\norm{\zeta_1 uv}_{\donsp{\beta+\gamma}G}\bigr)\\
&\le&C\bigl(\norm{\eta_2 u}_{E}\norm{\eta_1 v}_{F}
+\norm{\zeta_1 u}_{\donsp\beta E}\norm{\zeta_0 v}_{\donsp\gamma F}\bigr)\\
&\le&C\bigl(\norm{\eta_2 u}_{E}+\norm{\zeta_1 u}_{\donsp\beta E}\bigr)
\bigl(\norm{\eta_1 v}_{F}+\norm{\zeta_0 v}_{\donsp\gamma F}\bigr)\\
&\le&C\norm{u}_{\consp\beta E}\norm{v}_{\consp\gamma F}
\end{eqnarray*}
for any $u\in\consp\beta E$ and $v\in\consp\gamma F$.
\end{proof}

\begin{rem}
Suppose multiplication defines a continuous bilinear map $E\times F\to G$ for some model spaces $E$, $F$ and $G$. As a straightforward consequence of Definition \ref{defnofElver1} and the Leibniz rule we immediately have that multiplication also defines a continuous bilinear map $E^l\times F^l\to G^l$ for any $l\in\NZ$. 

Let $E'$ denote the closure of $E\cap C^\infty_\loc$ in $E$ with the induced norm (n.b.~we have $E_0\subseteq E'\subseteq E$ although both inclusions could be strict in general). If $u\in E\cap C^\infty_\loc$ and $v\in G_0^*$ then $(uv,f)_{\R^n}=(v,fu)_{\R^n}$ and $\norm{uf}_G\le C\norm{u}_E\norm{f}_F$ for all $f\in C^\infty_0$. Using \eqref{defofdualnorm} we thus get 
\[
\norm{uv}_{F_0^*}
\;=\;\sup_{0\ne f\in C^\infty_0}\frac{\abs{(uv,f)_{\R^n}}}{\norm{f}_{F}}
\;\le\;C\norm{u}_E\!\sup_{\substack{f\in C^\infty_0\\ uf\ne 0}}\frac{\abs{(v,uf)_{\R^n}}}{\norm{uf}_G}
\;=\;C\norm{u}_E\vtsp\norm{v}_{G_0^*}.
\]
Taking completions we finally have that multiplication defines a continuous bilinear map $E'\times G_0^*\to F_0^*$.
\end{rem}

\bigskip

For any model space $E$, Condition (A2) ensures that multiplication defines a continuous bilinear map $\Cnt\infty\times E\to E$. Now the topology on $\Cnt\infty$ is defined by the collection of semi-norms $\{\norm\cdot_{\Cnt l}\,|\,l\in\NZ\}$. Using Lemma \ref{techlemmfromfana} we can thus find $l\in\NZ$ and a constant $C$ such that $\norm{\phi u}_E\le C\norm{\phi}_{\Cnt l}\norm{u}_E\le C\norm{\phi}_{\Cnt{l+1}}\norm{u}_E$ for all $u\in E$ and $\phi\in\Cnt\infty$. The fact that the closure of $\Cnt\infty$ in $\Cnt l$ includes $\Cnt{l+1}$ now completes the following result.

\begin{lem}
\label{lemaboutmultibyCl}
For all sufficiently large $l\in\NZ$, multiplication defines a continuous bilinear map $\Cnt l\times E\to E$. 
\end{lem}

We conclude this section with an immediate consequence of Proposition \ref{multionallsp} and Lemma \ref{lemaboutmultibyCl} which will be used later in dealing with symbols.

\begin{prop}
\label{multibysymprop}
Suppose $E$ is a model space and $\beta,\gamma\in\R$. Then, for all sufficiently large $l\in\NZ$, multiplication defines continuous bilinear maps $\clysp\gamma{\Cnt{l}}\times\clysp\beta E\to\clysp{\beta+\gamma}E$, $\donsp\gamma{\Cnt{l}}\times\donsp\beta E\to\donsp{\beta+\gamma}E$ and $\consp\gamma{\Cnt{l}}\times\consp\beta E\to\consp{\beta+\gamma}E$.
\end{prop}

\subsubsection{Compactness}
Let $E$ and $F$ be model spaces. We say \emph{$E$ is locally compact in $F$} if multiplication by any $\phi\in C^\infty_0$ defines a compact map $E\to F$; in particular it follows that we have a local inclusion $E_\loc\hookrightarrow F_\loc$.

\begin{prop}
\label{genformofcpt}
Suppose $E$, $F$ and $G$ are model spaces with $E$ locally compact in $G$ and for which multiplication defines a continuous bilinear map $E\times F\to G$. Let $\beta,\gamma\in\R$ and $v\in\consp\gamma{F_0}$. Then multiplication by $v$ defines a compact map $\consp\beta E\to\consp{\beta+\gamma}G$. 
\end{prop}

\begin{proof}
Initially suppose $v\in C^\infty_0$. Let $\{u_i\}_{i\in\N}\subset\consp\beta E$ be a bounded sequence. By Lemma \ref{equivofnorms1} it follows that $\{v u_i\}_{i\in\N}$ is a bounded sequence in $E$. By local compactness we can thus find a subsequence $\{v u_{i(j)}\}_{j\in\N}$ which is convergent in $G$. Lemma \ref{equivofnorms1} then implies this subsequence must also be convergent in $\consp{\beta+\gamma}G$. It follows that multiplication by $v$ defines a compact map $\consp\beta E\to\consp{\beta+\gamma}G$.

Now let $v$ be an arbitrary element of $\consp\gamma{F_0}$. Let $\epsilon>0$ and, using Lemma~\ref{denofCinft0incertsp}, choose $v_\epsilon\in C^\infty_0$ with $\norm{v-v_\epsilon}_{\consp\gamma F}<\epsilon$. Proposition \ref{multionallsp} then implies multiplication by $v-v_\epsilon$ defines a map in $\mathscr L(\consp\beta E,\consp{\beta+\gamma}G)$ with norm at most $C\epsilon$, where $C$ is independent of $\epsilon$. The result now follows from the fact that the set of compact maps in $\mathscr L(\consp\beta E,\consp{\beta+\gamma}G)$ is closed (see Theorem~III.4.7 in \cite{Ka} for example).
\end{proof}

\subsubsection{Some results relating to symbols}

\label{someresrelsymsssec}

Let $l\in\NZ$. Clearly $\phi_0\in C^\infty(\cly)\subset\Cnt l(\cly)$, whilst $\phi_i$ is simply a translation of $\phi_0$ for any $i\in\Z$. Condition (A3) for $\Cnt l(\cly)=\clysp0{\Cnt l}$ then implies $\norm{\phi_i}_{\clysp0{\Cnt l}}$ is independent of $i\in\Z$. However $\zeta_i=\phi_i\circ\Theta^{-1}$ (by definition) so \eqref{defofdonspnorm} now implies $\norm{\zeta_i}_{\donsp0{\Cnt l}}$ is also independent of $i\in\Z$. Since $1\in\consp0{\Cnt l}$ (by Proposition \ref{equivnormmodspCnt}) Lemma \ref{equivofnorms2} now completes the following result.

\begin{lem}
\label{lemnoneoflsec}
Let $l\in\NZ$. Then $\norm{\zeta_i}_{\donsp0{\Cnt l}}$ and $\norm{\zeta_i}_{\consp0{\Cnt l}}$ are bounded uniformly for $i\in\N$. 
\end{lem}

Proposition \ref{multionallsp} and the fact that multiplication defines a continuous bilinear map $\Cnt l\times\Cnt l\to\Cnt l$ shows that multiplication also defines a continuous bilinear map $\consp0{\Cnt l}\times\consp\beta{\Cnt l}\to\consp\beta{\Cnt l}$ for any $\beta\in\R$. Given $u\in\consp\beta{\Cnt l}$, Lemmas \ref{equivofnorms2} and \ref{lemnoneoflsec} now imply $\norm{\zeta_i u}_{\consp\beta{\Cnt l}}$ and $\norm{\zeta_i u}_{\donsp\beta{\Cnt l}}$ are bounded uniformly for $i\in\N$.

\begin{lem}
\label{lemntwooflsec}
Let $l\in\NZ$, $\beta\in\R$ and $u\in\consp\beta{\Cnt l}$. Then the following are equivalent.
\begin{list}{}{\listinit}
\item[\emph{(i)}\hf]
$\norm{\zeta_i u}_{\consp\beta{\Cnt l}}\to0$ as $i\to\infty$.
\item[\emph{(ii)}\hf]
$\norm{\zeta_i u}_{\donsp\beta{\Cnt l}}\to0$ as $i\to\infty$.
\item[\emph{(iii)}\hf]
$\dfo_x^\alpha u(x)=o(\abs x^{-\beta-\abs\alpha})$ as $x\to\infty$ for each $\abs\alpha\le l$.
\end{list}
\end{lem}

\begin{proof}
The equivalence of (i) and (ii) follows from Lemma \ref{equivofnorms2}.

\begin{claim}[\emph{(iii)$\implies$(ii)}.]
Let $i\in\NZ$. Using \eqref{equivnormsdonspCnt} and the Leibniz rule we have 
\begin{eqnarray*}
\norm{\zeta_i u}_{\donsp\beta{\Cnt l}}
&\le&C\sum_{\abs\alpha\le l}\,\sup_{x\in\don}
\,\abs x^{\beta+\abs\alpha}\,\abs{\dfo_x^\alpha(\zeta_i u)(x)}\\
&\le&C\vtsp\biggl(\,\sum_{\abs\alpha\le l}\,\sup_{x\in\don}
\,\abs x^{\abs\alpha}\,\abs{\dfo_x^\alpha\zeta_i(x)}\biggr)\,
\biggl(\,\sum_{\abs\alpha\le l}\,\sup_{x\in\supp(\zeta_i)}
\,\abs x^{\beta+\abs\alpha}\,\abs{\dfo_x^\alpha u(x)}\biggr)\\
&\le&C\,\norm{\zeta_i}_{\donsp0{\Cnt l}}
\sum_{\abs\alpha\le l}\,\sup_{x\in\supp(\zeta_i)}
\,\abs x^{\beta+\abs\alpha}\,\abs{\dfo_x^\alpha u(x)},
\end{eqnarray*}
where the constants are independent of $i$. With the help of Lemma \ref{lemnoneoflsec} and the fact that $\supp(\zeta_i)\subset\{\abs x>e^{i-1}\}$ we now get (iii)$\implies$(ii).
\end{claim}

\begin{claim}[\emph{(ii)$\implies$(iii)}.]
Suppose (iii) is not satisfied. Thus we can find some multi-index $\alpha$ with $\abs\alpha\le l$ and a sequence of points $\{x_j\}_{j\in\N}$ with $x_j\to\infty$ such that
\begin{equation}
\label{eqniiiantsat}
\abs{x_j}^{\beta+\abs\alpha}\,\abs{\dfo_x^\alpha u(x_j)}
\,\ge\,C\,>\,0
\quad\text{for all $j\in\N$.}
\end{equation}
Choose a sequence $\{j(i)\}_{i\in\N}$ with $j(i)\to\infty$ as $i\to\infty$ and $\abs{x_{j(i)}}>e^i$. Therefore $\zeta_i=1$ on neighbourhood of $x_{j(i)}$ so \eqref{equivnormsdonspCnt} gives
\[
\norm{\zeta_i u}_{\donsp\beta{\Cnt l}}
\ \ge\ C\sum_{\abs\alpha\le l}\,\sup_{x\in\don}
\,\abs x^{\beta+\abs\alpha}\,\abs{\dfo_x^\alpha(\zeta_i u)(x)}
\ \ge\ C\abs{x_{j(i)}}^{\beta+\abs\alpha}\,\abs{\dfo_x^\alpha u(x_{j(i)})}.
\]
The fact that (ii) is not satisfied now follows from \eqref{eqniiiantsat}.
\end{claim}
\end{proof}

By Lemma \ref{denofCinft0incertsp} we know that, for any $\beta\in\R$, $\consp\beta{\Cnt l_0}$ is the separable subspace of $\consp\beta{\Cnt l}$ obtained by taking the closure of $C^\infty_0$ with respect to norm $\norm\cdot_{\consp\beta{\Cnt l}}$. Elements of $\consp\beta{\Cnt l_0}$ can be given an alternative characterisation as follows.

\begin{prop}
\label{altcharofsepsp}
Suppose $l\in\NZ$, $\beta\in\R$ and $u\in\consp\beta{\Cnt l}$. Then $u\in\consp\beta{\Cnt l_0}$ iff $\dfo_x^\alpha u(x)=o(\abs x^{-\beta-\abs\alpha})$ as $x\to\infty$ for each multi-index $\alpha$ with $\abs\alpha\le l$.
\end{prop}

\begin{rem}
\label{usefullocfor1}
It follows that $\Lambda^{-s}\in\consp\beta{\Cnt l_0}$ for any $l\in\NZ$ and $\beta,s\in\R$ with $\beta<s$.
\end{rem}

\begin{proof}
By Lemma \ref{lemntwooflsec} it suffices to show $u\in\consp\beta{\Cnt l_0}$ iff $\norm{\zeta_i u}_{\consp\beta{\Cnt l}}\to0$ as $i\to\infty$. 

\smallskip

Let $u\in\consp\beta{\Cnt l_0}$ and $\epsilon>0$. Thus we can find $u_\epsilon\in C^\infty_0$ with $\norm{u-u_\epsilon}_{\consp\beta{\Cnt l}}\le\epsilon$. Now suppose $i\in\N$ is sufficiently large so that $\supp(u_\epsilon)\subseteq\{\abs x< e^{i-1}\}$. Since $\zeta_i=0$ on the latter set we have $\zeta_iu=\zeta_i(u-u_\epsilon)$. Lemma \ref{lemnoneoflsec} and the continuity of multiplication as a bilinear map $\consp0{\Cnt l}\times\consp\beta{\Cnt l}\to\consp\beta{\Cnt l}$ now imply $\norm{\zeta_i u}_{\consp\beta{\Cnt l}}\le C\epsilon$ for some $C$ which is independent of $\epsilon$. It follows that $\norm{\zeta_i u}_{\consp\beta{\Cnt l}}\to0$ as $i\to\infty$. 

\smallskip

On the other hand, suppose $\norm{\zeta_i u}_{\consp\beta{\Cnt l}}\to0$ as $i\to\infty$. Let $\epsilon>0$ and choose $I\in\N$ so that $\norm{\zeta_I u}_{\consp\beta{\Cnt l}}\le\epsilon$. Now Lemma \ref{equivofnorms1} gives us a constant $C_1$ such that 
\begin{equation}
\label{conseqofeqivn1}
\norm{\eta_{I+1} v}_{\consp\beta{\Cnt l}}
\;\le\;C_1\norm{\eta_{I+1} v}_{\Cnt l}
\end{equation}
for all $v\in\consp\beta{\Cnt l}$. Since $\eta_{I}u\in\Cnt l$ with $\supp(\eta_Iu)\subset\{\abs x< e^{I+1}\}$ we can find $u_\epsilon\in C^\infty_0$ with $\supp(u_\epsilon)\subset\{\abs x< e^{I+1}\}$ and $\norm{\eta_Iu-u_\epsilon}_{\Cnt l}\le\epsilon/C_1$. It follows that $\eta_{I+1}(\eta_I u-u_\epsilon)=\eta_I u-u_\epsilon$ so \eqref{conseqofeqivn1} gives $\norm{\eta_I u-u_\epsilon}_{\consp\beta{\Cnt l}}\le\epsilon$. However $u=\eta_Iu+\zeta_Iu$ so $\norm{u-u_\epsilon}_{\consp\beta{\Cnt l}}\le2\epsilon$, completing the result.
\end{proof}

From Definition \ref{defofpsym} and Proposition \ref{equivnormmodspCnt} it is clear that we have 
\begin{equation}
\label{eqnthatsayswssymlive}
\ssym\gamma\;\subset\;\consp\gamma{\Cnt l}
\end{equation}
for any $\gamma\in\R$ and $l\in\NZ$. Together with Proposition \ref{multibysymprop} we then get the following.

\begin{prop}
\label{contofmultbyssym}
If $E$ is a model space and $\beta,\gamma\in\R$ then multiplication by any $p\in\ssym\gamma$ defines a continuous map $\consp\beta E\to\consp{\beta+\gamma}E$
\end{prop}

\begin{rem}
\label{intoremscomtsec}
Suppose $\gamma\in\R$ and $p\in\ssym\gamma$ with principal part $r^{-\gamma}a(\omega)$. For any $f\in\genXYfns$ define a function $p_f$ by 
\[
p_f(x)\,=\,f(x)p(x)+\bigl(1\negvtsp-\negvtsp f(x)\bigr)\vtsp r^{-\gamma}a(\omega).
\]
If $f=1$ on a neighbourhood of 0 it is clear that $p_f\in\ssym\gamma$ with the same principal part as $p$. On the other hand, if $f=0$ on a neighbourhood of 0 then Lemma \ref{equivofnorms1} and \eqref{equivnormsdonspCnt} imply $p_f$ is contained in $\donsp\gamma{\Cnt l}$ for any $l\in\NZ$.
\end{rem}

Finally, Condition (ii) of Definition \ref{defofpsym} and Lemma \ref{lemntwooflsec} give the following result.

\begin{lem}
\label{lemnthreeoflsec}
Suppose $\gamma\in\R$ and $p\in\ssym\gamma$ with principal part $r^{-\gamma}a(\omega)$. Then, for any $l\in\NZ$,
\[
\lim_{i\to\infty}\,\bignorm{\zeta_i
\bigl(p-a(\omega)r^{-\gamma}\bigr)}_{\donsp\gamma{\Cnt l}}
\;=\;0\;=\;\lim_{i\to\infty}\,\bignorm{\zeta_i
\bigl(p-a(\omega)r^{-\gamma}\bigr)}_{\consp\gamma{\Cnt l}}\vtsp.
\]
\end{lem}

\subsection{Admissible spaces}

\label{admissspace}

\begin{defn}
For any $m\in\R$ let $\psym{m}$ denote the set of functions on $a\in C^\infty_\loc(\R^n\!\times\!\R^n)$ which satisfy estimates of the form
\begin{equation}
\label{normestforpsyms}
\bigabs{\dfo^\alpha_x\dfo^{\alpha'}_\xi\vtsp a(x,\xi)}
\;\le\;C_{\alpha,\alpha'}\Lambda^{m-\abs{\alpha'}}(\xi)
\end{equation}
for all multi-indices $\alpha$ and $\alpha'$. The best constants in \eqref{normestforpsyms} provide $\psym{m}$ with a collection of semi-norms making it into a locally convex space. 

For any $a\in\psym{m}$ we shall use $a(x,\dfo_x)$ to denote the pseudo-differential operator defined by the symbol $a(x,\xi)$. The set of all pseudo-differential operators of order $m$ (i.e.~the set of all operators defined by symbols in $\psym{m}$) shall be denoted by $\pdo{m}$.
\end{defn}

\begin{rem}
\label{remaboutconv}
For any $m\in\R$ the pairing $(a,u)\mapsto a(x,\dfo_x)u$ defines a continuous bilinear map $\psym m\times\mathscr S\to\mathscr S$ and a bilinear map $\psym m\times\mathscr S'\to\mathscr S'$ which is separately continuous in each variable. If $m,l\in\R$ it can also be shown that the composition of operators in $\pdo m$ and $\pdo l$ gives an operator in $\pdo{m+l}$, whilst the adjoint of an operator in $\pdo m$ is again in $\pdo m$. Furthermore the corresponding symbol maps $\psym m\times\psym l\to\psym{l+m}$ and $\psym m\to\psym m$ are continuous bilinear and continuous anti-linear respectively. Further details can be found in Section 18.1 of \cite{Hor3} for example.
\end{rem}

If $a\in\psym{m}$ for some $m\in\R$ then we can write
\begin{equation}
\label{intformofpdo}
a(x,\dfo_x)u
\;=\;\int K_a(x,y)\vtsp u(y)\d^ny
\end{equation}
for all $u\in\mathscr S'$, where $K_a\in\mathscr S'(\R^n\times\R^n)$ is the Schwartz kernel of $a(x,\dfo_x)$; that is, $K_a(x,y)=(2\pi)^{-n}\widehat a(x,y-x)$ where $\widehat a$ is the Fourier transform of $a(x,\xi)$ with respect to the second variable (see Section 18.1 of \cite{Hor3} for example).

\begin{lem}
\label{smoothingmapslem}
We have $K_a\in\mathscr S(\R^n\times\R^n)$ iff $a\in\mathscr S(\R^n\times\R^n)$. Furthermore, in this case $a(x,\dfo_x)$ defines a continuous map $\mathscr S'\to\mathscr S$. 
\end{lem}

\begin{proof}
The first part of the result follows from the fact that the Fourier transform defines an isomorphism $\mathscr S\to\mathscr S$. Now suppose $K_a\in\mathscr S(\R^n\times\R^n)$ and choose any multi-index $\alpha$ and $s\in\R$. Then \eqref{intformofpdo} gives
\begin{equation}
\label{estforcontwereq}
\bigabs{\Lambda^s(x)\dfo_x^\alpha\vtsp a(x,\dfo_x)u}
\;\le\;\bigabs{\bigl(\Lambda^s(x)\dfo_x^\alpha K_a(x,\cdot),\,u\bigr)_{\R^n}}
\end{equation}
for all $u\in\mathscr S'$. The assumption that $K_a\in\mathscr S(\R^n\times\R^n)$ immediately implies $\Lambda^s(x)\dfo_x^\alpha K_a(x,\cdot)$ forms a bounded subset of $\mathscr S$ as $x$ varies over $\R^n$. The second part of the result now follows from \eqref{estforcontwereq} and the continuity of the dual pairing of $\mathscr S$ and $\mathscr S'$. 
\end{proof}

As a corollary of this result we have the following.

\begin{cor}
\label{corinformforapptov}
Let $s\in\R$ and $\phi_1,\phi_2\in C^\infty$ with $\supp(\phi_1)\cap\supp(\phi_2)=\emptyset$ and either $\phi_1\in C^\infty_0$ or $\phi_2\in C^\infty_0$. Then the pseudo-differential operator $\phi_1(x)\Lambda^s(\dfo_x)\phi_2(x)$ defines a continuous map $\mathscr S'\to\mathscr S$.
\end{cor}

\begin{proof}
The Schwartz kernel of the operator $\phi_1(x)\Lambda^s(\dfo_x)\phi_2(x)$ is 
\[
K(x,y)\;=\;(2\pi)^{-n}\phi_1(x)\psi(y-x)\phi_2(y)
\]
where $\psi\in\mathscr S'$ is the Fourier transform of $\Lambda^s$. By standard properties of the Fourier transform of symbols (see Proposition VI.4.1 in \cite{St} for example) $\psi$ is smooth and rapidly decaying away from 0, along with all its derivatives. On the other hand $\phi_1$ and $\phi_2$ have disjoint supports, at least one of which is compact. It follows that $K\in\mathscr S(\R^n\times\R^n)$. Lemma \ref{smoothingmapslem} now completes the result.
\end{proof}

\medskip

Using the mapping properties of pseudo-differential operators we can now single out a special class of model spaces which will provide the natural function space setting for our main results.

\begin{defn}
An \emph{admissible space} is a model space $E$ satisfying the following additional condition.
\begin{list}{}{\listinit}
\item[(B)\hf]
We have a continuous bilinear map $\psym{0}\times E\to E$ which sends $(a,u)$ to $a(x,\dfo_x)u$.
\end{list}
\end{defn}

\begin{rem}
Condition (A2) for a model space is a special case of Condition (B).
\end{rem}

\begin{ex}
A rich class of admissible spaces is provided by the Besov spaces $B^s_{pq}$ for $s\in\R$ and $p,q\in[1,\infty]$, and the Triebel-Lizorkin spaces $F^s_{pq}$ for $s\in\R$ and $p,q\in[1,\infty]$ with $q\ne1$ if $p=\infty$; see \cite{Tr1} and \cite{Tr2} for the definitions of these spaces and the justification Conditions (A1) to (A4) and (B)\footnote{In the cited literature the continuity of the bilinear map in Condition (B) is only explicitly established for the second variable. However the full continuity of this map can be easily obtained from the proof of the relevant result (Theorem 6.2.2 in \cite{Tr2}) by using symbol norms of $a$ to make simple estimates of the constants appearing therein.}.

The class of Besov spaces and Triebel-Lizorkin spaces include a large number of the `standard' function spaces as follows (see \cite{Tr1}).
\begin{description}
\item
$\quad F^s_{p2}=\Sob{p}{s}$ (the Sobolev or Bessel-potential spaces) for $s\in\R$ and $p\in(1,\infty)$.
\item
$\quad F^0_{p2}=h_p$ (the local Hardy spaces) for $p\in[1,\infty)$.
\item
$\quad F^0_{\infty2}=bmo$ (the inhomogeneous version of BMO).
\item
$\quad B^s_{\infty\infty}=\Zyg{s}$ (the Zygmund spaces) for $s>0$.
\item
$\quad B^s_{pp}=F^s_{pp}=W^s_p$ (the Slobodeckij spaces) for $s\in\R^+\!\setminus\!\N$ and $p\in[1,\infty)$.
\item
$\quad B^s_{pq}=\Lambda^s_{p,q}$ (the Lipschitz spaces) for $s>0$, $p\in[1,\infty)$ and $q\in[1,\infty]$.
\end{description}
The most notable omissions from this list are the spaces $\Cnt l$ for $l\in\NZ$; these spaces do not satisfy Condition (B) (n.b.\ although the Zygmund space $\Zyg{s}$ coincides with the H\"older space $\Hol l\sigma$ whenever $l\in\NZ$, $\sigma\in(0,1)$ and $s=l+\sigma$, we only have a \emph{strict} inclusion $\Zyg l\subset\Cnt l$ when $l\in\N$).
\end{ex}

\begin{defn}
\label{defnofElver2}
Suppose $E$ is an admissible space and let $s\in\R$. We define $E^s$ to be the set of $u\in\mathscr S'$ for which $\Lambda^s(\dfo_x)u\in E$. Furthermore we give this set a norm $\norm\cdot_{E^s}$ defined by $\norm{u}_{E^s}=\norm{\Lambda^s(\dfo_x)u}_E$.
\end{defn}

\begin{rem}
If $E=\Sob{p}{s}$ for some $s\in\R$ and $p\in(1,\infty)$ then $E^{s'}=\Sob{p}{s+s'}$ for any $s'\in\R$. Likewise, if $E=\Zyg{s}$ for some $s>0$ then $E^{s'}=\Zyg{s+s'}$ for any $s'\in\R$, provided $s+s'>0$. 
\end{rem}

\begin{prop}
\label{mappingcoroforadmissspaces}
Suppose $s,m\in\R$. Then the assignment $(a,u)\mapsto a(x,\dfo_x)u$ defines a continuous bilinear map $\psym{m}\times E^s\to E^{s-m}$. In particular, given any $s\in\R$ and multi-index $\alpha$, the differential operator $\dfo_x^\alpha$ defines a continuous map $E^{s+\abs\alpha}\to E^s$. 
\end{prop}

\begin{proof}
Given $a\in\psym{m}$ define a new symbol $b\in\psym{0}$ as the symbol of the pseudo-differential operator $b(x,\dfo_x)=\Lambda^{s-m}(\dfo_x)a(x,\dfo_x)\Lambda^{-s}(\dfo_x)$. Standard results on the calculus of pseudo-differential operators (see Remark \ref{remaboutconv}) imply that the assignment $a\mapsto b$ defines a continuous map $\psym{m}\to\psym{0}$. On the other hand $\Lambda^{s-m}(\dfo_x)a(x,\dfo_x) u=b(x,\dfo_x)\Lambda^s(\dfo_x) u$ for any $u\in\mathscr S'$. The result now follows from Definition \ref{defnofElver2} and Condition (B) on the admissible space $E$. 
\end{proof}

\begin{prop}
\label{basicpropadmisssp}
Any admissible space(s) $E$ (and $F$) satisfy the following properties.
\begin{list}{}{\listinit}
\item[\emph{(i)}\hf]
If $s\in\R$ then $E^s$ is again an admissible space.
\item[\emph{(ii)}\hf]
If $s,s'\in\R$ then $(E^s)^{s'}=E^{s+s'}$.
\item[\emph{(iii)}\hf]
If $s=l\in\NZ$ then $E^s=E^l$ up to equivalent norms (where $E^l$ is as given by Definition \ref{defnofElver1}).
\item[\emph{(iv)}\hf]
We have a continuous inclusion $E\hookrightarrow F$ iff we have a continuous inclusion $E^s\hookrightarrow F^s$ for any $s\in\R$.
\item[\emph{(v)}\hf]
We have a local inclusion $E_\loc\hookrightarrow F_\loc$ iff we have a local inclusion $E^s_\loc\hookrightarrow F^s_\loc$ for any $s\in\R$.
\item[\emph{(vi)}\hf]
The spaces $E_0$ and its dual $E_0^*$ are admissible spaces. 
\item[\emph{(vii)}\hf]
If $s\in\R$ then $(E^s)_0=(E_0)^s$ and $(E^s_0)^*=(E^*_0)^{-s}$.
\end{list}
\end{prop}

\begin{proof}
\begin{claim}[Part (i).]
The only non trivial Conditions are (B) and (A4) (n.b.~Condition (A2) is covered by Condition (B)). The former is established in Proposition \ref{mappingcoroforadmissspaces} whilst the latter uses the fact that conjugation of a pseudo-differential operator by a diffeomorphism which is linear outside a compact region gives another pseudo-differential operator of the same order (see Theorem 18.1.17 in \cite{Hor3} for example\footnote{Technically this only gives a local version of what we need; however it is easy to see how the proof can be modified to give the conjugation result as stated.}).
\end{claim}

\begin{claim}[Part (ii).]
This is an easy consequence of the identity $\Lambda^s(\dfo_x)\Lambda^{s'}(\dfo_x)=\Lambda^{s+s'}(\dfo_x)$.
\end{claim}

\begin{claim}[Part (iii).]
Suppose $u\in E^s$ and let $\alpha$ be a multi-index with $\abs\alpha\le l=s$. Thus $\Lambda^s(\dfo_x)u\in E$ whilst $\xi^\alpha\Lambda^{-s}(\xi)$ defines a symbol in $\psym{0}$ so $\dfo_x^\alpha\Lambda^{-s}(\dfo_x)$ defines a continuous map on $E$ (by Condition (B) for the admissible space $E$). Hence $\dfo_x^\alpha u=\dfo_x^\alpha\Lambda^{-s}(\dfo_x)\Lambda^s(\dfo_x)u\in E$ and we have a norm estimate of the form
\[
\norm{\dfo_x^\alpha u}_E
\;\le\; C\norm{\Lambda^s(\dfo_x)u}_E
\;=\; C\norm{u}_{E^s}\,.
\]
The existence of a continuous inclusion $E^s\hookrightarrow E^l$ now follows from the definition of $E^l$. 

\smallskip

On the other hand we can write 
\begin{equation}
\label{redformofLam2l}
\Lambda^{2s}(\dfo_x)
\;=\;\Lambda^{2l}(\dfo_x)
\;=\;\Bigl(1+\sum_{i=1}^n\dfo_i^2\Bigr)^l
\;=\;\sum_{\abs\alpha\le l} a^\alpha \dfo_x^\alpha \dfo_x^\alpha
\end{equation}
for some constants $a^\alpha\in\C$. Now suppose $u\in E^l$. Therefore $\dfo_x^\alpha u\in E$ with $\norm{\dfo_x^\alpha u}_E\le C_{1,\alpha}\norm{u}_{E^l}$ (see Remark \ref{remnofElver1}) whilst $\Lambda^{-s}(\xi)\xi^\alpha\in\psym{0}$ so $\Lambda^{-s}(\dfo_x)\dfo_x^\alpha$ defines a continuous map on $E$. Hence $\Lambda^{-s}(\dfo_x)\dfo_x^\alpha\dfo_x^\alpha u\in E$ and 
\[
\norm{\Lambda^{-s}(\dfo_x)\dfo_x^\alpha\dfo_x^\alpha u}_E
\;\le\;C_{2,\alpha}\norm{u}_{E^l}. 
\]
With the help of \eqref{redformofLam2l} if follows that 
\[
\Lambda^s(\dfo_x)u
\;=\;\sum_{\abs\alpha\le l} a^\alpha \Lambda^{-s}(\dfo_x)\dfo_x^\alpha \dfo_x^\alpha u
\;\in\;E
\]
and $\norm{u}_{E^s}=\norm{\Lambda^s(\dfo_x)u}_E\le C\norm{u}_{E^l}$. This completes the proof of part (iii).
\end{claim}

\begin{claim}[Part (iv).]
The fact that a continuous inclusion $E\hookrightarrow F$ induces a continuous inclusion $E^s\hookrightarrow F^s$ is trivial. Part (ii) now gives the converse. 
\end{claim}

\begin{claim}[Part (v).]
Let $\phi\in C^\infty_0$ and choose $\phi_1,\phi_2\in C^\infty_0$ with $\phi_1\succ\phi$ and $\phi_2=1$ on a neighbourhood of $\supp(\phi_1)$. Setting $\phi_3=1-\phi_2$ it follows that $\phi_3\in C^\infty$ and $\supp(\phi_1)\cap\supp(\phi_3)=\emptyset$. Now let $u\in E^s_\loc$. Therefore $\phi u\in E^s$ or, equivalently, $\Lambda^s(\dfo_x)\phi u\in E$. Using the local inclusion $E_\loc\hookrightarrow F_\loc$ it follows that $\phi_2\Lambda^s(\dfo_x)\phi u\in F$ and 
\begin{equation}
\label{normestpart1locincres}
\bignorm{\phi_2\Lambda^s(\dfo_x)\phi u}_F
\;\le\;C_1\norm{\phi u}_{E^s},
\end{equation}
where $C_1$ may depend on $\phi$ (through $\phi_2$) but not on $u$. On the other hand Corollary \ref{corinformforapptov} implies $\phi_3\Lambda^s(\dfo_x)\phi_1$ defines a continuous map $\mathscr S'\to\mathscr S$. Condition (A1) for the model spaces $E^s$ and $F$ then shows that $\phi_3\Lambda^s(\dfo_x)\phi_1$ defines a continuous map $E^s\to F$. Therefore $\phi_3\Lambda^s(\dfo_x)\phi u=\phi_3\Lambda^s(\dfo_x)\phi_1\vtsp \phi u\in F$ and 
\begin{equation}
\label{normestpart2locincres}
\bignorm{\phi_3\Lambda^s(\dfo_x)\phi u}_F
\;\le\;C_2\norm{\phi u}_{E^s},
\end{equation}
where $C_1$ may depend on $\phi$ (through $\phi_1$ and $\phi_3$) but not on $u$. Since $\phi_2+\phi_3=1$, \eqref{normestpart1locincres} and \eqref{normestpart2locincres} combine to establish the existence of a local inclusion $E^s_\loc\hookrightarrow F^s_\loc$. Part (ii) now gives the converse.
\end{claim}

\begin{claim}[Part (vi).]
Condition (B) for $E_0$ follows from the same Condition for $E$ and the fact that pseudo-differential operators preserve the set $\mathscr S$ (which is dense in $E_0$). For the dual space we can use Lemma \ref{EmodimpliesE0stmod} and the fact that the map which sends a pseudo-differential operator $a(x,\dfo_x)\in\pdo{0}$ to its adjoint induces an (anti-linear) isomorphism on $\psym{0}$ (see Remark \ref{remaboutconv}).
\end{claim}

\begin{claim}[Part (vii).]
Clearly $(E_0)^s\subseteq E^s$ whilst $\Lambda^s(\dfo_x):(E_0)^s\to E_0$ is an isomorphism which preserves the set $\mathscr S$. The fact that $\mathscr S$ is dense in $E_0$ now implies the same is true for $(E_0)^s$, with the first identity following immediately. The second identity can be obtained from \eqref{defofdualnorm} and the expression $(u,v)_{\R^n}=\bigl(\Lambda^s(\dfo_x)u,\vtsp\Lambda^{-s}(\dfo_x)v\bigr)_{\R^n}$ which is valid for all $u\in\mathscr S$ and $v\in\mathscr S'$. 
\end{claim}
\end{proof}

\bigskip

For the proofs of the next three results let $\psi$ denote a choice of function in $C^\infty_0$ with $\psi=1$ in a neighbourhood of 0 and $\Ran\psi=[0,1]$. Also, for each $j\in\N$, define $\psi_j\in C^\infty_0$ by $\psi_j(x)=\psi(x/j)$. 

\begin{lem}
\label{approxdenofEinEdel}
Suppose $f\in E$ for some admissible space $E$. Then we can find a sequence $\{f_j\}_{j\in\N}\subset\bigcap_{s\in\R} E^s$ with $f_j\to f$ in $E^{-\delta}$ for any $\delta>0$.
\end{lem}

\begin{proof}
For each $j\in\N$ set $f_j=\psi_j(\dfo_x)f$. Since $\psi_j(\xi)\in\psym m$ for any $m\in\R$, Proposition \ref{mappingcoroforadmissspaces} gives $f_j\in E^s$ for any $s\in\R$. On the other hand, a straightforward check shows $\Lambda^{-\delta}(\xi)\psi_j(\xi)\to\Lambda^{-\delta}(\xi)$ in $\psym 0$ for any $\delta>0$. Condition (B) for the admissible space $E$ then implies $\Lambda^{-\delta}(\dfo_x)f_j\to\Lambda^{-\delta}(\dfo_x)f$ in $E$ for any $\delta>0$; by definition this means $f_j\to f$ in $E^{-\delta}$.
\end{proof}

Although the space $\consp\beta E$ need not contain $C^\infty_0$ as a dense subset for a general admissible space $E$, Lemma \ref{approxdenofEinEdel} leads to the following slightly weaker result.

\begin{lem}
\label{bestapproxresinnonsepcase}
Suppose $f\in\consp\beta E$ for some admissible space $E$ and $\beta\in\R$. Then we can find a sequence $\{f_i\}_{i\in\N}\subset C^\infty_0$ such that $f_i\to f$ in $\consp{\beta-\epsilon}{E^{-\delta}}$ for any $\epsilon,\delta>0$.
\end{lem}

\begin{proof}
\begin{claim}[Claim: Given $f\in\clysp\beta E$ we can find a sequence $\{f_j\}_{j\in\N}\subset\bigcap_{s\in\R}\clysp\beta{E^s}$ with $f_j\to f$ in $\clysp\beta{E^{-\delta}}$ for any $\delta>0$.]
Since multiplication by $e^{\beta t}$ defines an isomorphism $\clysp\beta{E^s}\to\clysp0{E^s}=E^s(\cly)$ for any $s\in\R$, we may prove the claim assuming $\beta=0$. Consider the notation of Remark \ref{enlargechirem} and define $\widetilde{\chi}'_i=\chi_i\circ\Psi_i^{-1}\in C^\infty$ for each $i\in I$. Now Lemma \ref{cordequivnorm} implies $g_i:=(\Psi_i^{-1})^*(\chi_if)\in E$ so Lemma \ref{approxdenofEinEdel} gives us a sequence $\{g_{ij}\}_{j\in\N}\subset\bigcap_{s\in\R}E^s$ with $g_{ij}\to g_i$ in $E^{-\delta}$ for any $\delta>0$. Condition (A2) and the fact that $\chi_i'\succ\chi_i$ then give $\widetilde{\chi}'_ig_{ij}\to\widetilde{\chi}'_ig_i=g_i$ in $E^{-\delta}$ for any $\delta>0$. Setting 
\[
f_j\;=\;\sum_{i\in I}\Psi_i^*(\widetilde{\chi}_i'g_{ij})
\]
for any $j\in\N$, Lemma \ref{cordequivnorm} implies $\{f_j\}_{j\in\N}\subset\bigcap_{s\in\R}E^s(\cly)$ whilst
\[
f_j\;\longrightarrow\;\sum_{i\in I}\Psi_i^*g_i
\;=\;\sum_{i\in I}\chi_i f
\;=\;f
\]
in $E^{-\delta}(\cly)$ for any $\delta>0$. This completes the claim.
\end{claim}

\medskip

Let $f\in\consp\beta E$. By Lemma \ref{equivofnorms1} we have $\eta_0f\in E$ and $\zeta_0f\in\donsp\beta E$ so Lemma \ref{approxdenofEinEdel} and the above Claim (coupled with \eqref{Thetaisaniso}) give us sequences $\{g_j\}_{j\in\N}\subset\bigcap_{s\in\R} E^s$ and $\{h_j\}_{j\in\N}\subset\bigcap_{s\in\R}\donsp\beta{E^s}$ with $g_j\to\eta_0 f$ in $E^{-\delta}$ and $h_j\to\zeta_0 f$ in $\donsp\beta{E^{-\delta}}$ for any $\delta>0$. For any $j\in\N$ define further functions by $f'_j=\eta_1 g_j+\zeta_{-1} h_j$ and $f_j=\psi_jf'_j$. Another application of Lemma \ref{equivofnorms1} now shows $\{f'_j\}_{j\in\N}\subset\bigcap_{s\in\R}\consp\beta{E^s}$ and $f'_j\to(\eta_1\eta_0+\zeta_{-1}\zeta_0)f=f$ in $\consp\beta{E^{-\delta}}$ for any $\delta>0$. 

By Remark \ref{easyincrem} and Lemma \ref{ElintoCntglob} we can find $k\in\NZ$ so that we have a continuous inclusion $\consp\beta{E^{k+l}}\hookrightarrow\consp\beta{\Cnt l}$ for all $l\in\NZ$. It follows that $\{f_j'\}_{j\in\N}\subset\bigcap_{l\in\N}\consp\beta{\Cnt l}\subset C^\infty_\loc$. Therefore $\{f_j\}_{j\in\N}\subset C^\infty_0$. Now let $\epsilon,\delta>0$. A straightforward application of Proposition \ref{equivnormmodspCnt} shows $\psi_j\to1$ in $\consp{-\epsilon}{C^l}$ for any $l\in\NZ$. Proposition \ref{multibysymprop} and the convergence $f'_j\to f$ in $\consp\beta{E^{-\delta}}$ then implies $f_j\to f$ in $\consp{\beta-\epsilon}{E^{-\delta}}$.
\end{proof}

\begin{lem}
\label{denseofSinSp}
The set $C^\infty_0$ is dense in both $\mathscr S$ and $\mathscr S'$.
\end{lem}

\begin{proof}
For the density of $C^\infty_0$ in $\mathscr S$ see Proposition VI.1.3 in \cite{Y}. Now let $f\in\mathscr S'$ and, for each $j\in\N$, set $f_j=\psi_j(x)\psi_j(\dfo_x)f$. Thus $f_j$ has compact support (contained in $\supp(\psi_j)$) whilst the symbol $\psi_j(x)\psi_j(\xi)$ is contained in $C^\infty_0(\R^n\!\times\!\R^n)$ so Lemma \ref{smoothingmapslem} gives $f_j\in\mathscr S$. Therefore $\{f_j\}_{j\in\N}\subset C^\infty_0$. 

A straightforward check shows $\Lambda^{-1}(x)\psi_j(x)\psi_j(\dfo_x)\to\Lambda^{-1}(x)$ in $\psym 1$ so Remark \ref{remaboutconv} gives $\Lambda^{-1}(x)f_j\to\Lambda^{-1}(x)f$ in $\mathscr S'$. Since the operator (of multiplication by) $\Lambda(x)$ defines a continuous map on $\mathscr S'$ we now get $f_j\to f$ in $\mathscr S'$. 
\end{proof}

The next result follows from Lemma \ref{denseofSinSp} and an argument similar to that used to prove the Claim in the Proof of Lemma \ref{bestapproxresinnonsepcase}.

\begin{lem}
\label{denseofZSinZSp}
For any $\beta\in\R$ the set $\clysp\beta{\mathscr S}$ is dense in $\clysp\beta{\mathscr S'}$. 
\end{lem}

The compactness results given in Section \ref{subsecBasicFunSpProp} can be refined for admissible spaces. We begin with two technical lemmas, the first of which is essentially the Ascoli-Arzel\`a Theorem (see Section III.3 of \cite{Y}).

\begin{lem}
\label{EssAscoliArzela}
Any bounded subset of $\mathscr S$ is pre-compact.
\end{lem}

\begin{lem}
\label{uselemforcptadmisssp}
If $E$ is an admissible space and $s>0$ then $E^s$ is locally compact in $E$. 
\end{lem}

\begin{proof}
Let $\phi\in C^\infty_0$ and choose a sequence $\{\psi_i\}_{i\in\N}\subset C^\infty_0$ such that $\psi_i\to\Lambda^{-s}$ in $\consp0{C^\infty}$ (which is possible by Remark \ref{usefullocfor1} since $s>0$). It follows that the sequence of symbols $\phi(x)\psi_i(\xi)$ converges to $\phi(x)\Lambda^{-s}(\xi)$ in $\psym{0}$ and so the pseudo-differential operator $\phi(x)\psi_i(\dfo_x)$ converges to $\phi(x)\Lambda^{-s}(\dfo_x)$ in $\mathscr L(E,E)$ (by Condition (B) for the admissible space $E$). However the map $E^s\to E$ given by multiplication by $\phi$ can be written as $\phi(x)\Lambda^{-s}(\dfo_x)\vtsp\Lambda^s(\dfo_x)$ where $\Lambda^s(\dfo_x)$ acts as an isomorphism $E^s\to E$. Using the fact that the set of compact operators is closed in $\mathscr L(E,E)$ it therefore suffices to show that the operator $\phi(x)\psi_i(\xi)$ defines a compact map $E\to E$ for any $i\in\N$.

Now $\phi(x)\psi_i(\xi)\in C^\infty_0(\R^n\times\R^n)\subset\mathscr S(\R^n\times\R^n)$ so Lemma \ref{smoothingmapslem} implies $\phi(x)\psi_i(\dfo_x)$ defines a continuous map $\mathscr S'\to\mathscr S$. On the other hand Condition (A1) for the admissible space $E$ gives us continuous inclusions $E\hookrightarrow\mathscr S'$ and $\mathscr S\hookrightarrow E$. By composing these maps and using Lemma \ref{EssAscoliArzela} it follows that $\phi(x)\psi_i(\dfo_x)$ defines a compact map $E\to E$. 
\end{proof}

Lemma \ref{uselemforcptadmisssp} and Proposition \ref{genformofcpt} (with $E$, $F$ and $G$ replaced by $E^s$, $F$ and $E$ respectively) immediately leads to the following useful result.

\begin{prop}
\label{genuseformofcpt}
Let $\beta,\gamma,s\in\R$ with $s>0$ and suppose, for an admissible space $E$ and a model space $F$, multiplication defines a continuous bilinear map $E\times F\to E$ (or, more generally, $E^s\times F\to E$). Then multiplication by any $\phi\in\consp\gamma{F_0}$ defines a compact map $\consp\beta{E^s}\to\consp{\beta+\gamma}E$. 
\end{prop}

\begin{cor}
\label{specformofcpt}
Suppose $E$ be an admissible space and $\beta,\gamma,s\in\R$ with $\gamma>\beta$ and $s>0$. Then the inclusion $\consp\gamma{E^s}\hookrightarrow\consp\beta E$ is compact. 
\end{cor}

\begin{proof}
By Lemma \ref{lemaboutmultibyCl} we know that multiplication defines a continuous map $E\times\Cnt l\to E$ for all sufficiently large $l\in\NZ$. Furthermore $\beta-\gamma<0$ so $1\in\consp{\beta-\gamma}{\Cnt l_0}$ by Remark \ref{usefullocfor1}. The result now follows directly from an application of Proposition \ref{genuseformofcpt}.
\end{proof}

\medskip

For any admissible space $E$, $\beta\in\R$, $k\in\N$ and vector of non-negative integers $\kappa=(\kappa_1,\dots,\kappa_k)$, we define the following product spaces;
\begin{gather*}
E^\kappa\,=\,\prod_{i=1}^k E^{\kappa_i},
\quad
\clyvsp\beta\kappa E\,=\,\prod_{i=1}^k\clysp\beta{E^{\kappa_i}},\\
\donvsp\beta\kappa E\,=\,\prod_{i=1}^k\donsp{\beta-\kappa_i}{E^{\kappa_i}}
\quad\text{and}\quad
\convsp\beta\kappa E\,=\,\prod_{i=1}^k\consp{\beta-\kappa_i}{E^{\kappa_i}}.
\end{gather*}
Let $r^{\pm\kappa}$ denote the matrix operator given as multiplication by $\diag(r^{\pm\kappa_1},\dots,r^{\pm\kappa_k})$. Thus \eqref{Thetaisaniso} and Proposition \ref{basicchangbetaiso} imply that we have an isomorphism
\begin{equation}
\label{icbtdahal1}
\Theta^*\vtsp r^{-\kappa}
\;:\;\donvsp\beta\kappa E
\;\xrightarrow{r^{-\kappa}}\;\prod_{i=1}^k\donsp\beta{E^{\kappa_i}}
\;\xrightarrow{\Theta^*}\;\clyvsp\beta\kappa E
\end{equation}
with inverse
\begin{equation}
\label{icbtdahal2}
r^{\kappa}\vtsp(\Theta^{-1})^*
\;:\;\clyvsp\beta\kappa E
\;\xrightarrow{(\Theta^{-1})^*}\;\prod_{i=1}^k\donsp\beta{E^{\kappa_i}}
\;\xrightarrow{r^{\kappa}}\;\donvsp\beta\kappa E.
\end{equation}
Finally let $\pi_\kappa$ denote the projection defined on any of the above product spaces by
\[
(\pi_\kappa u)_i
\,=\,\begin{cases}
u_i&\text{if $\kappa_i>0$}\\
0&\text{if $\kappa_i=0$.}
\end{cases}
\]
Clearly $\pi_\kappa$ commutes with multiplication by a scalar or a diagonal matrix.

\subsection{Elliptic operators}

Unless otherwise stated we shall assume that the coefficients of any differential operator on $\R^n$ are contained in $C^\infty$. In particular, Condition (A2) and Remark \ref{remnofElver1} imply that any differential operator $A(x,\dfo_x)$ on $\R^n$ of order $m$ defines a continuous map $E^m\to E$ for any model space $E$. 

\begin{defn}
\label{defnofellingenset}
Let $\mathcal A(x,\dfo_x)$ be a $k\times k$ system of differential operators on $\R^n$ with entries $A_{ij}(x,\dfo_x)$ for $i,j=1,\dots,k$. We say that $\mathcal A$ is \emph{elliptic} on some open set $U\subseteq\R^n$ if the following conditions are satisfied.
\begin{list}{}{\listinit}
\item[(i)\hf]
There exist vectors $\mu=(\mu_1,\dots,\mu_k)$ and $\nu=(\nu_1,\dots,\nu_k)$ of non-negative integers such that $\order A_{ij}=\mu_j-\nu_i$ (with $A_{ij}=0$ whenever $\mu_j-\nu_i<0$). Furthermore $\min_i\nu_i=0$. 
\item[(ii)\hf]
Let $a_{ij}(x,\xi)$ denote the principal symbol of $A_{ij}(x,\xi)$ (so $a_{ij}(x,\xi)$ is a homogeneous polynomial of degree $\mu_j-\nu_i$ in $\xi$) and define $\detsym{\mathcal A}(x,\xi)$ to be the determinant of the $k\times k$ matrix with entries $a_{ij}(x,\xi)$. Then $\detsym{\mathcal A}(x,\xi)\ne0$ for any $(x,\xi)\in U\times\don$.
\end{list}
We say that $\mathcal A$ is \emph{uniformly elliptic} on $U$ if Condition (ii) can be replaced by the following stronger condition.
\begin{list}{}{\listinit}
\item[(ii)$'$\hf]
Set $\sumord{\mu\negvtsp}{\negvtsp\nu}=\sum_{i=1}^k(\mu_i\negvtsp-\negvtsp\nu_i)$ and let $\detsym{\mathcal A}(x,\xi)$ be defined as in Condition~(ii) above. Then $\abs{\detsym{\mathcal A}(x,\xi)}\ge C\abs\xi^{\sumord\mu\nu}$ for all $(x,\xi)\in U\times\R^n$.
\end{list}
We shall refer to the pair $(\mu,\nu)$ as the \emph{order} of $\mathcal A$.
\end{defn}

In the case of scalar operators (i.e.~when $k=1$) this definition of ellipticity clearly reduces to the usual one. Furthermore $\nu_1=0$ and $\mu_1=m$, the usual order of the operator. On the other hand, if $\mathcal A$ is elliptic on some open set $U$ then $\mathcal A$ is uniformly elliptic on any bounded open set $V$ with $\Closure(V)\subseteq U$.

\begin{rem}
If $n\ge3$ ellipticity implies that $\sumord{\mu\negvtsp}{\negvtsp\nu}$ is even (say $2l$) and the polynomial equation $\detsym{\mathcal A}(x,\xi+t\eta)$, $\xi,\eta\in\don$, has exactly $l$ roots with $\im t>0$ (and hence exactly $l$ roots with $\im t<0$). These conditions need not be satisfied when $n=2$ and must be imposed as extra assumptions when working with boundary value problems. However we do \emph{not} need to impose such assumptions here.
\end{rem}

\begin{rem}
\label{rem267plusabit}
Suppose $\mathcal A$ and $\mathcal B$ are $k\times k$ systems of differential operators of order $(\mu,\nu)$ on $U\subseteq\R^n$. Furthermore suppose the coefficients of $\mathcal B$ are bounded (pointwise on $U$) by some function $b(x)$. Then we have
\[
\bigabs{\detsym{\mathcal A+\mathcal B}(x,\xi)-\detsym{\mathcal A}(x,\xi)}
\;\le\;C\bigl(b(x)+b^k(x)\bigr)\vtsp\abs{\xi}^{\sumord\mu\nu}
\]
for all $(x,\xi)\in U\times\R^n$, where $C$ may depend on $\mathcal A$ but not on $\mathcal B$, $b$ or $(x,\xi)$. It follows that if $\mathcal A$ is uniformly elliptic on $U$ and $b$ is bounded by a sufficiently small constant (depending on $\mathcal A$) then $\mathcal A+\mathcal B$ is also uniformly elliptic on $U$.
\end{rem}

Suppose $\mathcal A$ is an operator of order $(\mu,\nu)$ which is uniformly elliptic on an open set $U\subseteq\R^n$. Clearly $\mathcal A$ defines a continuous map $E^\mu\to E^\nu$ for any admissible space $E$. Condition (B) for admissible spaces and the calculus of pseudo-differential operators also allows us to derive the following regularity result and elliptic estimates.

\begin{thm}
\label{stanapriorireg}
Suppose $\mathcal A$ is a $k\times k$ system of differential operators on $\R^n$ of order $(\mu,\nu)$ which is uniformly elliptic on some open set $U\subseteq\R^n$. Suppose further that $\chi_1,\chi_2\in C^\infty$ with $\chi_1\prec\chi_2$ and $\supp(\chi_1)\subseteq U$. If $u\in\mathscr S'$ satisfies $\chi_2\mathcal A u\in E^\nu$ and $\pi_\mu\chi_2 u\in(E^{-1})^\mu$ for some admissible space $E$, then we also have $\chi_1 u\in E^\mu$. Furthermore
\[
\norm{\chi_1 u}_{E^\mu}
\;\le\;C\bigl(\norm{\chi_2\mathcal A u}_{E^\nu}+\norm{\pi_\mu \chi_2 u}_{(E^{-1})^\mu}\bigr)
\]
for any such $u$.
\end{thm}

\begin{proof}
Without loss of generality we may assume $\chi_2=1$ on a neighbourhood of $\supp(\chi_1)$ (if this were not the case we could simply replace $\chi_2$ with $\chi_3\chi_2$ where $\chi_3\in C^\infty$ is chosen so that $\chi_3\chi_2=1$ on $\chi_2^{-1}(1/2,+\infty)\,$). It follows that we can find $\chi\in C^\infty$ with $\chi_1\prec\chi\prec\chi_2$ and $\supp(\chi)\subseteq U$. 

\smallskip

For each $i,j\in\{1,\dots,k\}$ let $A_{ij}(x,\dfo_x)$ denote the $ij$th entry of $\mathcal A(x,\dfo_x)$ and $A_{ij}(x,\xi)\in\psym{\mu_j-\nu_i}$ its (full) symbol. Also let $\mathcal A(x,\xi)$ denote the $k\times k$ matrix with entries $A_{ij}(x,\xi)$ and $A(x,\xi)$ the determinant of this matrix; in particular $A(x,\xi)\in\psym{m}$ where $m=\sumord{\mu\negvtsp}{\negvtsp\nu}$. Now, by the definition of uniform ellipticity, we can find $B(x,\xi)\in\psym{-m}$ such that $\chi_1(x)A(x,\xi)B(x,\xi)\equiv\chi_1(x)\pmod{\psym{-1}}$. Let $A^\dagger_{ij}(x,\xi)$ denote the $ij$th cofactor of the matrix $\mathcal A(x,\xi)$. Hence $A^\dagger_{ij}(x,\xi)\in\psym{m-\mu_j+\nu_i}$. Also set $B_{ij}(x,\xi)=(-1)^{i+j}B(x,\xi)A^\dagger_{ji}(x,\xi)\in\psym{\nu_j-\mu_i}$, define $B_{ij}(x,\dfo_x)\in\pdo{\nu_j-\mu_i}$ to be the pseudo-differential operator with symbol $B_{ij}(x,\xi)$ and let $\mathcal B$ denote the $k\times k$ system of pseudo-differential operators with entries $B_{ij}$. 

If $i,j,l\in\{1,\dots,k\}$ the differential operator $B_{il}\chi_1A_{lj}$ is contained in $\pdo{\mu_j-\mu_i}$ and has symbol $\chi_1(x)B_{il}(x,\xi)A_{lj}(x,\xi)\pmod{\psym{\mu_j-\mu_i-1}}$. Therefore the $ij$th entry of $\mathcal B\chi_1\mathcal A$ is contained in $\pdo{\mu_j-\mu_i}$ and, modulo an element of $\psym{\mu_j-\mu_i-1}$, has symbol
\[
\sum_{l=1}^k\chi_1(x)B_{il}(x,\xi)A_{lj}(x,\xi)
\;=\;\chi_1(x)\delta_{ij}A(x,\xi)B(x,\xi)
\;\equiv\;\chi_1(x)\delta_{ij}.
\]
Hence 
\begin{equation}
\label{whatBCaredo}
\mathcal B\chi_1\mathcal A
\;=\;\chi_1+\mathcal C
\end{equation} 
where $\mathcal C$ is a $k\times k$ matrix of pseudo-differential operators whose $ij$th entry is contained in $\pdo{\mu_j-\mu_i-1}$. In particular $\mathcal B$ and $\mathcal C$ define continuous maps 
\begin{equation}
\label{howBCmapthng}
\mathcal B:E^\nu\to E^\mu
\quad\text{and}\quad
\mathcal C:(E^{-1})^\mu\to E^\mu.
\end{equation}

\smallskip

Using the assumption of uniform ellipticity on $\mathcal A$ we can choose $\xi_0\in\R^n$ so that $A(x,\xi_0)=\det\mathcal A(x,\xi_0)$ is bounded away from 0 uniformly for $x\in U$. For each $i,j\in\{1,\dots,k\}$ set $D_{ij}(x)=(-1)^{i+j}\chi(x)A^\dagger_{ji}(x,\xi_0)/A(x,\xi_0)\in C^\infty$ and let $\mathcal D$ denote the $k\times k$ matrix of multiplication operators whose $ij$th entry is given by 
\begin{equation}
\label{defoftheopcalD}
\begin{cases}
D_{ij}(x)&\text{if $\mu_i=0$}\\
0&\text{if $\mu_i>0$.}
\end{cases}
\end{equation}
Now, as matrices, $\mathcal D(x)\chi(x)\mathcal A(x,\xi_0)=\chi(x)(I-\pi_\mu)$ (where $\pi_\mu$ is simply the diagonal matrix with a 1 in the $i$th position if $\mu_i>0$ and 0 otherwise). Since $(I-\pi_\mu)^2=I-\pi_\mu$ while $\mathcal A(x,\xi)(I-\pi_\mu)$ is independent of $\xi$ (n.b.\ if $\mu_j=0$ then $A_{ij}$ is a zeroth order differential operator) we have $\mathcal D\chi\mathcal A(I-\pi_\mu)=\chi(I-\pi_\mu)$ as differential operators. Coupled with the relation $\chi\prec\chi_2$ and the fact that $\mathcal A$ is a differential operator we then get 
\begin{equation}
\label{whatDisdo}
\chi
\;=\;\chi\pi_\mu+\mathcal D\chi\mathcal A(I-\pi_\mu)
\;=\;\chi\vtsp\pi_\mu\chi_2+\mathcal D\chi\vtsp\chi_2\mathcal A
-\mathcal D\chi\mathcal A\vtsp\pi_\mu\chi_2.
\end{equation}
On the other hand, Condition (A2), \eqref{defoftheopcalD} and the fact that $\nu_j\ge0$ clearly imply $\mathcal D$ defines a continuous map
\begin{equation}
\label{howDmapthng}
\mathcal D:(E^{-1})^\nu\to(E^{-1})^\mu.
\end{equation}

\smallskip

Using \eqref{whatBCaredo} and \eqref{whatDisdo}, the relations $\chi_1\prec\chi\prec\chi_2$ and the fact that $\mathcal A$ is a differential operator, we get 
\begin{equation}
\label{fancyidforchi1}
\chi_1
\;=\;\mathcal B\chi_1\mathcal A\chi-\mathcal C\chi
\;=\;\mathcal B\chi_1\vtsp\chi_2\mathcal A
-\mathcal C\chi\vtsp\pi_\mu\chi_2
-\mathcal C\mathcal D\chi\vtsp\chi_2\mathcal A
+\mathcal C\mathcal D\chi\mathcal A\vtsp\pi_\mu\chi_2.
\end{equation}
We also observe that $\mathcal A$ defines a continuous map
\begin{equation}
\label{howAmapthng}
\mathcal A:(E^{-1})^\mu\to(E^{-1})^\nu.
\end{equation}
The first part of the result now follows if we apply \eqref{fancyidforchi1} to $u$ and use Condition (A2), the continuous inclusion $E\hookrightarrow E^{-1}$ and the mapping properties given by \eqref{howBCmapthng}, \eqref{howDmapthng} and \eqref{howAmapthng}. On the other hand we can combine norm estimates for the various continuous maps to get
\begin{eqnarray*}
\norm{\chi_1 u}_{E^\mu}
&\le&C\bigl(\norm{\chi_2\mathcal Au}_{E^\nu}+\norm{\pi_\mu\chi_2u}_{(E^{-1})^\mu}\bigr.\\
&&\qquad\qquad{}+\norm{\mathcal D\chi\vtsp\chi_2\mathcal A u}_{(E^{-1})^\mu}
+\norm{\mathcal D\chi\mathcal A\vtsp\pi_\mu\chi_2 u}_{(E^{-1})^\mu}\bigr)\\
&\le&C\bigl(\norm{\chi_2\mathcal Au}_{E^\nu}+\norm{\pi_\mu\chi_2u}_{(E^{-1})^\mu}
+\norm{\chi_2\mathcal A u}_{(E^{-1})^\nu}
+\norm{\mathcal A\vtsp\pi_\mu\chi_2 u}_{(E^{-1})^\nu}\bigr)\\
&\le&C\bigl(\norm{\chi_2\mathcal Au}_{E^\nu}+\norm{\pi_\mu\chi_2u}_{(E^{-1})^\mu}\bigr)
\end{eqnarray*}
for all such $u$.
\end{proof}

Suppose $\psi$ is a diffeomorphism on $\R^n$ which is linear outside some compact set. By considering the standard rules for transforming principal symbols of differential operators it is clear that a system of differential operators $\mathcal A$ is elliptic on an open set $U$ iff the operator $(\psi^{-1})^*\mathcal A\psi^*$ is elliptic on the open set $\psi(U)$. It follows that the definition of ellipticity can be applied to systems of differential operators on a smooth manifold. The next result is used to justify a later remark (n.b.~$S^{n-1}$ can be replaced by any compact manifold without boundary).

\begin{thm}
\label{equivtocptmanifolds}
Suppose $\mathcal A$ is a $k\times k$ system of differential operators on $S^{n-1}$ of order $(\mu,\nu)$. If $\mathcal A$ is elliptic and $u\in\mathscr D'(S^{n-1})$ satisfies $\mathcal Au\in E^\nu(S^{n-1})\cap F^\nu(S^{n-1})$ and $\pi_\mu u\in F^\mu(S^{n-1})$ for some admissible spaces $E$ and $F$ (on $\R^{n-1}$) then we also have $u\in E^\mu(S^{n-1})$. Furthermore
\[
\norm{u}_{E^\mu(S^{n-1})}
\;\le\;C\bigl(\norm{\mathcal A u}_{E^\nu(S^{n-1})}+\norm{\pi_\mu u}_{F^\mu(S^{n-1})}\bigr)
\]
for all such $u$.
\end{thm}

\begin{proof}
We can obtain the result with $F=E^{-1}$ by applying Theorem \ref{stanapriorireg} on coordinate charts and patching the conclusion together with the help of observations similar to Lemma \ref{cordequivnorm} and Remark \ref{enlargechirem}. Induction, Lemma \ref{incforfunsponSnm1}, Corollary \ref{EFlocincpositivever} and Proposition \ref{basicpropadmisssp} then complete the result for general $F$. 
\end{proof}

\section{Model problems}

\label{modprobsec}

\subsection{Elliptic operators on $\boldsymbol{\cly}$ and a priori estimates}

\label{modprobinfclynv}

\begin{defn}
A differential operator $\mathcal B=\mathcal B(t,\omega,\dfo_t,\dfo_\omega)$ on $\cly$ is said to be \emph{uniform} if its coefficients are contained in $C^\infty(\cly)$.
\end{defn}

Suppose $B$ is a uniform scalar operator on $\cly$ of order $m$ and $E$ is an admissible space. Using Lemma \ref{diffoponclysplem} and Proposition \ref{multibysymprop} (n.b.~$C^\infty(\cly)=\bigcap_{l\in\NZ}\clysp0{\Cnt l}$) it is clear that $B$ defines a continuous map $\clysp\beta{E^m}\to\clysp\beta E$ for any $\beta\in\R$. We also have that 
\begin{equation}
\label{conjbyebetat}
e^{\beta t} B(t,\omega,\dfo_t,\dfo_\omega)
\;=\;B(t,\omega,\dfo_t+i\beta,\dfo_\omega)(e^{\beta t}\,\cdot\,)
\end{equation}
for any $\beta\in\R$, where $B(t,\omega,\dfo_t+i\beta,\dfo_\omega)$ is once again a uniform operator on $\cly$. 
This observation essentially allows us to take $\beta=0$ when studying the properties of the map $B:\clysp\beta{E^m}\to\clysp\beta E$ for a uniform operator $B$.

\medskip

Let $\mathcal B$ be a uniform operator on $\cly$ of order $(\mu,\nu)$ with entries $B_{ij}$ for $i,j\in\{1,\dots,k\}$. As above let $\detsym{\mathcal B}(x,\xi)$ denote the determinant of the $k\times k$ matrix formed from the principal symbols of the operators $B_{ij}$. We shall say that $\mathcal B$ is a \emph{uniform elliptic operator on $\cly$} if we can find $C>0$ such that $\abs{\detsym{\mathcal B}(x,\xi)}\ge C\abs\xi^{\sumord\mu\nu}_\cly$ for all $(x,\xi)\in T^*\cly$, where $\abs\cdot_\cly$ denotes the norm on the fibres of the cotangent bundle $T^*\cly$ given by the Riemannian metric on $\cly$.

For $i,j\in\{1,\dots,k\}$, $B_{ij}$ is a uniform scalar operator on $\cly$ of order $\mu_j-\nu_i$, so the above discussion implies that $B_{ij}$ defines a continuous map $\clysp\beta{E^{\mu_j}}\to\clysp\beta{E^{\nu_i}}$ for any $\beta\in\R$. It follows that $\mathcal B$ defines a continuous map $\clyvsp\beta\mu E\to\clyvsp\beta\nu E$. We shall denote this map by $\mathcal B^{(E,\beta)}$ or $\mathcal B^{(\beta)}$ when we need to make clear the spaces $\mathcal B$ is acting between.

\medskip

\begin{lem}
\label{ellregonclyvsp1}
Suppose $E$ is an admissible space and $\beta\in\R$. If $\mathcal Bu\in\clyvsp\beta\nu E$ for some $u\in\clyvsp\beta\mu{E^{-1}}$ then we also have $u\in\clyvsp\beta\mu E$. Furthermore
\[
\norm{u}_{\clyvsp\beta\mu E}
\;\le\;C\bigl(\norm{\mathcal Bu}_{\clyvsp\beta\nu E}+\norm{\pi_\mu u}_{\clyvsp\beta\mu{E^{-1}}}\bigr)
\]
for all $u\in\clyvsp\beta\mu E$.
\end{lem}

\begin{proof}
Suppose $\chi_1\in C^\infty(S^{n-1})$ with $\supp(\chi_1)\ne S^{n-1}$. Thus we can choose a chart $(\psi,U)$ of $S^{n-1}$ with $\supp(\chi_1)\subset U$.  Let $(\Psi,\R\times U)$ be the corresponding chart of $\cly$ and choose further functions $\chi_2,\chi_3\in C^\infty_0(U)$ with $\chi_3\succ\chi_2\succ\chi_1$. Also define $\widetilde\chi_i\in C^\infty$ by $\widetilde\chi_i=\chi_i\circ\Psi^{-1}$ (extended by 0) for $i=1,2,3$. 

Using the transformation properties of the principal symbol of a differential operator under diffeomorphisms (see Section 18.1 in \cite{Hor3} for example) it is straightforward to check that the operator $(\Psi^{-1})^*\chi_3\mathcal B\chi_3\Psi^*$ is uniformly elliptic on the open set $\widetilde\chi_3^{-1}(1/2,\infty)$. Theorem \ref{stanapriorireg} then gives us the following; if $\widetilde\chi_2(\Psi^{-1})^*\chi_3\mathcal B\chi_3\Psi^*v\in E^\nu$ for some $v\in(E^{-1})^\mu$ then $v\in E^\mu$ and we have an estimate of the form
\[
\norm{\widetilde\chi_1 v}_{E^\mu}
\;\le\;C\bigl(\norm{\widetilde\chi_2(\Psi^{-1})^*\chi_3\mathcal B\chi_3\Psi^*v}_{E^\nu}
+\norm{\pi_\mu\widetilde\chi_2 v}_{(E^{-1})^\mu}\bigr).
\]
Putting $v=(\Psi^{-1})^*(\chi_3 u)$ and using the fact that $\mathcal B$ is a differential operator we get
\[
\widetilde\chi_2(\Psi^{-1})^*\chi_3\mathcal B\chi_3\Psi^*v
\;=\;(\Psi^{-1})^*\chi_2\chi_3\mathcal B\chi_3^2u
\;=\;(\Psi^{-1})^*\chi_2\mathcal B u.
\]
Combining these observations with Lemma \ref{cordequivnorm} we now get the following; if $\chi_2\mathcal Bu\in\clyvsp0\nu E$ for some $u\in\clyvsp0\mu{E^{-1}}$ then $\chi_1 u\in\clyvsp0\mu E$ and we have an estimate of the form
\[
\norm{\chi_1 u}_{\clyvsp0\mu E}
\;\le\;C\bigl(\norm{\chi_2\mathcal Bu}_{\clyvsp0\nu E}
+\norm{\pi_\mu\chi_2 u}_{\clyvsp0\mu{E^{-1}}}\bigr)
\;\le\;C\bigl(\norm{\mathcal Bu}_{\clyvsp0\nu E}
+\norm{\pi_\mu u}_{\clyvsp0\mu{E^{-1}}}\bigr),
\]
the second inequality coming from Remark \ref{PropA2clysp} and the fact that $\chi_2\in C^\infty(\cly)$. The result now follows by allowing $\chi_1$ to vary of the elements of a suitable partition of unity on $S^{n-1}$.
\end{proof}

We need to make a couple of general observations before we proceed to the next regularity result.

\begin{rem}
\label{amigettopadlem1}
By Proposition \ref{multibysymprop} multiplication defines a continuous bilinear map $\clysp0{C^l}\times\clysp\beta E\to\clysp\beta E$ for all sufficiently large $l\in\NZ$. Using the translational invariance of the norm on $\clysp0{C^l}=C^l(\cly)$ and the fact that $\phi_{ij}=\phi_{j+1}-\phi_i$ we can bound $\norm{\phi_{ij}}_{\clysp0{C^l}}$ independently of $i,j\in\Z\cup\{\pm\infty\}$ with $i\le j$. If follows that multiplication by $\phi_{ij}$ defines a continuous map $\clysp\beta E\to\clysp\beta E$ whose operator norm can be bounded independently of $i,j$.
\end{rem}

\begin{rem}
\label{amigettopadlem2}
Suppose $\mathcal B$ is a model operator on $\cly$ of order $(\mu,\nu)$ and $\phi\in C^\infty(\cly)$. Then the $ij$th entry of the commutator $[\mathcal B,\phi]$ is the scalar model operator $[B_{ij},\phi]$ of order at most $\mu_j-\nu_i-1$; in particular $[B_{ij},\phi]=0$ whenever $\nu_i\ge\mu_j$ whilst $[\mathcal B,\phi]$ defines a continuous map $\clyvsp\beta\mu{E^{-1}}\to\clyvsp\beta\nu E$. Furthermore the operator norm of this map can be estimated by the $\clysp0{C^l}$ norm of $\phi$ and the coefficients of $\mathcal B$ for some sufficiently large $l\in\NZ$.
\end{rem}

\begin{lem}
\label{ellregonclyvsp2}
Suppose $E$ is an admissible space, $\beta\in\R$, $i,j\in\Z\cup\{\pm\infty\}$ with $i\le j$ and $l\in\NZ$. If $\phi_{i-l\,j+l}\mathcal Bu\in\clyvsp\beta\nu E$ and $\phi_{i-l\,j+l}u\in\clyvsp\beta\mu{E^{-l}}$ for some $u\in\mathscr D'(\cly)$ then we also have $\phi_{ij}u\in\clyvsp\beta\mu E$. Furthermore 
\[
\norm{\phi_{ij}u}_{\clyvsp\beta\mu E}
\;\le\;C\bigl(\norm{\phi_{i-l\,j+l}\mathcal Bu}_{\clyvsp\beta\nu E}
+\norm{\pi_\mu\phi_{i-l\,j+l} u}_{\clyvsp\beta\mu{E^{-l}}}\bigr)
\]
for all such $u$, where $C$ is independent of $i,j$.
\end{lem}

\begin{proof}
Induction clearly reduces the result to the case $l=1$. Now suppose we have $\phi_{i-1\,j+1}\mathcal Bu\in\clyvsp\beta\nu E$ and $\phi_{i-1\,j+1}u\in\clyvsp\beta\mu{E^{-1}}$ for some $u\in\mathscr D'(\cly)$. Therefore $\phi_{ij}u\in\clyvsp\beta\mu{E^{-1}}$ whilst $\phi_{ij}\prec\phi_{i-1\,j+1}$ and $\mathcal B$ is a differential operator so
\begin{equation}
\label{intermedeqntogetBphiijuincsp}
\mathcal B\phi_{ij}u
\;=\;\phi_{ij}\phi_{i-1\,j+1}\mathcal Bu+[\mathcal B,\phi_{ij}]\phi_{i-1\,j+1}u.
\end{equation}
By Remark \ref{amigettopadlem1} we know that multiplication by $\phi_{ij}$ defines a continuous map $\clyvsp\beta\nu E\to\clyvsp\beta\nu E$ with norm bounded independently of $i,j$. On the other hand Remark \ref{amigettopadlem2} implies $[\mathcal B,\phi_{ij}]$ maps $\clyvsp\beta\mu{E^{-1}}\to\clyvsp\beta\nu E$ continuously with operator norm bounded independently of $i,j$. Combining these observations with \eqref{intermedeqntogetBphiijuincsp} and our original hypothesis we thus get $\mathcal B\phi_{ij}u\in\clyvsp\beta\nu E$ with a norm estimate of the form
\[
\norm{\mathcal B\phi_{ij}u}_{\clyvsp\beta\nu E}
\;\le\;C\bigl(\norm{\phi_{i-1\,j+1}\mathcal Bu}_{\clyvsp\beta\nu E}
+\norm{\pi_\mu\phi_{i-1\,j+1}u}_{\clyvsp\beta\mu{E^{-1}}}\bigr),
\]
where $C$ is independent of $i,j$. Lemma \ref{ellregonclyvsp1} (applied to $\phi_{ij}u$) now implies $\phi_{ij}u\in\clyvsp\beta\mu E$ and 
\[
\norm{\phi_{ij}u}_{\clyvsp\beta\mu E}
\;\le\;C\bigl(\norm{\phi_{i-1\,j+1}\mathcal Bu}_{\clyvsp\beta\nu E}
+\norm{\pi_\mu\phi_{i-1\,j+1}u}_{\clyvsp\beta\mu{E^{-1}}}
+\norm{\pi_\mu\phi_{ij}u}_{\clyvsp\beta\mu{E^{-1}}}\bigr),
\]
where $C$ is independent of $i,j$. A further application of Remark \ref{amigettopadlem1} clearly completes the result.
\end{proof}

\subsection{Adjoint operators}

Suppose $\mathcal B$ be a uniform elliptic operator on $\cly$ of order $(\mu,\nu)$ with entries $B_{ij}$ for $i,j\in\{1,\dots,k\}$. Let $\mathcal B^*$ denote the formal adjoint of $\mathcal B$ (with respect to the volume measure $dt\d S^{n-1}$ on $\cly$); that is $\mathcal B^*$ is the $k\times k$ system of differential operators with entries 
\begin{equation}
\label{enteriesoffadj}
(B^*)_{ij}(t,\omega,\dfo_t,\dfo_\omega)
\;=\;\bigl(B_{ji}(t,\omega,\dfo_t,\dfo_\omega)\bigr)^*
\end{equation} 
for $i,j=1,\dots,k$. Define $m\in\NZ$ and vectors of integers $\overline{\mu}$ and $\overline{\nu}$ by 
\begin{equation}
\label{defofadjorders}
m=\max_i\mu_i
\;\quad\text{and}\;\quad
\overline{\mu}_i=m-\nu_i,\quad
\overline{\nu}_i=m-\mu_i
\quad\text{for $i=1,\dots,k$.}
\end{equation}
Therefore $\min_i\overline{\mu}_i\ge\min_i\overline{\nu}_i=0$ (n.b.~if the first inequality was not valid then we would have an $i$ such that $\mu_j-\nu_i<0$ for all $j\in\{1,\dots,k\}$; the $i$th row of $\mathcal B$ would then be 0, contradicting the assumption that $\mathcal B$ is elliptic). Furthermore 
\[
\order{(B^*)_{ij}}
\,=\,\order{B_{ji}}
\,=\,\mu_i-\nu_j
\,=\,\overline{\mu}_j-\overline{\nu}_i
\]
whilst it is easily seen that $\detsym{\mathcal B^*}(x,\xi)=\overline{\detsym{\mathcal B}(x,\xi)}$ for all $(x,\xi)\in T^*\cly$. It follows that $\mathcal B^*$ is a uniform elliptic operator on $\cly$ of order $(\overline{\mu},\overline{\nu})$. 

\bigskip

Let $E$ be an admissible space and set $F=E_0^*$. By Propositions \ref{dualforclysps} and \ref{basicpropadmisssp} we know that $F$ is also an admissible space whilst
\[
(\clyvsp\beta\mu{E_0})^*
\,=\,\prod_{i=1}^k\clysp{-\beta}F^{-\mu_i}
\,=\,\clyvsp{-\beta}{\overline{\nu}}{F^{-m}}
\;\quad\text{and}\;\quad
(\clyvsp\beta\nu{E_0})^*
\,=\,\prod_{i=1}^k\clysp{-\beta}F^{-\nu_i}
\,=\,\clyvsp{-\beta}{\overline{\mu}}{F^{-m}}
\]
for any $\beta\in\R$. Therefore the adjoint of the map $\mathcal B^{(E_0,\beta)}:\clyvsp\beta\mu{E_0}\to\clyvsp\beta\nu{E_0}$ is the map 
\[
\bigl(\mathcal B^{(E_0,\beta)}\bigr)^*
\,=\,(\mathcal B^*)^{(F^{-m},-\beta)}
\,:\,\clyvsp{-\beta}{\overline{\mu}}{F^{-m}}\to\clyvsp{-\beta}{\overline{\nu}}{F^{-m}}
\quad\;\text{where $F=E_0^*$.}
\]

\subsection{Operator pencils}

\label{oppensec}

\begin{defn}
A \emph{model operator on $\cly$} is a differential operator $\mathcal B=\mathcal B(\omega,\dfo_t,\dfo_\omega)$ on $\cly$ whose coefficients are independent of $t$. If $\mathcal B$ is also elliptic we shall refer to it as a \emph{model elliptic operator on $\cly$}.
\end{defn}

Clearly any model operator on $\cly$ is a uniform differential operator whilst any model elliptic operator is also a uniform elliptic operator. Now let $\mathcal B$ be a model elliptic operator on $\cly$ of order $(\mu,\nu)$. For $i,j\in\{1,\dots,k\}$ we can write the $ij$th entry of $\mathcal B$ in the form 
\begin{equation}
\label{formofcptofoppen1}
B_{ij}(\omega,\dfo_t,\dfo_\omega)
\;=\;\sum_{l=0}^{\mu_j-\nu_i} B_{ij}^{\mu_j-\nu_i-l}(\omega,\dfo_\omega)\,\dfo_t^l
\end{equation}
where $B_{ij}^l(\omega,\dfo_\omega)$ is a differential operator on $S^{n-1}$ of order at most $l$.

For each $\lambda\in\C$ let $\mathfrak B(\lambda)$ denote the differential operator on $S^{n-1}$ given by $\mathfrak B(\lambda)(\omega,\dfo_\omega)=\mathcal B(\omega,\lambda,\dfo_\omega)$; that is $\mathfrak B(\lambda)$ is the $k\times k$ system of differential operators with entries given by 
\begin{equation}
\label{formofcptofoppen2}
B_{ij}(\omega,\lambda,\dfo_\omega)
\;=\;\sum_{l=0}^{\mu_j-\nu_i} B_{ij}^{\mu_j-\nu_i-l}(\omega,\dfo_\omega)\,\lambda^l.
\end{equation}
Given an admissible space $E$ on $\R^{n-1}$, the operator $B_{ij}^l$ defines a continuous map $E^{\mu_j}(S^{n-1})\to E^{\nu_i}(S^{n-1})$ for $l=0,\dots,\mu_j-\nu_i$. Thus $\mathfrak B$ defines an operator valued function (or \emph{operator pencil}) 
\begin{equation}
\label{defofgenoppen}
\mathfrak B:\C\to\mathscr L\bigl(E^\mu(S^{n-1}),\vtsp E^\nu(S^{n-1})\bigr).
\end{equation}
The \emph{spectrum} $\fspec{\mathfrak B}$ of the operator pencil $\mathfrak B$ is defined to be the set of all $\lambda\in\C$ for which $\mathfrak B(\lambda):E^\mu(S^{n-1})\to E^\nu(S^{n-1})$ is not an isomorphism. Suppose $\lambda_0\in\fspec{\mathfrak B}$. Any non-zero element in the kernel of $\mathfrak B(\lambda_0)$ is called an eigenfunction of $\mathfrak B$ and the dimension of the kernel is the \emph{geometric multiplicity} of $\lambda_0$. A collection of functions $\phi_0,\dots,\phi_{M-1}\in E^\mu(S^{n-1})$ is called a \emph{Jordan chain} corresponding to $\lambda_0\in\fspec{\mathfrak B}$ iff $\phi_0$ is an eigenfunction corresponding to $\lambda_0$ and the meromorphic function $\Phi(\lambda)$ defined by 
\begin{equation}
\label{yabdef1}
\Phi(\lambda)
\;=\;\sum_{j=0}^{M-1}\frac{\phi_j}{(\lambda-\lambda_0)^{M-j}}
\end{equation}
satisfies
\begin{equation}
\label{yabdef2}
\mathfrak B(\lambda)\Phi(\lambda)
\,=\,O(1)
\quad\text{as $\lambda\to\lambda_0$.}
\end{equation}
It can be seen that this is equivalent to the condition
\begin{equation}
\label{yabdef3}
\sum_{j=0}^{M'}\frac1{j!}\,\vtsp\partial_\lambda^{\vtsp j}\mathfrak B(\lambda_0)\,\phi_{M'-j}
\;=\;0
\end{equation}
for $M'=0,\dots,M-1$. The functions $\phi_1,\dots,\phi_{M-1}$ are called generalised eigenfunctions. The \emph{algebraic multiplicity} of $\lambda_0$ is defined to be the dimension of the space of all functions of the form \eqref{yabdef1} which satisfy \eqref{yabdef2}.

\begin{rem}
\label{soothnessofeigenstuff}
It is straightforward to check that $\mathfrak B(\lambda)$ is an elliptic operator of order $(\mu,\nu)$ on $S^{n-1}$ for any $\lambda\in\C$ (n.b.~the principal symbol of each entry of $\mathfrak B(\lambda)$ is independent of $\lambda$). Standard elliptic regularity results (see Theorem \ref{equivtocptmanifolds}) then imply the spectrum, eigenfunctions and generalised eigenfunctions of the operator pencil \eqref{defofgenoppen} are independent of the admissible space $E$. Furthermore the eigenfunctions and generalised eigenfunctions are all contained in $C^\infty(S^{n-1})$.
\end{rem}

Using the a priori estimates for the elliptic operator $\mathcal B$ on $\cly$ and general results for Fredholm operator pencils we can obtain general information about the structure of the spectrum of $\mathfrak B$. The following result was obtained in \cite{AN} (see Chapter V) and \cite{AV} for scalar operators, and appears in Section 1.2.1 of \cite{NP2} for systems of operators (see also \cite{GGK} or Appendix A to \cite{KM} for details about the algebraic multiplicities of the isolated points in $\fspec{\mathfrak B}$).

\begin{thm}
\label{fspecthm}
The spectrum of $\mathfrak B$ consists of isolated points of finite algebraic multiplicity. Furthermore there exists constants $C_1,C_2>0$ such that 
\[
\fspec{\mathfrak B}
\;\subset\;\bigl\{\lambda\in\C\;\big|\;\abs{\re\lambda}\le C_1\abs{\im\lambda}+C_2\bigr\}.
\]
\end{thm}

\begin{rem}
\label{fspecthmrem}
Theorem \ref{fspecthm} implies $\im\fspec{\mathfrak B}$ (the projection of $\fspec{\mathfrak B}$ onto the imaginary axis) consists of isolated points. Furthermore, given $\gamma\in\im\fspec{\mathfrak B}$, the total algebraic multiplicity of all those $\lambda\in\fspec{\mathfrak B}$ with $\im\lambda=\gamma$ is finite.
\end{rem}

\medskip

A system of Jordan chains
\begin{equation}
\label{consysJc}
\bigl\{\phi_{j,0},\dots,\phi_{j,M_j-1}\;\big|\;
j=1,\dots,J\bigr\}
\end{equation}
is called \emph{canonical} if $\{\phi_{1,0},\dots,\phi_{J,0}\}$ forms a basis for $\Ker\mathfrak B(\lambda_0)$ (i.e.~the geometric eigenspace of $\lambda_0$), $M_1\ge\dots\ge M_J$ and $M_1+\dots+M_J=M$, the algebraic multiplicity of $\lambda_0$. It follows that $J$ is simply the geometric multiplicity of $\lambda_0$ whilst $M_1,\dots,M_J$ are known as the \emph{partial algebraic multiplicities} of $\lambda_0$. It is well known (see \cite{GGK} or Appendix A to \cite{KM} for example) that a canonical system of Jordan chains exists for any $\lambda_0\in\fspec{\mathfrak B}$ of finite algebraic multiplicity.

Define functions $u_{j,m}\in C^\infty(\cly)_\loc$ by
\begin{equation}
\label{powerexpsol}
u_{j,m}(t,\omega)
\;=\;e^{i\lambda_0 t}\sum_{l=0}^m\frac1{l!}\,(it)^l\,\phi_{j,m-l}(\omega)
\end{equation}
for $j=1,\dots,J$, $m=0,\dots,M_j-1$. This collection of functions forms a basis for the set of `power-exponential' solutions of the equation $\mathfrak B(\lambda_0)u=0$.

\medskip

Let $\mathfrak B^*$ denote the operator pencil associated to the model elliptic operator $\mathcal B^*$ (the formal adjoint of $\mathcal B$). For any $\lambda\in\C$ and $i,j\in\{1,\dots,k\}$ the $ij$th entry of $\mathfrak B^*(\lambda)$ is thus given by $(B^*)_{ij}(\omega,\lambda,\dfo_\omega)$, where $(B^*)_{ij}(\omega,\dfo_t,\dfo_\omega)$ is the operator given by~\eqref{enteriesoffadj}. On the other hand, by taking formal adjoints of \eqref{formofcptofoppen1} and \eqref{formofcptofoppen2} (on $\cly$ and $S^{n-1}$ respectively) we have 
\[
(B^*)_{ij}(\omega,\lambda,\dfo_\omega)
\;=\;\bigl(B_{ji}(\omega,\overline{\lambda},\dfo_\omega)\bigr)^*\,.
\]
Therefore $\mathfrak B^*(\lambda)=\bigl(\mathfrak B(\overline{\lambda})\bigr)^*$, the formal adjoint of the operator $\mathfrak B(\overline{\lambda})$ on $S^{n-1}$. It follows that $\lambda_0\in\fspec{\mathfrak B}$ iff $\overline{\lambda_0}\in\fspec{\mathfrak B^*}$ with full agreement of geometric, algebraic and partial algebraic multiplicities. In fact, given a canonical system of Jordan chains \eqref{consysJc} corresponding to $\lambda_0\in\fspec{\mathfrak B}$ we can find a unique canonical system of Jordan chains 
\begin{equation}
\label{adjconsysJc}
\bigl\{\psi_{j,0},\dots,\psi_{j,M_j-1}\;\big|\;
j=1,\dots,J\bigr\}
\end{equation}
corresponding to $\overline{\lambda_0}\in\fspec{\mathfrak B^*}$ which satisfy the biorthogonality conditions
\begin{equation}
\label{biorthogcond}
\sum_{l=0}^m\sum_{l'=0}^{m'}
\frac1{(l+l'+1)!}\,\bigipd{\partial_\lambda^{\vtsp l+l'+1}\mathfrak B(\lambda_0)
\phi_{j,\vtsp m-l}}{\psi_{j',\vtsp m'-l'}}_{S^{n-1}}
\;=\;\delta_{j,j'}\delta_{M_j-1-m, m'}
\end{equation}
for $j,j'=1,\dots,J$, $m=0,\dots,M_j-1$ and $m'=0,\dots,M_{j'}-1$ (see \cite{GGK} or Appendix A to \cite{KM} for example; n.b.~using \eqref{yabdef3} and a similar expression for the $\psi_{j,m}$'s it is possible to rewrite the biorthogonality conditions in a number of different ways --- in particular, it is enough to consider \eqref{biorthogcond} for $m=M_j-1$, $j,j'=1,\dots,J$ and $m'=0,\dots,M_{j'}-1$).

As above, the power exponential solutions of $\mathfrak B^*(\overline{\lambda_0})v=0$ can be expanded in terms of a basis given by the functions 
\begin{equation}
\label{adjpowerexpsol}
v_{j,m}(t,\omega)
\;=\;e^{i\overline{\lambda_0}t}\sum_{l=0}^m\frac1{l!}\,(it)^l\,\psi_{j,m-l}(\omega)
\end{equation}
for $j=1,\dots,J$, $m=0,\dots,M_j-1$.

\subsection{Isomorphisms}

\label{subsecisooncly}

The study of Fredholm properties of a certain class of elliptic operators on $\R^n$ (to be introduced in Section \ref{secclassofop}) can be largely reduced to the study of model elliptic operators on $\cly$. Theorems \ref{mainresclysp} and \ref{mainresclyspdetail} provide the key results concerning the latter. These (or closely related) results have been established by a number of authors, at least for the admissible spaces $E=\Sob pk$ with $p\in(1,\infty)$ and $k\in\Z$, and $E=\Zyg{s}_0$ with $s\in\R^+\!\setminus\!\N$ (in particular see \cite{AN}, \cite{AV}, \cite{Ko} and \cite{NP1}; see also \cite{KM} where such results are obtained as part of a general theory of differential equations with operator coefficients). Below we give a general argument to derive Theorems \ref{mainresclysp} and \ref{mainresclyspdetail} for an arbitrary admissible space $E$ from the case $E=\Sob 20$.

\begin{thm}
\label{mainresclysp}
Let $E$ be an admissible space, let $\mathcal B$ be a model elliptic operator on $\cly$ of order $(\mu,\nu)$ and define $\mathfrak B$ to the associated operator pencil. If $\beta\in\R\!\setminus\!\im\fspec{\mathfrak B}$ then $\mathcal B:\clyvsp\beta\mu E\to\clyvsp\beta\nu E$ is an isomorphism.
\end{thm}

\begin{thm}
\label{mainresclyspdetail}
Let $E$, $\mathcal B$ and $\mathfrak B$ be as in Theorem \ref{mainresclysp}. Suppose $\beta_1,\beta_2\in\R\!\setminus\!\im\fspec{\mathfrak B}$ with $\beta_1<\beta_2$. Set $\Sigma=\{\lambda\in\fspec{\mathfrak B}\,|\,\im\lambda\in[\beta_1,\beta_2]\}$ and, for each $\lambda\in\Sigma$, let $\{\phi^\lambda_{j,0},\dots,\phi^\lambda_{j,M_{\lambda,j}-1}\,|\,j=1,\dots,J_\lambda\}$ and $\{\psi^\lambda_{j,0},\dots,\psi^\lambda_{j,M_{\lambda,j}-1}\,|\,j=1,\dots,J_\lambda\}$ be canonical systems of Jordan chains corresponding to $\lambda\in\fspec{\mathfrak B}$ and $\overline{\lambda}\in\fspec{\mathfrak B^*}$ respectively, which satisfy the biorthogonality condition \eqref{biorthogcond}. Define $u^\lambda_{j,m}$ and $v^\lambda_{j,m}$ to be the corresponding power exponential functions given by \eqref{powerexpsol} and \eqref{adjpowerexpsol} respectively. Now let $f\in\clyvsp{\beta_1\!}\nu E\cap\clyvsp{\beta_2\!}\nu E$ and let $u^{(i)}\in\clyvsp{\beta_i\!}\mu E$ be the corresponding solutions of $\mathcal Bu=f$ for $i=1,2$. Then we have
\begin{equation}
\label{expnofdiff}
u^{(1)}(t,\omega)-u^{(2)}(t,\omega)
\;=\;\sum_{\lambda\in\Sigma}\sum_{j=1}^{J_\lambda}\sum_{m=0}^{M_{\lambda,j}-1}
c^\lambda_{j,m}(f)\, u^\lambda_{j,M_{\lambda,j}-1-m}(t,\omega)
\end{equation}
where the coefficient functions are given by
\begin{equation}
\label{formofcoeffs}
c^\lambda_{j,m}(f)
\;=\;\bigipd{f}{\vtsp iv^\lambda_{j,m}}_{\cly}
\end{equation}
for $\lambda\in\Sigma$, $j=1,\dots,J_\lambda$ and $m=0,\dots,M_{\lambda,j}-1$.
\end{thm}

\medskip

We begin by using a straightforward regularity argument to obtain Theorem \ref{mainresclysp} for the spaces $\Sob 2s$, $s\in\R$ from the case $s=0$.

\begin{prop}
\label{propofpartresHs}
Theorem \ref{mainresclysp} holds with $E=\Sob 2s$ for any $s\in\R$.
\end{prop}

\begin{proof}
Theorem \ref{mainresclysp} is known to hold for $E=\Sob 20=\Sob 2{{0}}$; see Theorem 3.1.1 in \cite{NP2} for example. Duality and induction now reduce our task to proving that if Theorem~\ref{mainresclysp} holds for $E=\Sob 2s$, $s\in\R$, then it also holds for $E=\Sob 2{s+\delta}$ where $\delta\in[0,1]$. 

Denote the maps $\mathcal B:\clyvsp\beta\mu{\Sob 2s}\to\clyvsp\beta\nu{\Sob 2s}$ and $\mathcal B:\clyvsp\beta\mu{\Sob 2{s+\delta}}\to\clyvsp\beta\nu{\Sob 2{s+\delta}}$ by $\mathcal B_0$ and $\mathcal B_\delta$ respectively, and assume $\mathcal B_0$ is an isomorphism. Since $\clyvsp\beta\mu{\Sob 2{s+\delta}}\subseteq\clyvsp\beta\mu{\Sob 2s}$, $\mathcal B_\delta$ must also be injective. Now let $f\in\clyvsp\beta\mu{\Sob 2{s+\delta}}\subseteq\clyvsp\beta\mu{\Sob 2s}$ and set $u=(\mathcal B_0)^{-1}f$. Therefore $u\in\clyvsp\beta\mu{\Sob 2s}\subseteq\clyvsp\beta\mu{\Sob 2{s+\delta-1}}$ whilst $\mathcal Bu=f\in\clyvsp\beta\nu{\Sob 2{s+\delta}}$. Lemma \ref{ellregonclyvsp1} then implies $u\in\clyvsp\beta\mu{\Sob 2{s+\delta}}$, thereby establishing the surjectivity of $\mathcal B_\delta$. The Open Mapping Theorem now shows that $\mathcal B_\delta$ is an isomorphism.
\end{proof}

Let $\mathscr O^{0,1}$ denote the set of scalar differential operators on $\cly$ of order 0 or 1 whose first order coefficients are independent of $t$ and whose zeroth order coefficients are linear in $t$. In what follows we use the notation $[P,Q]=PQ-QP$ for the commutator of operators $P$ and $Q$.

\begin{rem}\begin{claim}[(i)]
\label{remaboutcomm}
If $P\in\mathscr O^{0,1}$ and $\chi\in C^\infty(S^{n-1})$ then $[P,\chi]\in C^\infty(S^{n-1})$; that is, the operator $[P,\chi]$ is given as multiplication by some function in $C^\infty(\cly)$ which is independent of $t$.
\end{claim}

\begin{claim}[(ii)]
If $m\in\N$ and $P_1,\dots,P_m\in\mathscr O^{0,1}$ then $\mathcal B':=[P_1,\dots,[P_m,\mathcal B]\dots]$ is a $k\times k$ matrix of differential operators on $\cly$ whose coefficients are independent of $t$ and whose $ij$th entry has order at most $\mu_j-\nu_i-m'$ where 
\[
m'\;=\;\sum_{k=1}^m\,(1-\order{P_k}).
\]
It follows that $\mathcal B'$ defines a continuous map $\clyvsp\beta\mu E\to\clyvsp\beta\nu{E^{m'}}$ for any admissible space $E$ and $\beta\in\R$. 
\end{claim}
\end{rem}

Let $\beta\in\R\!\setminus\!\im\fspec{\mathfrak B}$ and define $\mathcal C$ to be the continuous map $\clysp\beta{\mathscr S}\to\clysp\beta{\mathscr S'}$ obtained by taking the restriction of the inverse of the isomorphism $\mathcal B:\clyvsp\beta\mu{\Sob 20}\to\clyvsp\beta\nu{\Sob 20}$. For any $m\in\NZ$ let $\mathscr O^m_{\mathcal C}$ denote the set of operators obtained by taking finite sums of operators of the form 
\begin{equation}
\label{genformofterminOmC}
\chi_1\mathcal C\mathcal Q_1\mathcal C\mathcal Q_2\dots\mathcal C\mathcal Q_l\mathcal C\chi_2
\end{equation}
where $\chi_1,\chi_2\in C^\infty(S^{n-1})$, $l\in\NZ$ and, for $i=1,\dots,l$, we can write
\[
\mathcal Q_i
\;=\;[P_{i1},\dots,[P_{im_i},\mathcal B]\dots]
\]
for some $m_i\in\NZ$ and $P_{i1},\dots,P_{im_i}\in\mathscr O^{0,1}$ with 
\[
\sum_{i=1}^l\sum_{k=1}^{m_i}\,(1-\order P_{ik})\;\le\;m.
\]

\begin{rem}
\label{remofpartresHscor}
For any $s\in\R$, Lemma \ref{denofCinft0incertsp} and Proposition \ref{propofpartresHs} imply that $\mathcal C$ has a unique extension to a continuous map $\clyvsp\beta\nu{\Sob 2s}\to\clyvsp\beta\mu{\Sob 2s}$. This observation can be coupled with Remark \ref{remaboutcomm}(ii) to show that the operator given by the formal expression \eqref{genformofterminOmC} defines a continuous map $\clyvsp\beta\nu{\Sob 2s}\to\clyvsp\beta\mu{\Sob 2{s+m}}$ for any $s\in\R$. Taking linear combinations it follows that any element of $\mathscr O^m_{\mathcal C}$ has similar mapping properties.
\end{rem}

\begin{lem}
\label{commlemm}
Let $\chi_1,\chi_2\in C^\infty(S^{n-1})$, $m\in\N$ and $P_1,\dots,P_m\in\mathscr O^{0,1}$. Set
\[
m'\;=\;\sum_{k=1}^m\,(1-\order{P_k})
\;\in\;\NZ.
\]
Then the operator 
\begin{equation}
\label{commversofOpC}
[P_1,\dots,[P_m,\chi_1\mathcal C\chi_2]\dots]
\end{equation}
has a unique extension to a continuous map $\clyvsp\beta\nu{\Sob 2s}\to\clyvsp\beta\mu{\Sob 2{s+m'}}$ for any $s\in\R$. 
\end{lem}

The operator $\mathcal C$ defines a continuous map $\clysp\beta{\mathscr S}\to\clysp\beta{\mathscr S'}$ whilst any $P_i$ defines continuous maps $\clysp\beta{\mathscr S}\to\clysp\beta{\mathscr S}$ and $\clysp\beta{\mathscr S'}\to\clysp\beta{\mathscr S'}$. It follows that the operator \eqref{commversofOpC} can be initially defined as a continuous map $\clysp\beta{\mathscr S}\to\clysp\beta{\mathscr S'}$.

\begin{proof}[Proof of Lemma \ref{commlemm}]
By Lemma \ref{denofCinft0incertsp} and Remark \ref{remofpartresHscor} it suffices to prove that the operator \eqref{commversofOpC} is contained in $\mathscr O^{m'}_{\mathcal C}$. Since we obviously have $\chi_1\mathcal C\chi_2\in\mathscr O^0_{\mathcal C}$ the following claim completes the result.

\begin{claim}[{Claim: If $\mathcal R\in\mathscr O^m_{\mathcal C}$ and $P\in\mathscr O^{0,1}$ then $[P,\mathcal R]\in\mathscr O^{m+1-\order P}_{\mathcal C}$.}] 
Clearly it suffices to prove the claim assuming $\mathcal R$ is given by \eqref{genformofterminOmC}. Then
\begin{eqnarray}
\lefteqn{[P,\mathcal R]
\,\ =\,\ [P,\chi_1]\mathcal C\mathcal Q_1\dots\mathcal Q_l\mathcal C\chi_2
\,\ +\,\ \sum_{i=1}^{l+1}\chi_1\mathcal C\mathcal Q_1\dots
\mathcal Q_{i-1}[P,\mathcal C]\mathcal Q_i\dots\mathcal Q_l\mathcal C\chi_2}\quad
\nonumber\\
\label{eqnthatihcnumasdag}
&&+\;{}\sum_{i=1}^l\chi_1\mathcal C\mathcal Q_1\dots
\mathcal C[P,\mathcal Q_i]\mathcal C\dots\mathcal Q_l\mathcal C\chi_2
\,\ +\,\ \chi_1\mathcal C\mathcal Q_1\dots\mathcal Q_l\mathcal C[P,\chi_2].\qquad
\end{eqnarray}
Now, by Remark \ref{remaboutcomm}(i), $[P,\chi_1]=\chi_1'$ and $[P,\chi_2]=\chi_2'$ for some $\chi_1',\chi_2'\in C^\infty(S^{n-1})$. On the other hand $[P,\mathcal C]=-\mathcal C[P,\mathcal B]\mathcal C$. Therefore each term in \eqref{eqnthatihcnumasdag} is contained in $\mathscr O^{m+1-\order P}_{\mathcal C}$.
\end{claim}
\end{proof}

Suppose $(\psi,U)$ is a chart for $S^{n-1}$ and let $(\Psi,\R\times U)$ be the corresponding chart for $\cly$. Also choose $\chi_1,\chi_2\in C^\infty_0(U)$; we shall use the same letters to denote the extensions of these functions to $\cly$ which are independent $t$. Define a $k\times k$ matrix of operators by $\mathcal D:=(\Psi^{-1})^*e^{\beta t}\chi_1\mathcal C\chi_2e^{-\beta t}\Psi^*$.

\begin{rem}
\label{movtoanfro}
Suppose $R$ is a scalar operator on $\cly$ which is supported on $\R\times U$ (in the sense that $R=\chi R\chi$ for some $\chi\in C^\infty_0(U)$). If $E$ and $F$ are admissible spaces then Lemma \ref{cordequivnorm} and \eqref{defofnormonclysp} show that $R$ defines a continuous map $R:\clysp\beta E\to\clysp\beta F$ iff $(\Psi^{-1})^*e^{\beta t}Re^{-\beta t}\Psi^*$ defines a continuous map $E\to F$. Clearly similar statements hold when $E$ or $F$ is replaced by $\mathscr S$ or $\mathscr S'$, or when $R$ is replaced by a matrix of operators.
\end{rem}

\begin{lem}
\label{Cisapdolem}
For each $i,j\in\{1,\dots,n\}$, the $ij$th entry of the operator $\mathcal D$ is contained in $\pdo{\nu_j-\mu_i}$.
\end{lem}

\begin{proof}
Since we know that $\mathcal D$ defines a continuous map $\mathscr S\to\mathscr S'$ (see Remark \ref{movtoanfro}) we can apply a characterisation of pseudo-differential operators given in \cite{B}. Let $m\in\N$ and suppose, for $k=1,\dots,m$, $Q_k$ is either the operator $\dfo_i$ or (multiplication by) $x_i$ for some $i\in\{1,\dots,n\}$. Let $m'\in\NZ$ be the number of $Q_k$ which are of zero order. Choose $\chi_0\in C^\infty(U)$ with $\chi_0\,\succ\,\chi_1,\chi_2$ and define $\widetilde\chi_i\in C^\infty$ by $\widetilde\chi_i=\chi_i\circ\Psi^{-1}$ (extended by 0) for $i=0,1,2$. Now the fact that $Q_1,\dots,Q_m$ are differential operators implies
\begin{equation}
\label{eqnthaiscallst2}
[Q_1,\dots,[Q_m,\mathcal D]\dots]
\;=\;\bigl[\widetilde\chi_0Q_1\widetilde\chi_0,\dots
\bigl[\widetilde\chi_0Q_m\widetilde\chi_0,\mathcal D\bigr]\dots\bigr].
\end{equation}
For $k=1,\dots,m$ set $P_k=e^{-\beta t}\Psi^*\widetilde\chi_0Q_k\widetilde\chi_0(\Psi^{-1})^*e^{\beta t}$. Therefore $P_k\in\mathscr O^{0,1}$ whilst \eqref{eqnthaiscallst2} gives 
\[
[Q_1,\dots,[Q_m,\mathcal D]\dots]
\;=\;(\Psi^{-1})^*e^{\beta t}[P_1,\dots,[P_m,\chi_1\mathcal C\chi_2]\dots]e^{-\beta t}\Psi^*.
\]
Combining this with Lemma \ref{commlemm} and Remark \ref{movtoanfro} it follows that the $ij$th entry of $[Q_1,\dots,[Q_m,\mathcal D]\dots]$ defines a continuous map $\Sob 2{s+\nu_j}\to\Sob 2{s+\mu_i+m'}$ for any $s\in\R$. Theorem 2.9 in \cite{B} now completes the result.
\end{proof}

\begin{proof}[Proof of Theorem \ref{mainresclysp}]
Let $\{\chi_i\}_{i\in I}$ be a finite partition of unity for $S^{n-1}$ such that $\supp(\chi_i)\cup\supp(\chi_j)\ne S^{n-1}$ for any $i,j\in I$. Thus we can write
\begin{equation}
\label{partunidecompofC}
\mathcal C
\;=\;\sum_{i,j\in I}\chi_i\mathcal C\chi_j\,.
\end{equation}
Now let $i,j\in I$ and choose a chart $(\psi,U)$ of $S^{n-1}$ with $\supp(\chi_i)\cup\supp(\chi_j)\subset U$. Let $(\Psi,\R\times U)$ be the corresponding chart of $\cly$. Set $\mathcal D_{ij}=(\Psi^{-1})^*e^{\beta t}\chi_1\mathcal C\chi_2e^{-\beta t}\Psi^*$, initially defined as a continuous map $\mathscr S\to\mathscr S'$. By Lemma \ref{Cisapdolem} $\mathcal D_{ij}\in\pdo{\nu_j-\mu_i}$ so $\mathcal D_{ij}$ defines a continuous map $\mathscr S\to\mathscr S$ which extends to give continuous maps $\mathscr S'\to\mathscr S'$ and $E^\nu\to E^\mu$ for any admissible space $E$. Furthermore any two of these extensions agree on functions common to their domains. Using \eqref{partunidecompofC} and Remark \ref{movtoanfro} it follows that $\mathcal C$ defines a continuous map $\mathcal C^{(\mathscr S)}:\clysp\beta{\mathscr S}\to\clysp\beta{\mathscr S}$ which has continuous extensions $\mathcal C^{(\mathscr S')}:\clysp\beta{\mathscr S'}\to\clysp\beta{\mathscr S'}$ and $\mathcal C^{(E)}:\clyvsp\beta\nu E\to\clyvsp\beta\mu E$ for any admissible space $E$. Furthermore
\begin{equation}
\label{eqngivingstabofinv}
\mathcal C^{(E)}f\,=\,\mathcal C^{(\mathscr S')}f
\quad\text{for any $f\in\clyvsp\beta\nu E$.}
\end{equation}
From the definition of $\mathcal C$ it clearly follows that the compositions $\mathcal C^{(\mathscr S)}\mathcal B^{(\mathscr S)}$ and $\mathcal B^{(\mathscr S)}\mathcal C^{(\mathscr S)}$ are both equal to the identity on $\clysp\beta{\mathscr S}$; that is, $\mathcal B^{(\mathscr S)}$ is an isomorphism with inverse $\mathcal C^{(\mathscr S)}$. On the other hand $\mathcal B^{(\mathscr S')}$ and $\mathcal C^{(\mathscr S')}$ are extensions of $\mathcal B^{(\mathscr S)}$ and $\mathcal C^{(\mathscr S)}$ so the compositions $\mathcal C^{(\mathscr S')}\mathcal B^{(\mathscr S')}$ and $\mathcal B^{(\mathscr S')}\mathcal C^{(\mathscr S')}$ must also give the identity map when restricted to $\clysp\beta{\mathscr S}$. Lemma \ref{denseofZSinZSp} then implies $\mathcal B^{(\mathscr S')}$ is an isomorphism with inverse $\mathcal C^{(\mathscr S')}$. Now let $E$ be any admissible space. Using \eqref{eqngivingstabofinv} and the obvious fact that $\mathcal B^{(E)}u=\mathcal B^{(\mathscr S')}u$ for any $u\in\clyvsp\beta\mu E$, we finally have that $\mathcal B^{(E)}$ is an isomorphism with inverse $\mathcal C^{(E)}$. 
\end{proof}

\medskip

We now turn our attention to Theorem \ref{mainresclyspdetail}, firstly establishing that the coefficient functions $c^\lambda_{j,m}$ given by \eqref{formofcoeffs} are in fact well defined. 

\begin{lem}
\label{thediffisalwayscnts}
Consider the notation of Theorems \ref{mainresclysp} and \ref{mainresclyspdetail} and let $\lambda\in\Sigma$, $j\in\{1,\dots,J_\lambda\}$ and $m\in\{0,\dots,M_{\lambda,j}-1\}$. Then the function $c^\lambda_{j,m}$ given by \eqref{formofcoeffs} defines a continuous map $\clysp{\beta_1}{\mathscr S'}\cap\clysp{\beta_2}{\mathscr S'}\to\C$. Using the continuous inclusions $\clyvsp{\beta_j\!}\nu E\hookrightarrow\clysp{\beta_j}{\mathscr S'}$ for $j=1,2$ (see Remark \ref{clyspforCinfty}), it follows that $c^\lambda_{j,m}$ also defines a continuous map $\clyvsp{\beta_1\!}\nu E\cap\clyvsp{\beta_2\!}\nu E\to\C$. 
\end{lem}

The topology on $\clyvsp{\beta_1\!}\nu E\cap\clyvsp{\beta_2\!}\nu E$ is smallest making the inclusion $\clyvsp{\beta_1\!}\nu E\cap\clyvsp{\beta_2\!}\nu E\hookrightarrow\clyvsp{\beta_j}\nu E$ continuous for $j=1,2$. A similar remark applies to $\clysp{\beta_1}{\mathscr S'}\cap\clysp{\beta_2}{\mathscr S'}$.

\begin{proof}
Using \eqref{adjpowerexpsol} we can write 
\[
v^\lambda_{j,m}(t,\omega)
\;=\;e^{i\overline{\lambda}t}\sum_{l=0}^m\frac1{l!}\,(it)^l\,\psi^\lambda_{j,m-l}(\omega),
\]
where $\psi^\lambda_{j,0},\dots,\psi^\lambda_{j,m}\in C^\infty(S^{n-1})$ (the fact that these functions are $\C^k$ valued is unimportant in what follows and won't be mentioned explicitly). Now, by assumption, $\im\lambda\in(\beta_1,\beta_2)$ so $\re(i\overline\lambda-\beta_2)<0$ and $\re(i\overline\lambda-\beta_1)>0$. Since $\supp(\phi_0)\subseteq(-1,+\infty)$ and $\supp(1-\phi_0)\subseteq(-\infty,1)$ it follows that 
\[
\phi_0\vtsp v^\lambda_{j,m}\;\in\;\clysp{-\beta_2}{\mathscr S}
\quad\text{and}\quad
(1\negvtsp-\negvtsp\phi_0)\vtsp v^\lambda_{j,m}
\;\in\;\clysp{-\beta_1}{\mathscr S}. 
\]
On the other hand, for any $f\in\clysp{\beta_1}{\mathscr S'}\cap\clysp{\beta_2}{\mathscr S'},$ \eqref{formofcoeffs} can be rewritten as
\[
c^\lambda_{j,m}(f)
\;=\;\bigipd{f}{i\phi_0\vtsp v^\lambda_{j,m}}_{\cly}
+\bigipd{f}{i(1\negvtsp-\negvtsp\phi_0)\vtsp v^\lambda_{j,m}}_{\cly}.
\]
The result now follows from the continuity of the dual pairings $\clysp{\beta_j}{\mathscr S'}\times\clysp{-\beta_j}{\mathscr S}\to\C$ for $j=1,2$ (see Remark \ref{propsA1A4innewform}).
\end{proof}

\begin{proof}[Proof of Theorem \ref{mainresclyspdetail}]
From the proof of Theorem \ref{mainresclysp} we know that, for $i=1,2$, the inverse of $\mathcal B^{(E,\beta_i)}$ is the restriction of the inverse of $\mathcal B^{(\mathscr S'\!,\beta_i)}$ to $\clyvsp\beta\nu E$. It therefore suffices to prove the result with $E$ replaced by $\mathscr S'$. Now Theorem \ref{mainresclyspdetail} holds with $E=\Sob 20$; see Theorems 3.1.4 and 3.2.1 in \cite{NP2} for example. It follows that \eqref{expnofdiff} is valid for all $f\in\clysp{\beta_1}{\mathscr S}\cap\clysp{\beta_2}{\mathscr S}$. Lemma \ref{thediffisalwayscnts} together with the following claim thus completes the proof.

\begin{claim}[Claim: The set $\clysp{\beta_1}{\mathscr S}\cap\clysp{\beta_2}{\mathscr S}$ is dense in $\clysp{\beta_1}{\mathscr S'}\cap\clysp{\beta_2}{\mathscr S'}$.] 
Let $f\in\clysp{\beta_1}{\mathscr S'}\cap\clysp{\beta_2}{\mathscr S'}$. Thus $(1-\phi_0)f=\phi_{-\infty0}f\in\clysp{\beta_1}{\mathscr S'}$ and $\phi_0f\in\clysp{\beta_2}{\mathscr S'}$. By Lemma \ref{denseofZSinZSp} we can then find sequences $\{f_{ji}\}_{i\in\N}\subset\clysp{\beta_j}{\mathscr S}$ for $j=1,2$ such that $f_{1i}\to\phi_{-\infty0}f$ in $\clysp{\beta_1}{\mathscr S'}$ and $f_{2i}\to\phi_0f$ in $\clysp{\beta_2}{\mathscr S'}$. Remark \ref{previouslemforZS} then gives $\phi_{-\infty1}f_{1i}\in\clysp{\beta_1}{\mathscr S}\cap\clysp{\beta_2}{\mathscr S}$ and $\phi_{-\infty1}f_{1i}\to\phi_{-\infty1}\phi_{-\infty0}f=\phi_{-\infty0}f$ in $\clysp{\beta_1}{\mathscr S'}\cap\clysp{\beta_2}{\mathscr S'}$. Similarly $\phi_{-1}f_{2i}\in\clysp{\beta_1}{\mathscr S}\cap\clysp{\beta_2}{\mathscr S}$ and $\phi_{-1}f_{2i}\to\phi_0f$ in $\clysp{\beta_1}{\mathscr S'}\cap\clysp{\beta_2}{\mathscr S'}$. Setting $f_i=\phi_{-\infty1}f_{1i}+\phi_{-1}f_{2i}$ for each $i\in\N$, it follows that $\{f_i\}_{i\in\N}$ is a sequence in $\clysp{\beta_1}{\mathscr S}\cap\clysp{\beta_2}{\mathscr S}$ which converges to $f$ in $\clysp{\beta_1}{\mathscr S'}\cap\clysp{\beta_2}{\mathscr S'}$.
\end{claim}
\end{proof}

\subsection{Model problems on $\pmb{\R}^{\boldsymbol n}_{\boldsymbol *}$}

\label{seconmodprobondon}

\begin{defn}
A scalar differential operator $A(x,\dfo_x)$ on $\don$ of order $m$ will be called \emph{uniform} if we can write
\begin{equation}
\label{relofmodopsdoncly}
A(x,\dfo_x)\;=\;r^{-m}(\Theta^{-1})^*B(t,\omega,\dfo_t,\dfo_\omega)\Theta^*
\;=\;r^{-m}B(\ln r,\omega,r\dfo_r,\dfo_\omega)
\end{equation}
for some uniform scalar operator $B(t,\omega,\dfo_t,\dfo_\omega)$ on $\cly$.

Suppose $\mathcal A(x,\dfo_x)$ is a $k\times k$ system of differential operators on $\don$ of order $(\mu,\nu)$. We will call $\mathcal A$ \emph{a uniform (elliptic) operator on $\don$} if we can write 
\begin{equation}
\label{defofunifellondon}
\mathcal A(x,\dfo_x)
\;=\;r^\nu(\Theta^{-1})^*\mathcal B(t,\omega,\dfo_t,\dfo_\omega)\Theta^*r^{-\mu}
\;=\;r^\nu\mathcal B(\ln r,\omega,r\dfo_r,\dfo_\omega)r^{-\mu}
\end{equation}
for some uniform (elliptic) operator $\mathcal B(t,\omega,\dfo_t,\dfo_\omega)$ on $\cly$ (where $r^\nu$ and $r^{-\mu}$ are the isomorphisms defined at the end of Section \ref{admissspace}). 
\end{defn}

\begin{rem}
If $B(t,\omega,\dfo_t,\dfo_\omega)$ is a differential operator on $\cly$ and $m\in\R$ then \eqref{simpleconjuncov} and \eqref{conjbyebetat} give 
\begin{eqnarray}
r^{-m}(\Theta^{-1})^*B(t,\omega,\dfo_t,\dfo_\omega)\Theta^*
&=&(\Theta^{-1})^*e^{-mt}B(t,\omega,\dfo_t,\dfo_\omega)\Theta^*
\nonumber\\
\label{effofmovrmmth}
&=&(\Theta^{-1})^*B(t,\omega,\dfo_t-im,\dfo_\omega)\Theta^*r^{-m}.
\end{eqnarray}
Now suppose $\mathcal A$ is a uniform operator on $\don$ of order $(\mu,\nu)$ and let $\mathcal B$ be the associated uniform operator on $\cly$ given by \eqref{defofunifellondon}. Then, for each $i,j\in\{1,\dots,k\}$, the $ij$th entries of $\mathcal A$ and $\mathcal B$ satisfy the relationship
\[
A_{ij}(x,\dfo_x)
\ =\ r^{\nu_i}(\Theta^{-1})^*B_{ij}(t,\omega,\dfo_t,\dfo_\omega)\Theta^*r^{-\mu_j}
\ =\ r^{\nu_i-\mu_j}(\Theta^{-1})^*B_{ij}(t,\omega,\dfo_t+i\mu_j,\dfo_\omega)\Theta^*.
\]
Therefore $A_{ij}$ is a uniform scalar operator on $\don$ of order $\mu_j-\nu_i$.
\end{rem}

\medskip

Suppose $A$ is a uniform scalar operator on $\don$ of order $m$. Using \eqref{Thetaisaniso}, Proposition \ref{basicchangbetaiso} and the mapping properties of uniform scalar operators on $\cly$ (see Section \ref{modprobinfclynv}), it is clear that $A$ defines a continuous map $\donsp\beta{E^m}\to\donsp{\beta-m}E$ for any admissible space $E$ and $\beta\in\R$. On the other hand, if $\mathcal A$ is a uniform operator on $\don$ of order $(\mu,\nu)$, \eqref{icbtdahal1}, \eqref{icbtdahal2} and the mapping properties of uniform operators on $\cly$ imply that $\mathcal A$ defines a continuous map $\mathcal A:\donvsp\beta\mu E\to\donvsp\beta\nu E$ for any admissible space $E$ and $\beta\in\R$. We shall denote this map by $\mathcal A^{(E,\beta)}$ or $\mathcal A^{(\beta)}$ when we need to make clear the spaces $\mathcal A$ is acting between.

\medskip

\begin{prop}
\label{sec35afterdef320i}
A scalar differential operator $A$ on $\don$ of order $m$ is uniform iff the coefficient of $\dfo_x^\alpha$ is contained in $\donsp{m-\abs\alpha}{\Cnt l}$ for all $l\in\NZ$ and multi-indices $\alpha$ with $\abs{\alpha}\le m$. 
\end{prop}

\begin{proof}
It is straightforward to check that the set of uniform scalar operators on $\cly$ can be generated from $C^\infty(\cly)$ and the set of first order model operators by taking linear combinations of products. With the help of Lemma \ref{locrepofscamodop1} and the fact that scalar differential operators commute to leading order, we can thus write any uniform scalar operator $B$ of order $m$ on $\cly$ in the form
\begin{eqnarray*}
\lefteqn{
B(t,\omega,\dfo_t,\dfo_\omega)
\;=\;\sum_{\abs\alpha\le m}f_\alpha(t,\omega) 
\bigl(B_1(\omega,\dfo_t,\dfo_\omega)\bigr)^{\alpha_1}\dots
\bigl(B_n(\omega,\dfo_t,\dfo_\omega)\bigr)^{\alpha_n}}\quad\\
&&=\;\sum_{\abs\alpha\le m}g_\alpha(t,\omega)\vtsp e^{\abs{\alpha}t} 
\bigl(e^{-t}B_1(\omega,\dfo_t,\dfo_\omega)\bigr)^{\alpha_1}\dots
\bigl(e^{-t}B_n(\omega,\dfo_t,\dfo_\omega)\bigr)^{\alpha_n}
\end{eqnarray*}
where, for each multi-index $\alpha$ with $\abs\alpha\le m$, $f_\alpha,g_\alpha\in C^\infty(\cly)$. However $g\in C^\infty(\cly)$ iff $g\circ\Theta^{-1}\in\donsp0{\Cnt l}$ for all $l\in\NZ$, whilst $D_i=(\Theta^{-1})^*e^{-t}B_i\Theta^*$ for $i=1,\dots,n$. Proposition \ref{basicchangbetaiso} and \eqref{simpleconjuncov} thus complete the result.
\end{proof}

\begin{prop}
\label{sec35afterdef320ii}
Suppose $\mathcal A$ is a $k\times k$ systems of differential operators on $\don$. Then $\mathcal A$ is a uniform elliptic operator iff $\mathcal A$ is a uniform operator which is uniformly elliptic on $\don$ in the sense of Definition \ref{defnofellingenset}.
\end{prop}

\begin{proof}
Let $\mathcal B$ be the uniform operator on $\cly$ given by \eqref{defofunifellondon}. Using Proposition \ref{sec35afterdef320i} we thus have to show that 
\begin{equation}
\label{whatwehavetoshowequiv1}
\abs{\detsym{\mathcal A}(x,\xi)}
\;\ge\;C\abs{\xi}^{\sumord\mu\nu}
\quad{\text{for all $(x,\xi)\in\don\times\R^n$}}
\end{equation}
is equivalent to 
\begin{equation}
\label{whatwehavetoshowequiv2}
\abs{\detsym{\mathcal B}(y,\eta)}
\;\ge\;C\abs{\eta}_\cly^{\sumord\mu\nu}
\quad{\text{for all $(y,\eta)\in T^*\cly$,}}
\end{equation}
(where $\abs\cdot_{\cly}$ denotes the norm on the fibres of the cotangent bundle $T^*\cly$ given by the Riemannian metric on $\cly$).

For each $i,j\in\{1,\dots,k\}$ let $a_{ij}$ and $b_{ij}$ denote the principal symbols of the $ij$th entries of $\mathcal A$ and $\mathcal B$ respectively. Using \eqref{defofunifellondon} and the transformation properties of principal symbols under diffeomorphisms (see Section 18.1 in \cite{Hor3} for example) we have
\[
a_{ij}\bigl(\Theta(y),\xi\bigr)
\;=\;r^{\nu_i-\mu_j}\,b_{ij}\bigl(y,(\dfo\Theta(y))^*\xi\bigr)
\]
for all $y\in\cly$ and $\xi\in\R^n$, where $r=\abs{\Theta(y)}$ and $(\dfo\Theta(y))^*$ denotes the adjoint of the linear map $\dfo\Theta(y):T_y\cly\to\R^n$. By taking determinants of the $k\times k$ matrices with entries $a_{ij}$ and $b_{ij}$ we now get
\begin{equation}
\label{Iapptoonneoneq}
\detsym{\mathcal A}(\Theta(y),\xi)
\;=\;r^{-\sumord\mu\nu}\,\detsym{\mathcal B}\bigl(y,(\dfo\Theta(y))^*\xi\bigr)
\;=\;\detsym{\mathcal B}\bigl(y,r^{-1}(\dfo\Theta(y))^*\xi\bigr)
\end{equation}
for all $y\in\cly$ and $\xi\in\R^n$, where we have made use of the fact that $\detsym{\mathcal B}(y,\eta)$ is homogeneous of degree $\sumord{\mu\negvtsp}{\negvtsp\nu}$ in $\eta$. On the other hand, the definition of $\Theta$ gives
\[
r^{-1}\bigabs{(\dfo\Theta(y))^*\xi}_\cly\,=\,\abs{\xi}
\]
for all $y\in\cly$ and $\xi\in\R^n$. The equivalence of \eqref{whatwehavetoshowequiv1} and \eqref{whatwehavetoshowequiv2} follows from this identity and \eqref{Iapptoonneoneq}.
\end{proof}

\bigskip

The relationship between $\mathcal A$ and $\mathcal B$ given by \eqref{defofunifellondon}, in combination with \eqref{icbtdahal1} and \eqref{icbtdahal2}, allows us to obtain results for uniform elliptic operators on $\don$ from results for uniform elliptic operators on $\cly$. The next result comes directly from Lemma~\ref{ellregonclyvsp2} (with $j=+\infty$ and $l=1$).

\begin{lem}
\label{aprioriondonsp}
Suppose $\mathcal A$ is a uniform elliptic operator on $\don$ of order $(\mu,\nu)$, $E$ is an admissible space, $\beta\in\R$ and $i\in\Z\cup\{-\infty\}$. If $\zeta_{i-1}\mathcal Au\in\donvsp\beta\nu E$ and $\zeta_{i-1}u\in\donvsp\beta\mu{E^{-1}}$ for some $u\in\mathscr D'(\don)$ then we also have $\zeta_iu\in\donvsp\beta\mu E$. Furthermore
\[
\norm{\zeta_iu}_{\donvsp\beta\mu E}
\;\le\;C\bigl(\norm{\zeta_{i-1}\mathcal Au}_{\donvsp\beta\nu E}
+\norm{\pi_\mu\zeta_{i-1}u}_{\donvsp\beta\mu{E^{-1}}}\bigr)
\]
for all such $u$, where $C$ is independent of $i$.
\end{lem}

\medskip

We also need to transfer the results of Section \ref{subsecisooncly} from $\cly$ to $\don$. With this in mind we firstly make the following definitions.

\begin{defn}
We say that a uniform scalar operator $A$ on $\don$ is \emph{a model scalar operator on $\don$} if the operator $B$ given by \eqref{relofmodopsdoncly} is a model scalar operator on $\cly$ (i.e.~the coefficients of $B$ do not depend on $t$). Likewise, we say that a uniform elliptic operator $\mathcal A$ on $\don$ is \emph{a model elliptic operator on $\don$} if the operator $\mathcal B$ given by \eqref{defofunifellondon} is a model elliptic operator on $\cly$ (i.e.~the coefficients of $\mathcal B$ do not depend on $t$). We shall denote the associated operator pencil by $\mathfrak B_{\mathcal A}$ and put
\begin{equation}
\label{wwitoctwnal}
\Res{\mathcal A}\,=\,\im\fspec{\mathfrak B_{\mathcal A}}.
\end{equation}
\end{defn}

Theorem \ref{mainresclysp} can be immediately rewritten as follows.

\begin{thm}
\label{mressec21}
Let $E$ be an admissible space and let $\mathcal A$ be a model elliptic operator on $\cly$ of order $(\mu,\nu)$. If $\beta\in\R\!\setminus\!\Res{\mathcal A}$ then $\mathcal A:\donvsp\beta\mu E\to\donvsp\beta\nu E$ is an isomorphism.
\end{thm}

Finally Theorem \ref{mainresclyspdetail} can be coupled with basic properties of the power-ex\-po\-nen\-tial solutions given by their explicit form (see \eqref{powerexpsol}) and Remark \ref{soothnessofeigenstuff} to give the following results.

\begin{thm}
\label{mressec22}
Let $E$ and $\mathcal A$ be as in Theorem \ref{mressec21} and define $\mathfrak B_{\mathcal A}$ to the associated spectral pencil. Suppose $\beta_1,\beta_2\in\R\!\setminus\!\Res{\mathcal A}$ satisfy $\beta_1<\beta_2$. Set $\Sigma=\{\lambda\in\fspec{\mathfrak B_{\mathcal A}}\,|\,\im\lambda\in[\beta_1,\beta_2]\}$ and let $M$ denote the sum of the algebraic multiplicities of all the $\lambda\in\Sigma$. Then there exists a vector space $X_\Sigma\subset C^\infty_\loc(\don)$ such that;
\begin{list}{}{\listinit}
\item[\emph{(i)}\hf]
For each $f\in\donvsp{\beta_1\!}\nu E\cap\donvsp{\beta_2\!}\nu E$ we have $\bigl(\mathcal A^{(\beta_1)}\bigr)^{-1}\!f-\bigl(\mathcal A^{(\beta_2)}\bigr)^{-1}\!f\,\in\,X_\Sigma$.
\item[\emph{(ii)}\hf]
$\dim X_\Sigma=M$ and $\mathcal Au=0$ for each $u\in X_\Sigma$. 
\item[\emph{(iii)}\hf]
Given any $u\in X_\Sigma$ and $\zeta\in\XYfns$ we have $\zeta u\in\donvsp{\beta_1\!}\mu E$ and $(1-\zeta)u\in\donvsp{\beta_2\!}\mu E$.
\end{list}
\end{thm}

\begin{cor}
\label{corthatgul1}
If $[\beta_1,\beta_2]\cap\Res{\mathcal A}=\emptyset$ then $\bigl(\mathcal A^{(\beta_1)}\bigr)^{-1}\!f=\bigl(\mathcal A^{(\beta_2)}\bigr)^{-1}\!f$ for all $f\in\donvsp{\beta_1\!}\nu E\cap\donvsp{\beta_2\!}\nu E$.
\end{cor}

\section{Main results}

\label{mainressec}

\subsection{The class of operators}

\label{secclassofop}

\begin{defn}
\label{earlydefofadmiss}
Let $A(x,\dfo_x)$ be a scalar differential operator on $\R^n$ of order $m$. We call $A$ an \emph{admissible operator} if we can write
\begin{equation}
\label{formofadmissop1}
A(x,\dfo_x)
\;=\;\sum_{\abs\alpha\le m} p^{\alpha}(x)\vtsp\dfo_x^\alpha\,,
\end{equation}
where $p^\alpha\in\ssym{m-\abs\alpha}$ for each multi-index $\alpha$ with $\abs\alpha\le m$.
\end{defn}

Suppose $A$ be an admissible operator of order $m$. By Propositions \ref{lemoncontofpderiv} and \ref{contofmultbyssym} it is clear that $A$ defines a continuous map $\consp{\beta-m}{E^m}\to\consp\beta E$ for any admissible space $E$ and $\beta\in\R$. As in the Introduction, we define the \emph{principal part of $A$} to be the operator on $\don$ which is given by 
\[
A_0(x,\dfo_x)
\;=\;\sum_{\abs\alpha\le m} r^{\abs\alpha-m}a^{\alpha}(\omega)\vtsp\dfo_x^\alpha\,,
\]
where $a^\alpha\in C^\infty(S^{n-1})$ is the principal part of $p^\alpha$ for each multi-index $\alpha$. We can rewrite $A_0$ in the form 
\[
A_0(x,\dfo_x)
\;=\;r^{-m}\sum_{j=0}^m A_0^{m-j}(\omega,\dfo_\omega)\vtsp(r\dfo_r)^j\,,
\]
where, for $j=0,\dots,m$, $A_0^j(\omega,\dfo_\omega)$ is a differential operator on $S^{n-1}$ of order at most $j$. It follows that $A_0$ is a model operator on $\don$.

\begin{rem}
\label{rem3?onmyprep}
Suppose $A$ is a scalar admissible operator with principal part $A_0$. For any $f\in\genXYfns$ set $A_f=fA+(1-f)A_0$. If $f=1$ on a neighbourhood of 0 then Remark \ref{intoremscomtsec} implies $A_f$ is a scalar admissible operator with the same principal part as $A$. On the other hand, if $f=0$ on a neighbourhood of 0 then Remark \ref{intoremscomtsec} and Proposition \ref{sec35afterdef320i} imply $A_f$ is a uniform scalar operator on $\don$.
\end{rem}

\bigskip

We now define the class of elliptic differential operators on $\R^n$ which we would like to consider. 

\begin{defn}
\label{defnofadellop}
Let $\mathcal A$ be a $k\times k$ system of differential operators on $\R^n$ of order $(\mu,\nu)$. We call $\mathcal A$ an \emph{admissible elliptic operator on $\R^n$} if it satisfies the following conditions.
\begin{list}{}{\listinit}
\item[(i)\hf]
For $i,j\in\{1,\dots,k\}$ the $ij$th entry of $\mathcal A$ is a scalar admissible operator on $\R^n$ (of order $\mu_j-\nu_i$). 
\item[(ii)\hf]
$\mathcal A$ is a uniformly elliptic operator on $\R^n$ (in the sense of Definition \ref{defnofellingenset}).
\end{list}
For $i,j\in\{1,\dots,k\}$ let $(A_0)_{ij}$ denote the principal part of the $ij$th entry of $\mathcal A$. By the \emph{principal part of $\mathcal A$} we mean the $k\times k$ system of differential operators on $\don$ with entries $(A_0)_{ij}$; this operator will be denoted by $\mathcal A_0$. By the \emph{operator pencil associated to $\mathcal A$} we mean $\mathfrak B_{\mathcal A_0}$; this will also be denoted by $\mathfrak B_{\mathcal A}$. Finally, define $\Res{\mathcal A}\subset\R$ by
\[
\Res{\mathcal A}
\,=\,\Res{\mathcal A_0}
\,=\,\im\fspec{\mathfrak B_{\mathcal A}}.
\]
\end{defn}

\smallskip

Suppose $\mathcal A$ satisfies Condition (i) of Definition \ref{defnofadellop} and, for $i,j\in\{1,\dots,k\}$, let $A_{ij}$ denote the $ij$th entry of $\mathcal A$. From the discussion after Definition \ref{earlydefofadmiss} we know that $A_{ij}$ defines a continuous map $\consp{\beta-\mu_j}{E^{\mu_j}}\to\consp{\beta-\nu_i}{E^{\nu_i}}$ for any admissible space $E$ and $\beta\in\R$. It follows that $\mathcal A$ defines a continuous map $\convsp\beta\mu E\to\convsp\beta\nu E$. We shall denote this map by $\mathcal A^{(E,\beta)}$ or $\mathcal A^{(\beta)}$ when we need to make clear the spaces $\mathcal A$ is acting between.

\medskip

\begin{rem}
\label{verygenremAA0}
Suppose $\mathcal A$ satisfies Condition (i) of Definition \ref{defnofadellop}. For any $f\in\genXYfns$ set $\mathcal A_f=f\mathcal A+(1-f)\mathcal A_0$. Using Remark \ref{rem3?onmyprep} it is straightforward to check that $\mathcal A_f$ is a uniform operator on $\don$ if $f=0$ on neighbourhood of 0 whilst $\mathcal A_f$ satisfies Condition (i) of Definition \ref{defnofadellop} provided $f=1$ on a neighbourhood of 0 --- in this case, the principal part of $\mathcal A_f$ is simply $\mathcal A_0$. 
\end{rem}

\begin{rem}
\label{estofdiffofdetsym}
Suppose $\mathcal A$ satisfies Condition (i) of Definition \ref{defnofadellop}. Using Condition (ii) of Definition \ref{defofpsym} we can find a bounded non-negative function $b$ on $\R^n$ with $b(x)\to0$ as $\abs x\to\infty$ such that the coefficients of all the entries of the operator $\mathcal A-\mathcal A_0$ are bounded by $b(x)$ when $\abs x\ge 1$. If $f,g\in\genXYfns$ then $\mathcal A_f-\mathcal A_g=(f-g)(\mathcal A-\mathcal A_0)$ so Remark \ref{rem267plusabit} now leads to the estimate
\[
\bigabs{\detsym{\mathcal A_f}(x,\xi)-\detsym{\mathcal A_g}(x,\xi)}
\;\le\;C\vtsp b(x)\vtsp\abs{(f\negvtsp-\negvtsp g)(x)}\,\abs{\xi}^{\sumord{\mu}{\nu}}
\]
for all $\abs{x}\ge1$ and $\xi\in\R^n$, where $C$ may depend on $\mathcal A$ and $\mathcal A_0$ but not on $(x,\xi)$.
\end{rem}

\begin{rem}
\label{equivofcondii}
If an operator $\mathcal A$ satisfies Condition (i) of Definition \ref{defnofadellop} then Condition (ii) is equivalent to the following.
\begin{list}{}{\listinit}
\item[(ii)$'$\hf]
$\mathcal A$ is an elliptic operator on $\R^n$ whilst $\mathcal A_0$ is a model elliptic operator on $\don$.
\end{list}
In order to see this first observe that Condition (ii) means
\[
\abs{\detsym{\mathcal A}(x,\xi)}
\;\ge\;C\abs{\xi}^{\sumord\mu\nu}
\quad\text{for all $(x,\xi)\in\R^n\times\R^n$.}
\]
On the other hand $\mathcal A_0$ is a model operator (by Condition (i) and the discussion proceeding Remark \ref{rem3?onmyprep}) so Proposition \ref{sec35afterdef320ii} implies $\mathcal A_0$ is a model elliptic operator on $\don$ iff $\mathcal A_0$ is uniformly elliptic on $\don$ in the sense of Definition \ref{defnofellingenset}. Thus Condition (ii)$'$ means $\detsym{\mathcal A}(x,\xi)\ne0$ for any $(x,\xi)\in\R^n\times\don$ and 
\[
\abs{\detsym{\mathcal A_0}(x,\xi)}
\;\ge\;C\abs{\xi}^{\sumord\mu\nu}
\quad\text{for all $(x,\xi)\in\don\times\R^n$.}
\]
The equivalence of Conditions (ii) and (ii)$'$ now follows from Remark \ref{estofdiffofdetsym} (with $f=1$ and $g=0$) and the fact that $\detsym{\mathcal A_0}(x,\xi)$ is homogeneous of degree 0 in $x$. 
\end{rem}

\begin{rem}
\label{oneisneededinanopap}
If $\mathcal A$ is an admissible elliptic operator, the associated operator pencil $\mathfrak B_{\mathcal A}$ can be equivalently defined using the expression 
\begin{equation}
\label{jabllab}
\mathfrak B_{\mathcal A}(\lambda)\phi(\omega)
\;=\;r^{-i\lambda}r^{-\nu}\mathcal A_0\bigl(r^{i\lambda}r^\mu\phi(\omega)\bigr)
\end{equation}
for any $\lambda\in\C$ and function (or distribution) $\phi$ on $S^{n-1}$. In order to see this let $\mathcal B(\omega,\dfo_t,\dfo_\omega)$ denote the model elliptic operator on $\cly$ associated to $\mathcal A_0$ by \eqref{defofunifellondon} (see also Remark \ref{equivofcondii}). Thus
\[
\mathcal A_0(x,\dfo_x)
\;=\;r^\nu\mathcal B(\omega,r\dfo_r,\dfo_\omega)r^{-\mu}
\]
and so
\[
r^{-i\lambda}r^{-\nu}\mathcal A_0\bigl(r^{i\lambda}r^\mu\phi(\omega)\bigr)
\;=\;r^{-i\lambda}\mathcal B(\omega,r\dfo_r,\dfo_\omega)
\bigl(r^{i\lambda}\phi(\omega)\bigr)
\;=\;\mathcal B(\omega,\lambda,\dfo_\omega)\phi(\omega).
\]
However $\mathfrak B_{\mathcal A}(\lambda)(\omega,\dfo_\omega)=\mathcal B(\omega,\lambda,\dfo_\omega)$ by definition.
\end{rem}

\medskip

For the remainder of this section we suppose that $\mathcal A$ is an admissible elliptic operator on $\R^n$ with principal part $\mathcal A_0$.

\begin{rem}
\label{firremaboutranops}
By Remark \ref{verygenremAA0} we know that $\mathcal A_{\zeta_i}$ is a uniform operator on $\don$ for all $i\in\N$. On the other hand, Condition (ii)$'$ for $\mathcal A_0$ and Remark \ref{estofdiffofdetsym} (with $f=0$ and $g=\zeta_i$) imply $\mathcal A_{\zeta_i}$ is uniformly elliptic on $\don$ for all sufficiently large $i\in\N$. Coupled with Proposition~\ref{sec35afterdef320ii} it follows that $\mathcal A_{\zeta_i}$ is a uniform elliptic operator on $\don$ for all sufficiently large $i\in\N$. 
\end{rem}

\begin{rem}
\label{secremaboutranops}
By Remark \ref{verygenremAA0} we know that $\mathcal A_{\eta_i}$ satisfies Condition (i) of Definition \ref{defnofadellop} for all $i\in\N$. On the other hand, Condition (ii) for $\mathcal A$ and Remark \ref{estofdiffofdetsym} (with $f=1$ and $g=\eta_i=1-\zeta_i$) imply $\mathcal A_{\eta_i}$ is uniformly elliptic on $\R^n$ for all sufficiently large $i\in\N$. Thus $\mathcal A_{\eta_i}$ is an admissible elliptic operator on $\R^n$ for all sufficiently large $i\in\N$. Furthermore the coefficients of $\mathcal A_{\eta_i}$ and $\mathcal A_0$ agree on a neighbourhood of $\infty$. It follows that the principal part of $\mathcal A_{\eta_i}$ is simply $\mathcal A_0$, the principal part of $\mathcal A$.
\end{rem}

\medskip

Suppose $E$ is an admissible space and $\beta\in\R$. Then, for any $i\in\N$, the operator $\mathcal A_{\zeta_i}-\mathcal A_0$ defines a continuous map $\donvsp\beta\mu{E}\to\donvsp\beta\nu{E}$ (see Remark \ref{firremaboutranops}) whilst $\mathcal A_{\eta_i}-\mathcal A$ defines a continuous map $\convsp\beta\mu{E}\to\convsp\beta\nu{E}$ (see Remark \ref{secremaboutranops}). 

\begin{lem}
\label{formofdecayofdiff}
We have $\norm{\mathcal A_{\zeta_i}\negvtsp-\negvtsp\mathcal A_0}_{\mathscr L(\donvsp\beta\mu E,\vtsp\donvsp\beta\nu E)}\vtsp,\,\norm{\mathcal A_{\eta_i}\negvtsp-\negvtsp\mathcal A}_{\mathscr L(\convsp\beta\mu E,\vtsp\convsp\beta\nu E)}\longrightarrow 0$ as $i\to\infty$. 
\end{lem}

\begin{proof}
If $A$ is a scalar admissible operator of order $m$ on $\R^n$ then \eqref{estfordonspE1}, Propositions \ref{lemoncontofpderiv} and \ref{multibysymprop}, and Lemma \ref{lemnthreeoflsec} give
\begin{equation}
\label{embforoswalab}
\lim_{i\to\infty}\,\norm{\zeta_i(A\negvtsp-\negvtsp 
A_0)}_{\mathscr L\left(\donsp{\beta-m_1}{E^{m_1}},\vtsp\donsp{\beta-m_2}{E^{m_2}}\right)}
\ =\ 0\ =\ \lim_{i\to\infty}\,\norm{\zeta_i(A\negvtsp-\negvtsp 
A_0)}_{\mathscr L\left(\consp{\beta-m_1}{E^{m_1}},\vtsp\consp{\beta-m_2}{E^{m_2}}\right)}
\end{equation}
for any $m_1,m_2\in\Z$ with $m=m_1-m_2$ (n.b.~the operator $\zeta_i(A-A_0)$ can be interpreted as acting on either $\donsp{\beta-m_1}{E^{m_1}}$ or $\consp{\beta-m_1}{E^{m_1}}$ using Lemma \ref{equivofnorms1} and the expression $\zeta_i(A-A_0)=\zeta_i(A-A_0)\zeta_{i-1}$). On the other hand, we have 
\[
\mathcal A_{\zeta_i}\negvtsp-\negvtsp\mathcal A_0
\;=\;\zeta_i(\mathcal A\negvtsp-\negvtsp\mathcal A_0)
\;=\;\mathcal A\negvtsp-\negvtsp\mathcal A_{\eta_i}\vtsp.
\]
The result follows by applying \eqref{embforoswalab} to the individual entries of $\zeta_i(\mathcal A-\mathcal A_0)$.
\end{proof}

\subsection{General estimates and some regularity}

In this section we give some general estimates and regularity results for the admissible elliptic operator $\mathcal A$; here `general' refers to the fact that these results hold without restriction on $\beta$ (the index appearing in the weighted spaces $\consp\beta E$).

\begin{prop}
\label{regnorestwf}
Suppose $E$ is an admissible space, $l\in\N$ and $\beta\in\R$. If we have $\mathcal Au\in\convsp\beta\nu E$ for some $u\in\convsp\beta\mu{E^{-l}}$ then we also have $u\in\convsp\beta\mu E$. Furthermore
\begin{equation}
\label{estwfolfromearLem}
\norm{u}_{\convsp\beta\mu E}
\;\le\;C\bigl(\norm{\mathcal A u}_{\convsp\beta\nu E}
+\norm{\pi_\mu u}_{\convsp\beta\mu{E^{-l}}}\bigr)
\end{equation}
for all such $u$. 
\end{prop}

\begin{proof}
Using induction it clearly suffices to prove the result in the case $l=1$. Now, by Remark \ref{firremaboutranops}, we can choose $I\in\N$ sufficiently large so that the operator $\mathcal A_{\zeta_{I-2}}$ is a uniform elliptic operator on $\don$. However $\zeta_{I-1}\mathcal A_{\zeta_{I-2}}u=\zeta_{I-1}\mathcal Au$ whilst Lemma \ref{equivofnorms1} gives us $\zeta_{I-1}\mathcal Au\in\donvsp\beta\nu E$ and $\zeta_{I-1}u\in\donvsp\beta\mu{E^{-1}}$. Lemma \ref{aprioriondonsp} then implies $\zeta_Iu\in\donvsp\beta\mu E$ and
\[
\norm{\zeta_I u}_{\donvsp\beta\mu E}
\;\le\;C\bigl(\norm{\zeta_{I-1}\mathcal A u}_{\donvsp\beta\nu E}
+\norm{\pi_\mu\zeta_{I-1}u}_{\donvsp\beta\mu{E^{-1}}}\bigr). 
\]
On the other hand $\mathcal A$ is an elliptic operator (on $\R^n$) whilst Lemma \ref{equivofnorms1} gives us $\eta_{I+1}\mathcal Au\in E^\nu$ and $\eta_{I+1}u\in(E^{-1})^\mu$. Theorem \ref{stanapriorireg} then implies $\eta_Iu\in E^\mu$ and 
\[
\norm{\eta_I u}_{E^\mu}
\;\le\;C\bigl(\norm{\eta_{I+1}\mathcal A u}_{E^\nu}
+\norm{\pi_\mu\eta_{I+1}u}_{(E^{-1})^\mu}\bigr) 
\]
(n.b.~$\eta_I\prec\eta_{I+1}$). The result now follows from Lemma \ref{equivofnorms1}.
\end{proof}

\medskip

The next result gives the strongest form of the elliptic regularity which can be achieved for an arbitrary $\beta\in\R$.

\begin{thm}
\label{regnorest}
Let $E$ and $F$ be admissible spaces and $\beta,\gamma\in\R$. Suppose that either $\beta<\gamma$ or $\beta\le\gamma$ and there exists a continuous inclusion $F^l\hookrightarrow E$ for some $l\in\NZ$. If we have $\mathcal Au\in\convsp\beta\nu E\cap\convsp\gamma\nu F$ for some $u\in\convsp\gamma\mu F$ then we also have $u\in\convsp\beta\mu E$. Furthermore
\begin{equation}
\label{estwfolfromearProp}
\norm{u}_{\convsp\beta\mu E}
\;\le\;C\bigl(\norm{\mathcal A u}_{\convsp\beta\nu E}
+\norm{\pi_\mu u}_{\convsp\gamma\mu F}\bigr)
\end{equation}
for all such $u$. 
\end{thm}

\begin{proof}
Using the hypothesis, Corollary \ref{EFlocincpositivever} and Proposition \ref{basicpropadmisssp}, we can choose $l\in\N$ sufficiently large so that we either have $\beta<\gamma$ and a local inclusion $F_\loc\hookrightarrow E^{-l}_\loc$ or $\beta=\gamma$ and a continuous inclusion $F\hookrightarrow E^{-l}$. Proposition \ref{genformincforconsp} or Remark \ref{easyincrem} then gives us a continuous inclusion $\convsp\gamma\mu F\hookrightarrow\convsp\beta\mu{E^{-l}}$. The result now follows from Proposition \ref{regnorestwf}.
\end{proof}

\subsection{A semi-Fredholm property, further regularity and stability of the index}

Throughout this section $\mathcal A$ is an admissible elliptic operator on $\R^n$. Let $\mathcal A_0$ and $\mathfrak B_{\mathcal A}$ denote the principal part of $\mathcal A$ and the associated operator pencil respectively. We start with a result which provides a key step for both Fredholm properties and further regularity results for $\mathcal A$.

\begin{prop}
\label{sdissprop1}
Suppose $E$ is an admissible space and $\beta_1,\beta_2\in\R$ satisfy $\beta_1\le\beta_2$ and $[\beta_1,\beta_2]\cap\Res{\mathcal A}=\emptyset$. If $\mathcal Au\in\convsp{\beta_1}\nu E\cap\convsp{\beta_2}\nu E$ for some $u\in\convsp{\beta_1}\mu E$ then we also have $u\in\convsp{\beta_2}\mu E$. Furthermore $\norm{u}_{\convsp{\beta_2}\mu E}\le C\norm{u}_{\convsp{\beta_1}\mu E}$ for all such $u$.
\end{prop}

For the proof of this result we make use of the operators $\mathcal A_{\zeta_i}=\zeta_i\mathcal A+(1-\zeta_i)\mathcal A_0$ for $i\in\N$. In particular Lemma \ref{formofdecayofdiff} implies $\mathcal A_{\zeta_i}\to\mathcal A_0$ in $\mathscr L(\donvsp\beta\mu E,\vtsp\donvsp\beta\nu E)$ as $i\to\infty$. Since the set of invertible elements in $\mathscr L(\donvsp\beta\mu E,\vtsp\donvsp\beta\nu E)$ is open, Theorem \ref{mressec21} immediately implies the following result.

\begin{lem}
\label{lemvermressec21}
If $\beta\in\R\!\setminus\!\Res{\mathcal A}$ then $\mathcal A_{\zeta_i}:\donvsp\beta\mu E\to\donvsp\beta\nu E$ is an isomorphism for all sufficiently large $i\in\N$.
\end{lem}

Using Neumann series we can also obtain the following extension to Corollary~\ref{corthatgul1}.

\begin{lem}
\label{lemvercorthatgul1}
If $\beta_1,\beta_2$ are as in Proposition \ref{sdissprop1} then $\bigl(\mathcal A_{\zeta_i}^{\vtsp(\beta_1)}\bigr)^{-1}\!f=\bigl(\mathcal A_{\zeta_i}^{\vtsp(\beta_2)}\bigr)^{-1}\!f$ for all $f\in\donvsp{\beta_1}\nu E\cap\donvsp{\beta_2}\nu E$ and sufficiently large $i\in\N$.
\end{lem}

\begin{proof}
From Lemma \ref{formofdecayofdiff} we know that for all sufficiently large $i\in\N$ we have
\[
\norm{\mathcal A_{\zeta_i}-\mathcal A_0}_{\mathscr L(\donvsp{\beta_j}\mu E,\vtsp\donvsp{\beta_j}\nu E)}
\;\le\;\frac12\norm{\mathcal A_0}_{\mathscr L(\donvsp{\beta_j}\mu E,\vtsp\donvsp{\beta_j}\nu E)}
\]
for $j=1,2$. Choose $i\in\N$ for which this is true and set $\mathcal P=\mathcal A_{\zeta_i}-\mathcal A_0$. It follows that we have a Neumann series expansion
\begin{equation}
\label{neumannseriesexpn}
\bigl(\mathcal A_{\zeta_i}^{\vtsp(\beta_j)}\bigr)^{-1}
\;=\;\sum_{l=0}^\infty\bigl(\bigl(\mathcal A_0^{(\beta_j)}\bigr)^{-1}\mathcal P\bigr)^l
\,\bigl(\mathcal A_0^{(\beta_j)}\bigr)^{-1}
\end{equation}
which is convergent in $\mathscr L(\donvsp{\beta_j}\nu E,\vtsp\donvsp{\beta_j}\mu E)$ for $j=1,2$. 

Let $f\in\donvsp{\beta_1}\nu E\cap\donvsp{\beta_2}\nu E$ and set $u=\bigl(\mathcal A_0^{(\beta_1)}\bigr)^{-1}\!f$. Then
\[
u\;=\;\bigl(\mathcal A_0^{(\beta_1)}\bigr)^{-1}\!f
\;=\;\bigl(\mathcal A_0^{(\beta_2)}\bigr)^{-1}\!f
\;\in\;\donvsp{\beta_1}\mu E\cap\donvsp{\beta_2}\mu E
\]
by Corollary \ref{corthatgul1}. Now let $l\in\NZ$ and suppose 
\begin{equation}
\label{indstepjldoo}
\bigl(\bigl(\mathcal A_0^{(\beta_1)}\bigr)^{-1}\mathcal P\bigr)^lu
\;=\;\bigl(\bigl(\mathcal A_0^{(\beta_2)}\bigr)^{-1}\mathcal P\bigr)^lu
\;\in\;\donvsp{\beta_1}\mu E\cap\donvsp{\beta_2}\mu E.
\end{equation}
Therefore
\[
\mathcal P\bigl(\bigl(\mathcal A_0^{(\beta_1)}\bigr)^{-1}\mathcal P\bigr)^lu
\;=\;\mathcal P\bigl(\bigl(\mathcal A_0^{(\beta_2)}\bigr)^{-1}\mathcal P\bigr)^lu
\;\in\;\donvsp{\beta_1}\nu E\cap\donvsp{\beta_2}\nu E,
\]
so Corollary \ref{corthatgul1} implies \eqref{indstepjldoo} also holds with $l$ replaced by $l+1$. The result now follows from \eqref{neumannseriesexpn} and induction.
\end{proof}

\begin{proof}[Proof of Proposition \ref{sdissprop1}]
Choose $I\in\N$ so that the conclusion of Lemma \ref{lemvercorthatgul1} holds with $i=I$. Now suppose $\mathcal Au\in\convsp{\beta_1}\nu E\cap\convsp{\beta_2}\nu E$ for some $u\in\convsp{\beta_1}\mu E$. Set $v=\zeta_Iu$ and $g=\mathcal A_{\zeta_I}v$. By Lemma \ref{equivofnorms1} we have $v\in\donvsp{\beta_1}\mu E$ and thus $g\in\donvsp{\beta_1}\nu E$ (by the mapping properties of uniform operators on $\don$; see also Remark \ref{verygenremAA0}). On the other hand $g=0$ on a neighbourhood of 0 whilst
\[
\zeta_{I+1} g
\;=\;\zeta_{I+1}\bigl(\zeta_I\mathcal A+(1\negvtsp-\negvtsp\zeta_I)\mathcal A_0\bigr)\zeta_I u
\;=\;\zeta_{I+1}\mathcal A u
\;\in\;\donvsp{\beta_1}\nu E\cap\donvsp{\beta_2}\nu E
\]
by the hypothesis on $\mathcal Au$ and Lemma \ref{equivofnorms1}. Remark \ref{newequivofnormsrem2} thus gives $g\in\donvsp{\beta_1}\nu E\cap\donvsp{\beta_2}\nu E$ and so $\zeta_I u=v=\bigl(\mathcal A_{\zeta_I}^{(\beta_1)}\bigr)^{-1}g=\bigl(\mathcal A_{\zeta_I}^{(\beta_2)}\bigr)^{-1}g$ by Lemma \ref{lemvercorthatgul1}. Using the fact that $\mathcal A_{\zeta_I}^{(\beta_j)}$ is an isomorphism for $j=1,2$ we also get a norm estimate of the form
\[
\norm{\zeta_I u}_{\donvsp{\beta_2}\mu E}
\;\le\;C\norm{\zeta_I u}_{\donvsp{\beta_1}\mu E}\,.
\]
Lemma \ref{equivofnorms1} now completes the result.
\end{proof}

\bigskip

The next result establishes a slightly weaker form of the Fredholm property for $\mathcal A$ in which `Fredholm' is replaced with `semi-Fredholm'; by the latter we mean an operator which has a closed range and for which either the kernel or the cokernel is finite dimensional. The key steps in establishing this semi-Fredholm property are provided by Corollary \ref{specformofcpt}, Proposition \ref{regnorestwf} and Proposition \ref{sdissprop1}.

\begin{thm}
\label{wfomres1}
Let $E$ be an admissible space and $\beta\in\R\!\setminus\!\Res{\mathcal A}$. Then the operator $\mathcal A:\convsp\beta\mu E\to\convsp\beta\nu E$ is semi-Fredholm with a finite dimensional kernel.
\end{thm}

\begin{proof}
Choose $I\in\N$ so that the conclusion to Lemma \ref{lemvermressec21} holds for $i=I$. Now suppose we have a sequence satisfying
\begin{equation}
\label{seqcritocon}
\{u_i\}_{i\in\N}\subset\convsp\beta\mu E,\quad
\norm{u_i}_{\convsp\beta\mu E}=1,\quad
\text{$\mathcal Au_i\to0$ in $\convsp\beta\nu E$.}
\end{equation}
Let $i\in\N$. Lemma \ref{equivofnorms1} gives us $\zeta_I\mathcal Au_i\in\donvsp\beta\nu E$ so $\mathcal A_{\zeta_I}^{-1}(\zeta_I\mathcal Au_i)\in\donvsp\beta\mu E$. Setting $v_i=\zeta_I\mathcal A_{\zeta_I}^{-1}(\zeta_I\mathcal Au_i)$ we then have $v_i\in\convsp\beta\mu E$. Furthermore, by combining norm estimates given by Lemma \ref{equivofnorms1} and the fact that $\mathcal A_{\zeta_I}:\donvsp\beta\mu E\to\donvsp\beta\nu E$ is an isomorphism, we get 
\begin{equation}
\label{sfwntgttw1}
\norm{v_i}_{\convsp\beta\mu E}
\;\le\; C\norm{\mathcal Au_i}_{\convsp\beta\nu E}
\;\le\; C\norm{u_i}_{\convsp\beta\mu E}\,.
\end{equation}
Also, $\zeta_{I+1}\prec\zeta_I$ so $\zeta_{I+1}\mathcal A\zeta_I=\zeta_{I+1}\mathcal A_{\zeta_I}$ and hence 
\begin{equation}
\label{sfwntgttw2}
\zeta_{I+1}\mathcal Av_i
\;=\;\zeta_{I+1}\mathcal A
\bigl(\zeta_I\,\mathcal A_{\zeta_I}^{-1}(\zeta_I\mathcal Au_i)\bigr)
\;=\;\zeta_{I+1}\mathcal A_{\zeta_I}
\bigl(\mathcal A_{\zeta_I}^{-1}(\zeta_I\mathcal Au_i)\bigr)
\;=\;\zeta_{I+1}\mathcal Au_i.
\end{equation}
Now define $w_i$ by $w_i=u_i-v_i$ for all $i\in\N$. By the second inequality in \eqref{sfwntgttw1} $\{w_i\}_{i\in\N}$ is a bounded sequence in $\convsp\beta\mu E$ whilst \eqref{sfwntgttw2} gives us $\zeta_{I+1}\mathcal Aw_i=0$ for all $i\in\N$. Choose $\gamma>\beta$ so that $[\beta,\gamma]\cap\Res{\mathcal A}=\emptyset$. Using Remark \ref{newequivofnormsrem2} and Proposition \ref{sdissprop1} it follows that $\{w_i\}_{i\in\N}$ is a bounded sequence in $\convsp\gamma\mu E$. However $\convsp\gamma\mu E\hookrightarrow\convsp\beta\mu{E^{-1}}$ is a compact map by Corollary \ref{specformofcpt}, so we can choose a subsequence $\{w_{i(j)}\}_{j\in\N}$ which is convergent in $\convsp\beta\mu{E^{-1}}$. On the other hand, the last part of \eqref{seqcritocon} and the first inequality in \eqref{sfwntgttw1} imply $u_i-w_i=v_i\to0$ in $\convsp\beta\mu E$. Since $\convsp\beta\mu E\hookrightarrow\convsp\beta\mu{E^{-1}}$ continuously we then get the convergence of $\{u_{i(j)}\}_{j\in\N}$ in $\convsp\beta\mu{E^{-1}}$. By combining Proposition \ref{regnorestwf} and the last part of \eqref{seqcritocon} it follows that $\{u_{i(j)}\}_{j\in\N}$ is convergent in $\convsp\beta\mu E$. 

Summarising, we have shown that any sequence satisfying \eqref{seqcritocon} has a subsequence which is convergent in $\convsp\beta\mu E$. A standard argument (see Proposition 19.1.3 in \cite{Hor3} or Theorems IV.5.9, IV.5.10 and IV.5.11 in \cite{Ka}) shows this implies $\mathcal A:\convsp\beta\mu E\to\convsp\beta\nu E$ has a finite dimensional kernel and a closed range. 
\end{proof}

If we are only interested in the kernel of the map $\mathcal A^{(E,\beta)}$ the restriction on $\beta$ in Theorem \ref{wfomres1} can be dropped.

\begin{thm}
\label{genfdkerresult}
For any admissible space $E$ and $\beta\in\R$ the map $\mathcal A:\convsp\beta\mu E\to\convsp\beta\nu E$ has a finite dimensional kernel.
\end{thm}

\begin{proof}
Choose $\gamma\in\R\!\setminus\!\Res{\mathcal A}$ with $\gamma\le\beta$. By \eqref{cntsincforXE} we have a continuous inclusion $\convsp\beta\mu E\hookrightarrow\convsp\gamma\mu E$ so $\Ker\mathcal A^{(\beta)}\subseteq\Ker\mathcal A^{(\gamma)}$. On the other hand $\gamma\in\R\!\setminus\!\Res{\mathcal A}$ so $\Ker\mathcal A^{(\gamma)}$ must be finite dimensional by Theorem \ref{wfomres1}. 
\end{proof}

\bigskip

We complete this section with two further consequences of Proposition \ref{sdissprop1}, the first of which gives a stronger form elliptic regularity than Theorem \ref{regnorest} but under some restrictions on $\beta$.

\begin{thm}
\label{changebetareggen}
Suppose $E$ and $F$ be admissible spaces and let $\beta,\gamma\in\R$ belong to the same component of $\R\!\setminus\!\Res{\mathcal A}$. If $\mathcal Au\in\convsp\beta\nu E\cap\convsp\gamma\nu F$ for some $u\in\convsp\gamma\mu F$ then we also have $u\in\convsp\beta\mu E$. Furthermore
\[
\norm{u}_{\convsp\beta\mu E}
\;\le\;C\bigl(\norm{\mathcal A u}_{\convsp\beta\nu E}
+\norm{\pi_\mu u}_{\convsp\gamma\mu F}\bigr)
\]
for all such $u$. 
\end{thm}

\begin{proof}
Choose any $\lambda\in\R$ which lies in the same component of $\R\!\setminus\!\Res{\mathcal A}$ as $\beta$ and $\gamma$, and satisfies $\lambda<\beta,\gamma$. Now $u\in\convsp\gamma\mu F$ and $\mathcal Au\in\convsp\beta\nu E\cap\convsp\gamma\nu F$ so we also have $\mathcal Au\in\convsp\lambda\nu E$ and $\norm{\mathcal Au}_{\convsp\lambda\nu E}\le C\norm{\mathcal Au}_{\convsp\beta\nu E}$ by \eqref{cntsincforXE}. Theorem \ref{regnorest} then gives $u\in\convsp\lambda\mu E$ and 
\[
\norm{u}_{\convsp\lambda\mu E}
\;\le\;C\bigl(\norm{\mathcal Au}_{\convsp\lambda\nu E}
+\norm{\pi_\mu u}_{\convsp\gamma\mu F}\bigr)
\;\le\;C\bigl(\norm{\mathcal Au}_{\convsp\beta\nu E}
+\norm{\pi_\mu u}_{\convsp\gamma\mu F}\bigr).
\]
Applying Proposition \ref{sdissprop1} (with $\beta_1=\lambda$ and $\beta_2=\beta$) now gives $u\in\convsp\beta\mu E$ and $\norm{u}_{\convsp\beta\mu E}\le C\norm{u}_{\convsp\lambda\mu E}$, completing the result.
\end{proof}

\medskip

For a semi-Fredholm map we can define a $\Z\cup\{\pm\infty\}$ valued index; in particular we have that a semi-Fredholm map is Fredholm iff its index is finite (something which will be established for $\mathcal A$ in Theorem \ref{sfomres1}). The next result establishes the stability of this index over particular ranges of the weighted spaces on which $\mathcal A$ is acting; in fact Theorem \ref{mressec5} (below) shows these ranges to be maximal. 

\begin{thm}
\label{stabofindex}
Suppose $E$ and $F$ be admissible spaces and let $\beta,\gamma\in\R$ belong to the same component of $\R\!\setminus\!\Res{\mathcal A}$. Then the semi-Fredholm maps $\mathcal A:\convsp\beta\mu E\to\convsp\beta\nu E$ and $\mathcal A:\convsp\gamma\mu F\to\convsp\gamma\nu F$ have the same index.
\end{thm}

\begin{proof}
We get $\Ker\mathcal A^{(E,\beta)}=\Ker\mathcal A^{(F,\gamma)}$ as a direct consequence of Theorem \ref{changebetareggen}. It therefore remains to show $\codim\bigl(\Ran\mathcal A^{(E,\beta)},\vtsp\convsp\beta\nu E\bigr)=\codim\bigl(\Ran\mathcal A^{(F,\gamma)},\vtsp\convsp\gamma\nu F\bigr)$. By symmetry, Corollary \ref{EFlocincpositivever}, Proposition \ref{basicpropadmisssp} and the obvious inclusion $E^l\hookrightarrow E$ for any $l\ge0$, it suffices to prove this equality under the assumption that either $\beta>\gamma$ and we have a local inclusion $E_\loc\hookrightarrow F_\loc$, or $\beta\ge\gamma$ and we have a continuous inclusion $E\hookrightarrow F$. Combining this assumption with Remark \ref{easyincrem} and Proposition \ref{genformincforconsp}, we get a continuous inclusion $\consp\beta E\hookrightarrow\consp\gamma F$. 

\smallskip

Suppose $V$ is a subspace of $\convsp\beta\nu E$ with $V\cap\Ran\mathcal A^{(E,\beta)}=0$.

\begin{claim}[Claim (i): We have $V\cap\Ran\mathcal A^{(F,\gamma)}=0$.] If $f\in\;V\cap\Ran\mathcal A^{(F,\gamma)}\subseteq\convsp\beta\nu E\cap\convsp\gamma\nu F$ then $f=\mathcal Au$ for some $u\in\convsp\gamma\mu F$. Theorem \ref{changebetareggen} then gives $u\in\convsp\beta\mu E$ and so $f=\mathcal Au\in V\cap\Ran\mathcal A^{(E,\beta)}=0$, completing the claim.
\end{claim}

If $\codim\bigl(\Ran\mathcal A^{(E,\beta)},\vtsp\convsp\beta\nu E\bigr)=\infty$ then $\codim\bigl(\Ran\mathcal A^{(F,\gamma)},\vtsp\convsp\gamma\nu F\bigr)=\infty$ by an easy application of Claim (i). Now suppose $\codim\bigl(\Ran\mathcal A^{(E,\beta)},\vtsp\convsp\beta\nu E\bigr)<\infty$ and choose a finite dimensional space $V\subset\convsp\beta\nu E$ which satisfies
\begin{equation}
\label{eqnforcompinfcdcase1}
V+\Ran\mathcal A^{(E,\beta)}\;=\;\convsp\beta\nu E
\quad\text{and}\quad
V\cap\Ran\mathcal A^{(E,\beta)}\;=\;0.
\end{equation}
Claim (i), \eqref{eqnforcompinfcdcase1} and the inclusions $V\subset\convsp\beta\nu E\subseteq\convsp\gamma\nu F$ give
\begin{equation}
\label{eqnforcompinfcdcase2}
V+\Ran\mathcal A^{(F,\gamma)}\;\subseteq\;\convsp\gamma\nu F
\quad\text{and}\quad
V\cap\Ran\mathcal A^{(F,\gamma)}\;=\;0.
\end{equation}
The following claim thus completes the proof.

\begin{claim}[Claim (ii): We have $\convsp\gamma\nu F\subseteq V+\Ran\mathcal A^{(F,\gamma)}$.]
Let $f\in\convsp\gamma\nu F$ and choose $\lambda<\gamma$ so that $\lambda$ and $\gamma$ lie in the same component of $\R\!\setminus\!\Res{\mathcal A}$. By Lemma \ref{bestapproxresinnonsepcase} we can find a sequence $\{f_i\}_{i\in\N}\subset C^\infty_0$ such that $f_i\to f$ in $\convsp\lambda\nu{F^{-1}}$. Now \eqref{eqnforcompinfcdcase1} and the inclusions $\convsp\beta\mu E\subseteq\convsp\gamma\mu F\subseteq\convsp\lambda\mu{F^{-1}}$ give
\begin{equation}
\label{prelimincofCinft}
f_i\;\in\;C^\infty_0
\;\subset\;\convsp\beta\nu E
\;=\;V+\Ran\mathcal A^{(E,\beta)}
\;\subseteq\;V+\Ran\mathcal A^{(F^{-1},\lambda)}.
\end{equation}
However Theorem \ref{wfomres1} implies $\Ran\mathcal A^{(F^{-1},\lambda)}$ is a closed subspace of $\convsp\lambda\nu{F^{-1}}$ whilst $V\subset\convsp\beta\nu E\subseteq\convsp\lambda\nu{F^{-1}}$. Thus $V+\Ran\mathcal A^{(F^{-1},\lambda)}$ is also a closed subspace of $\convsp\lambda\nu{F^{-1}}$, so \eqref{prelimincofCinft} implies $f\in V+\Ran\mathcal A^{(F^{-1},\lambda)}$. In turn this means we can find $g\in V$ and $u\in\convsp\lambda\mu{F^{-1}}$ which satisfy $f=\mathcal Au+g$. Since $f-g\in\convsp\gamma\nu F$ Theorem \ref{changebetareggen} then gives $u\in\convsp\gamma\mu F$. Therefore $f\in V+\Ran\mathcal A^{(F,\gamma)}$, completing the claim.
\end{claim}
\end{proof}

\subsection{Adjoint operators and the Fredholm property}

Suppose $\mathcal A$ is an admissible elliptic operator on $\R^n$ of order $(\mu,\nu)$ and let $\mathcal A^*$ denote the formal adjoint of $\mathcal A$ (with respect to the Lebesgue measure on $\R^n$). Also let $m$, $\overline{\mu}$ and $\overline{\nu}$ by defined as in \eqref{defofadjorders}. It follows easily from the definition of an admissible operator that $\mathcal A^*$ is a $k\negvtsp\times\negvtsp k$ system of differential operators on $\R^n$ or order $(\overline{\mu},\overline{\nu})$ which satisfies Condition (i) of Definition \ref{defnofadellop}. Furthermore $\detsym{\mathcal A^*}(x,\xi)=\overline{\detsym{\mathcal A}(x,\xi)}$ for all $(x,\xi)\in\R^n\times\R^n$, from which it is clear that $\mathcal A^*$ also satisfies Condition (ii) of Definition \ref{defnofadellop}. Hence $\mathcal A^*$ is an admissible elliptic operator on $\R^n$ of order $(\overline{\mu},\overline{\nu})$. 

Let $\mathcal A_0$ denote the principal part of $\mathcal A$ and let $\mathcal B(\omega,\dfo_t,\dfo_\omega)$ denote the model elliptic operator on $\cly$ associated to $\mathcal A_0$ by \eqref{defofunifellondon}; that is
\begin{equation}
\label{hwgahwga1}
\mathcal A_0(x,\dfo_x)
\,=\,r^\nu\mathcal B(\omega,r\dfo_r,\dfo_\omega)r^{-\mu}
\end{equation}
whilst $\mathfrak B_{\mathcal A}$, the operator pencil associated to $\mathcal A$, is defined by
\begin{equation}
\label{hwgahwga2}
\mathfrak B_{\mathcal A}(\lambda)(\omega,\dfo_\omega)
\,=\,\mathcal B(\omega,\lambda,\dfo_\omega).
\end{equation}
Now it is straightforward to check that the principal part of $\mathcal A^*$ is simply $\mathcal A_0^*$, the formal adjoint of $\mathcal A_0$ (with respect to the Lebesgue measure on $\don$). Let $\mathcal B^*$ denote the formal adjoint of $\mathcal B$ (with respect to the measure $dt\d S^{n-1}$ on $\cly$). Since the Lebesgue measure on $\don$ can be written as $d^nx=r^{n-1}\d r\d S^{n-1}$, we have
\[
(r\dfo_r)^*
\;=\;(r^{n-1})^{-1}\dfo_r\bigl(rr^{n-1}\,\cdot\,\bigr)
\;=\;r\dfo_r-in
\]
with respect to this measure. Combined with \eqref{effofmovrmmth} and \eqref{hwgahwga1} we now get
\begin{eqnarray*}
\mathcal A_0^*(x,\dfo_x)
&=&r^{-\mu}\mathcal B^*\bigl(\omega,r\dfo_r-in,\dfo_\omega\bigr)r^\nu\\
&=&r^m r^{-\mu}\mathcal B^*\bigl(\omega,r\dfo_r-i(n+m),\dfo_\omega\bigr)r^{-m}r^\nu\\
&=&r^{\overline{\nu}}\mathcal B^*\bigl(\omega,r\dfo_r-i(n+m),\dfo_\omega\bigr)
r^{-\overline{\mu}}\,.
\end{eqnarray*}
Thus the model operator on $\cly$ associated to $\mathcal A_0^*$ by \eqref{defofunifellondon} is $\mathcal B^*(\omega,\dfo_t-i(n+m),\dfo_\omega)$, and so the operator pencil associated to $\mathcal A^*$ is given by
\begin{equation}
\label{hwgahwga3}
\mathfrak B_{\mathcal A^*}(\lambda)(\omega,\dfo_\omega)
\,=\,\mathcal B^*\bigl(\omega,\lambda-i(n+m),\dfo_\omega\bigr).
\end{equation}
Taken together \eqref{hwgahwga2} and \eqref{hwgahwga3} complete the proof of the following result.

\begin{prop}
\label{relofadjfspec}
The operator pencils $\mathfrak B_{\mathcal A}$ and $\mathfrak B_{\mathcal A^*}$ associated to $\mathcal A$ and $\mathcal A^*$ respectively satisfy the relationship
\[
\mathfrak B_{\mathcal A^*}(\lambda)
\;=\;\bigl(\mathfrak B_{\mathcal A}(\,\overline{\lambda}+i(n+m))\bigr)^*
\]
for all $\lambda\in\C$, where the right hand side denotes the formal adjoint of the operator $\mathfrak B_{\mathcal A}(\,\overline{\lambda}+i(n+m))$ on $S^{n-1}$. It follows that $\lambda\in\fspec{\mathfrak B_{\mathcal A^*}}$ iff $\overline{\lambda}+i(n+m)\in\fspec{\mathfrak B_{\mathcal A}}$ with full agreement of geometric, algebraic and partial algebraic multiplicities. In particular 
\[
\Res{\mathcal A^*}\,=\,(n+m)-\Res{\mathcal A}.
\]
\end{prop}

\medskip

Let $E$ be an admissible space and set $F=E_0^*$. By Propositions \ref{dualpropforconsp} and \ref{basicpropadmisssp} we know that $F$ is also an admissible space whilst
\begin{eqnarray*}
&&(\convsp\beta\mu{E_0})^*
\,=\,\prod_{i=1}^k\consp{n-\beta+\mu_i}F^{-\mu_i}
\,=\,\convsp{n+m-\beta}{\overline{\nu}}{F^{-m}}\\
\text{and}&&(\convsp\beta\nu E_0)^*
\,=\,\prod_{i=1}^k\consp{n-\beta+\nu_i}F^{-\nu_i}
\,=\,\convsp{n+m-\beta}{\overline{\mu}}{F^{-m}}
\end{eqnarray*}
for any $\beta\in\R$. Therefore the adjoint of the map $\mathcal A^{(E_0,\beta)}:\convsp\beta\mu{E_0}\to\convsp\beta\nu{E_0}$ is the map
\begin{equation}
\label{adjofAbetaE}
\bigl(\mathcal A^{(E_0,\beta)}\bigr)^*
\,=\,(\mathcal A^*)^{(F^{-m}\negvtsp,\vtsp n+m-\beta)}
\,:\,\convsp{n+m-\beta}{\overline{\mu}}{F^{-m}}\to\convsp{n+m-\beta}{\overline{\nu}}{F^{-m}},
\quad F=E_0^*.
\end{equation}

\begin{rem}
Let $\beta\in\R$. By Theorem \ref{wfomres1} we have that $\mathcal A^{(E_0,\beta)}$ is semi-Fredholm provided $\beta\notin\Res{\mathcal A}$ and $(\mathcal A^*)^{(F^{-m}\negvtsp,\vtsp n+m-\beta)}$ is semi-Fredholm provided $n+m-\beta\notin\Res{\mathcal A^*}$. These conditions can be seen to be equivalent by Proposition \ref{relofadjfspec}, which is expected since the maps $\mathcal A^{(E_0,\beta)}$ and $(\mathcal A^*)^{(F^{-m}\negvtsp,\vtsp n+m-\beta)}$ are each others adjoints.
\end{rem}

\bigskip

\begin{thm}
\label{sfomres1}
Suppose $\mathcal A$ is an admissible elliptic operator on $\R^n$ of order $(\mu,\nu)$, let $E$ be an admissible space and choose $\beta\in\R\!\setminus\!\Res{\mathcal A}$. Then the map $\mathcal A:\convsp\beta\mu E\to\convsp\beta\nu E$ is Fredholm.
\end{thm}

\begin{proof}
Since Theorem \ref{wfomres1} shows that the map $\mathcal A^{(E,\beta)}$ is semi-Fredholm we simply have to establish that $\Index\mathcal A^{(E,\beta)}$ is finite. Furthermore, Theorem \ref{stabofindex} shows this index to be stable under changes in $E$. It thus suffices to prove the result for a particular admissible space $E$; we take this to be $E=\Sob 20$ so $E=E_0=E_0^*$. 

By \eqref{adjofAbetaE} we have $\bigl(\mathcal A^{(E,\beta)}\bigr)^*=(\mathcal A^*)^{(E^{-m}\negvtsp,\vtsp n+m-\beta)}$. Furthermore Proposition \ref{relofadjfspec} implies $n+m-\beta\notin\Res{\mathcal A^*}$. Theorem \ref{wfomres1} then shows the maps $\mathcal A^{(E,\beta)}$ and $\bigl(\mathcal A^{(E,\beta)}\bigr)^*$ are both semi-Fredholm with finite dimensional kernels. However we have the general identity
\[
\codim\bigl(\Ran\mathcal A^{(E,\beta)},\vtsp\convsp\beta\nu E\bigr)
\;=\;\dim\Ker\bigl(\mathcal A^{(E,\beta)}\bigr)^*
\]
for semi-Fredholm operators (see Theorem IV.5.13 in \cite{Ka} for example). Therefore
\begin{multline*}
\Index\mathcal A^{(E,\beta)}
\;=\;\dim\Ker\mathcal A^{(E,\beta)}
-\codim\bigl(\Ran\mathcal A^{(E,\beta)},\vtsp\convsp\beta\nu E\bigr)\\
\;=\;\dim\Ker\mathcal A^{(E,\beta)}-\dim\Ker(\mathcal A^*)^{(E^{-m}\negvtsp,\vtsp n+m-\beta)},
\end{multline*}
which is finite.
\end{proof}

\subsection{The change in index formula}

This section is devoted to establishing the following result which shows how the index of the Fredholm map $\mathcal A:\convsp\beta\mu E\to\convsp\beta\nu E$ varies when we change $\beta$ and $E$ over a greater range than is permitted in Theorem \ref{stabofindex}

\begin{thm}
\label{mressec5}
Suppose $E$ is an admissible space and $\beta_1,\beta_2\in\R\!\setminus\!\Res{\mathcal A}$ with $\beta_1\le\beta_2$. Set $\Sigma=\{\lambda\in\fspec{\mathfrak B_{\mathcal A}}\,|\,\im\lambda\in[\beta_1,\beta_2]\}$ and, for each $\lambda\in\Sigma$, let $m_\lambda$ denote the algebraic multiplicity of $\lambda$. Then we have 
\[
\Index\mathcal A^{(\beta_1)}
\;=\;\Index\mathcal A^{(\beta_2)}+\sum_{\lambda\in\Sigma}m_\lambda.
\]
\end{thm}

By Remark \ref{secremaboutranops} and Lemma \ref{formofdecayofdiff} $\mathcal A$ can be approximated arbitrarily closely (simultaneously in $\mathscr L\bigl(\convsp{\beta_j}\mu E,\vtsp\convsp{\beta_j}\nu E\bigr)$ for $j=1,2$) by another admissible elliptic operator $\mathcal A'$ whose coefficients agree with those of $\mathcal A_0$ on a neighbourhood of $\infty$. Since the set of Fredholm operators of a given index is open (in operator norm) and $\mathfrak B_{\mathcal A_0}$ is the spectral pencil associated to both $\mathcal A$ and $\mathcal A'$, it suffices to prove the result assuming the coefficients of $\mathcal A$ and $\mathcal A_0$ agree on a neighbourhood of $\infty$. Choose $I\in\N$ so that 
\begin{equation}
\label{nowdishyp}
\zeta_I(\mathcal A-\mathcal A_0)\,=\,0. 
\end{equation}

\bigskip

Set $M=\sum_{\lambda\in\Sigma}m_\lambda$ and let $\{w_1,\dots,w_M\}$ denote any basis of the vector space $X_\Sigma$ given by Theorem \ref{mressec22} for the operator $\mathcal A_0$. 

\begin{lem}
\label{wtwsautdtcfs2ts1}
Let $f\in\convsp{\beta_2}\nu E$ and suppose $\mathcal Au=f$ for some $u\in\convsp{\beta_1}\mu E$. Then
\begin{equation}
\label{largeRasyform}
u\,=\,v+\zeta_I\sum_{j=1}^M z_j w_j
\end{equation}
for some $v\in\convsp{\beta_2}\mu E$ and $z_1,\dots,z_M\in\C$.
\end{lem}

\begin{proof}
By Lemma \ref{equivofnorms1} $\zeta_{I-1} u\in\donvsp{\beta_1}\mu E$. Setting $g=\mathcal A_0(\zeta_{I-1}u)$ it follows that $g\in\donvsp{\beta_1}\nu E$ and $g=0$ on a neighbourhood of 0. On the other hand \eqref{nowdishyp} gives
\[
\zeta_I g
\;=\;\zeta_I\mathcal A_0(\zeta_{I-1} u)
\;=\;\zeta_I\mathcal Au
\;=\;\zeta_I f.
\]
Thus Lemma \ref{equivofnorms1} and Remark \ref{newequivofnormsrem2} imply $g\in\donvsp{\beta_2}\nu E$.

For $i=1,2$ set $u_i=\bigl(\mathcal A_0^{(\beta_i)}\bigr)^{-1}\!g\in\donvsp{\beta_i}\mu E$. By Theorem \ref{mressec22} it follows that 
\begin{equation}
\label{largeRasyforminter}
u_1-u_2
\;=\;\sum_{j=1}^M z_j w_j
\end{equation}
for some $z_1,\dots,z_M\in\C$. Now $u_2\in\donvsp{\beta_2}\mu E$ so, setting $v=\eta_I u+\zeta_I u_2$, we get $v\in\convsp{\beta_2}\mu E$ by Lemma \ref{equivofnorms1}. On the other hand $u_1=\zeta_{I-1}u$ (this follows from the definition of $g$) so $\zeta_I u_1=\zeta_I\zeta_{I-1} u=\zeta_I u$. Combining this with the definition of $v$ and \eqref{largeRasyforminter} we get \eqref{largeRasyform}.
\end{proof}

We have $\beta_1\le\beta_2$ so $\convsp{\beta_2}\mu E\subseteq\convsp{\beta_1}\mu E$ (see \eqref{cntsincforXE}) and hence $\Ker\mathcal A^{(\beta_2)}\subseteq\Ker\mathcal A^{(\beta_1)}$. Since both kernels are finite dimensional we can therefore choose $d\in\NZ$ and $u_1,\dots,u_d\in\Ker\mathcal A^{(\beta_1)}\subset\convsp{\beta_1}\mu E$ such that $u_1,\dots,u_d$ are linearly independent over $\convsp{\beta_2}\mu E$ and 
\begin{equation}
\label{diffofkers}
\Ker\mathcal A^{(\beta_1)}
\;=\;\Ker\mathcal A^{(\beta_2)}+\Span\{u_1,\dots,u_d\}.
\end{equation}
In particular
\begin{equation}
\label{diffofkerdims}
\dim\Ker\mathcal A^{(\beta_1)}
\;=\;\dim\Ker\mathcal A^{(\beta_2)}+d.
\end{equation}

By Lemma \ref{wtwsautdtcfs2ts1} we can find ${z}_{ij}\in\C$ for $i\in\{1,\dots,d\}$ and $j\in\{1,\dots,M\}$ such that 
\[
u_i\ \equiv\ \zeta_I\sum_{j=1}^M z_{ij} w_j
\pmod{\convsp{\beta_2}\mu E}
\]
for $i=1,\dots,d$. Furthermore, using the linear independence of $u_1,\dots,u_d$ over $\convsp{\beta_2}\mu E$, the $d\!\times\!M$ matrix with entries ${z}_{ij}$ must have rank $d$. Therefore $d\le M$ and we can choose the basis $\{w_1,\dots,w_M\}$ of $X_\Sigma$ so that
\begin{equation}
\label{nfotrbtusatws}
u_i\ \equiv\ \zeta_I w_i\pmod{\convsp{\beta_2}\mu E}
\end{equation}
for $i=1,\dots,d$. We will stick to such a choice of $\{w_1,\dots,w_M\}$ for the remainder of this section. 

\bigskip

Let $Y=\Ran\mathcal A^{(\beta_1)}\cap\convsp{\beta_2}\nu E$. 

\begin{lem}
\label{lemregY}
The space $Y$ is closed in $\convsp{\beta_2}\nu E$ and satisfies
\[
\codim(Y,\,\convsp{\beta_2}\nu E)=\codim(\Ran\mathcal A^{(\beta_1)},\,\convsp{\beta_1}\nu E).
\]
\end{lem}

\begin{proof}
Since $\convsp{\beta_2}\mu E\subseteq\convsp{\beta_1}\mu E$ we get $\Ran\mathcal A^{(\beta_2)}\subseteq Y\subseteq\convsp{\beta_2}\nu E$. However $\Ran\mathcal A^{(\beta_2)}$ is closed in $\convsp{\beta_2}\nu E$ with finite codimension (by Theorem \ref{sfomres1}) so the same must be true for $Y$. Now choose a finite dimensional space $V\subset\convsp{\beta_2}\nu E$ so that $\convsp{\beta_2}\nu E=V+Y$ and $V\cap Y=0$. 

\begin{claim}[Claim: $\convsp{\beta_1}\nu E=V+\Ran\mathcal A^{(\beta_1)}$.]
The inclusion $V\subset\convsp{\beta_2}\nu E\subseteq\convsp{\beta_1}\nu E$ immediately gives the reverse inclusion. Now let $f\in\convsp{\beta_1}\nu E$. Put $v=\zeta_{I-1}\bigl(\mathcal A_0^{(\beta_1)}\bigr)^{-1}(\zeta_I f)\in\convsp{\beta_1}\mu E$ and $g=f-\mathcal Av\in\convsp{\beta_1}\nu E$ (n.b.~Lemma \ref{equivofnorms1} and Theorem \ref{mressec21} ensure that $v$ is well defined). With the help of \eqref{nowdishyp} we thus have
\[
\zeta_{I+1}g
\,=\,\zeta_{I+1}f-\zeta_{I+1}\mathcal A_0
\bigl(\zeta_{I-1}\bigl(\mathcal A_0^{(\beta_1)}\bigr)^{-1}(\zeta_I f)\bigr)
\,=\,0.
\]
Hence $g\in\convsp{\beta_2}\nu E$ by Remark \ref{newequivofnormsrem2}. On the other hand, by the definition of $V$, $g=\mathcal Aw+h$ for some $w\in\convsp{\beta_1}\mu E$ and $h\in V$. Setting $u=v+w\in\convsp{\beta_1}\mu E$ we then have 
\[
f
\,=\,g+\mathcal Av
\,=\,h+\mathcal Av+\mathcal Aw
\,=\,h+\mathcal Au
\,\in\,V+\Ran\mathcal A^{(\beta_1)}.
\]
Therefore $\convsp{\beta_1}\mu E\subseteq V+\Ran\mathcal A^{(\beta_1)}$, completing the claim.
\end{claim}

Now $V\subset\convsp{\beta_2}\mu E$ so
\[
V\cap\Ran\mathcal A^{(\beta_1)}
\;=\;V\vtsp\cap\vtsp\convsp{\beta_2}\mu E\vtsp\cap\vtsp\Ran\mathcal A^{(\beta_1)}
\;=\;V\cap Y
\;=\;0.
\]
Combined with the above claim this now completes the result.
\end{proof}

From Lemma \ref{lemregY} and the observation that $\Ran\mathcal A^{(\beta_2)}\subseteq Y\subseteq\convsp{\beta_2}\nu E$ we get
\begin{eqnarray}
\lefteqn{\codim(\Ran\mathcal A^{(\beta_1)},\,\convsp{\beta_1}\nu E)
\;=\;\codim(Y,\,\convsp{\beta_2}\nu E)}
\nonumber\\
\label{diffofcokerdims}
&&\qquad\quad
=\;\codim(\Ran\mathcal A^{(\beta_2)},\,\convsp{\beta_2}\nu E)-\codim(\Ran\mathcal A^{(\beta_2)},\,Y).
\end{eqnarray}

\begin{proof}[Proof of Theorem \ref{mressec5}]
For $i=1,\dots,M-d$ set $f_i=\mathcal A(\zeta_I w_{i+d})\in\Ran\mathcal A^{(\beta_1)}\subseteq\convsp{\beta_1}\nu E$. Since $\zeta_{I+1}\mathcal A(\zeta_I w_{i+d})=\zeta_{I+1}\mathcal A_0w_{i+d}=0$ (by \eqref{nowdishyp} and Theorem \ref{mressec22}; recall that $w_{i+d}\in X_\Sigma$) we immediately get $f_i\in\convsp{\beta_2}\nu E$ using Remark \ref{newequivofnormsrem2}. Now let $W=\Span\{f_1,\dots,f_{M-d}\}\subset Y$. 

\begin{claim}[Claim (i): $f_1,\dots,f_{M-d}$ are linearly independent over $\Ran\mathcal A^{(\beta_2)}$.]
Suppose 
\[
\sum_{i=1}^{M-d}{z}_{i+d} f_i
\;=\;\mathcal Au
\]
for some $u\in\convsp{\beta_2}\mu E$ and ${z}_{d+1},\dots,{z}_M\in\C$. Thus
\[
\zeta_I\sum_{i=1}^{M-d}{z}_{i+d} w_{i+d}\;-\;u
\ \in\ \Ker\mathcal A^{(\beta_1)}
\]
and so, by \eqref{diffofkers}, there exists $v\in\Ker\mathcal A^{(\beta_2)}$ and ${z}_1,\dots,{z}_d$ such that 
\[
\sum_{i=1}^d{z}_i u_i\,+\,\zeta_I\sum_{i=1}^{M-d}{z}_{i+d} w_{i+d}
\ =\ u+v.
\]
Now $u+v\in\convsp{\beta_2}\mu E$ so \eqref{nfotrbtusatws} implies $\zeta_I w\in\convsp{\beta_2}\mu E$ where $w=\sum_{i=1}^M {z}_i w_i$. Lemma \ref{equivofnorms1} then implies $\zeta_I w\in\donvsp{\beta_2}\mu E$. On the other hand $w\in X_\Sigma$ so $(1-\zeta_I)w\in\donvsp{\beta_2}\mu E$ and $\mathcal A_0w=0$ by Theorem \ref{mressec22}. However $\mathcal A_0:\donvsp{\beta_2}\mu E\to\donvsp{\beta_2}\nu E$ is an isomorphism by Theorem \ref{mressec21}. Thus $w=\sum_{i=1}^M {z}_i w_i=0$. The fact that $\{w_1,\dots,w_M\}$ is a basis for $X_\Sigma$ now completes Claim (i).
\end{claim}

\begin{claim}[Claim (ii): $Y=W+\Ran\mathcal A^{(\beta_2)}$.]
The inclusion $W+\Ran\mathcal A^{(\beta_2)}\subseteq Y$ is trivial. Now let $f\in Y$ so $f\in\convsp{\beta_2}\nu E$ and $f=\mathcal Au$ for some $u\in\convsp{\beta_1}\mu E$. By Lemma \ref{wtwsautdtcfs2ts1} and \eqref{nfotrbtusatws} we can thus find $v\in\convsp{\beta_2}\mu E$ and ${z}_1,\dots,{z}_M\in\C$ such that 
\[
u\;=\;
v+\sum_{i=1}^d {z}_i u_i
\,+\,\sum_{i=1}^{M-d} {z}_{i+d}\,\zeta_I w_{i+d}\,.
\]
Using the fact that $u_1,\dots,u_d\in\Ker\mathcal A^{(\beta_1)}$ we then get
\[
f\ =\ \mathcal Au
\ =\ \mathcal Av+\sum_{i=1}^{M-d} {z}_{i+d} \mathcal A(\zeta_I w_{i+d})
\ =\ \mathcal Av+\sum_{i=1}^{M-d} {z}_{i+d} f_i
\ \in\ \Ran\mathcal A^{(\beta_2)}+W.
\]
Hence $Y\subseteq W+\Ran\mathcal A^{(\beta_2)}$, completing Claim (ii).
\end{claim}

By Claim (i) we have $\dim W=M-d$ and $W\cap\Ran\mathcal A^{(\beta_2)}=0$. Combined with Claim (ii) it follows that 
\begin{equation}
\label{leqnwnfciires}
\codim(\Ran\mathcal A^{(\beta_2)},\,Y)
\;=\;M-d.
\end{equation}
By combining \eqref{diffofkerdims}, \eqref{diffofcokerdims} and \eqref{leqnwnfciires} we finally get 
\begin{eqnarray*}
\Index\mathcal A^{(\beta_1)}
&=&\dim\Ker\mathcal A^{(\beta_1)}-\codim(\Ran\mathcal A^{(\beta_1)},\,\convsp{\beta_1}\nu E)\\
&=&\dim\Ker\mathcal A^{(\beta_2)}+d-\codim(\Ran\mathcal A^{(\beta_2)},\,\convsp{\beta_2}\nu E)+M-d\\
&=&\Index\mathcal A^{(\beta_2)}+M,
\end{eqnarray*}
completing the result. 
\end{proof}

\subsection{Self-adjoint operators}
Let $\mathcal A$ be an admissible elliptic operator on $\R^n$ of order $(\mu,\nu)$. If we can determine the set $\fspec{\mathfrak B_{\mathcal A}}$ together with the algebraic multiplicities of each point $\lambda\in\fspec{\mathfrak B_{\mathcal A}}$ Theorems \ref{stabofindex} and \ref{mressec5} allow us to compute the change in the index of the map $\mathcal A:\convsp\beta\mu E\to\convsp\beta\nu E$ whenever we change the admissible space $E$ or the parameter $\beta$ within the set $\R\!\setminus\!\Res{\mathcal A}$. In turn this allows us to compute the actual index for any such pair $(E,\beta)$ provided we know the index for one pair. In general this requires further computation; however if $\mathcal A$ is formally self-adjoint we can use symmetry to compute the index of $\mathcal A^{(E,\beta)}$ from knowledge of the spectrum of $\mathfrak B_{\mathcal A}$ alone.

\begin{thm}
\label{selfadjthm}
Suppose $\mathcal A$ is a formally self-adjoint admissible elliptic operator of order $(\mu,\nu)$. Then $\Res{\mathcal A}$ is symmetric about $(n+m)/2$. Furthermore, if $E$ is an admissible space and $\beta\in\R\!\setminus\!\Res{\mathcal A}$, then 
\begin{equation}
\label{symmofindex1}
\Index\mathcal A^{(n+m-\beta)}
\,=\,-\Index\mathcal A^{(\beta)}.
\end{equation}
In particular we either have $(n+m)/2\notin\Res{\mathcal A}$, in which case $\Index\mathcal A^{(n/2+m/2)}=0$, or $(n+m)/2\in\Res{\mathcal A}$, in which case the sum of the algebraic multiplicities of those $\lambda\in\fspec{\mathfrak B_{\mathcal A}}$ with $\im\lambda=(n+m)/2$ is even (say 2d for some $d\in\N$) and 
\[
\Index\mathcal A^{(n/2+m/2-\epsilon)}
\;=\;d\;=\;-\Index\mathcal A^{(n/2+m/2+\epsilon)}
\]
for all sufficiently small $\epsilon>0$. 
\end{thm}

\begin{proof}
For the sake of convenience we set $l=(n+m)/2$. The symmetry of $\Res{\mathcal A}$ about $l$ follows directly from Proposition \ref{relofadjfspec}. Now let $E$ be an admissible space and $\beta\in\R\!\setminus\!\Res{\mathcal A}$. Using \eqref{adjofAbetaE} and the assumption $\mathcal A^*=\mathcal A$, we have that the adjoint of the map $\mathcal A^{(E_0,\beta)}$ is $\mathcal A^{(F^{-m}\negvtsp,\vtsp 2l-\beta)}$ where $F=E_0^*$. Using Theorem \ref{stabofindex} together with the fact that the index of a Fredholm map changes sign when we take its adjoint (see Corollary IV.5.14 in \cite{Ka} for example) we now get 
\[
\Index\mathcal A^{(E,2l-\beta)}
\;=\;\Index\mathcal A^{(F^{-m}\negvtsp,\vtsp 2l-\beta)}
\;=\;-\Index\mathcal A^{(E_0,\beta)}
\;=\;-\Index\mathcal A^{(E,\beta)},
\]
which establishes \eqref{symmofindex1}.

If $l\notin\Res{\mathcal A}$ then we get $\Index\mathcal A^{(l)}=0$ by setting $\beta=l$ in \eqref{symmofindex1}. Now suppose $l\in\Res{\mathcal A}$. Since $\Res{\mathcal A}$ consists of isolated points (see Remark \ref{fspecthmrem}) it follows that we can find some $\delta>0$ such that $(l,l\!+\!\delta)\cap\Res{\mathcal A}=\emptyset$. Theorems \ref{stabofindex} and \ref{sfomres1} together with \eqref{symmofindex1} then imply the existence of $d\in\Z$ such that
\[
\Index\mathcal A^{(l-\epsilon)}
\;=\;d\;=\;-\Index\mathcal A^{(l+\epsilon)}
\]
for all $0<\epsilon<\delta$. By Theorem \ref{mressec5} we also know
\[
2d\;=\;\Index\mathcal A^{(l-\epsilon)}-\Index\mathcal A^{(l+\epsilon)}
\;=\;\sum_{\lambda\in\Sigma} m_\lambda,
\]
where $\Sigma=\bigl\{\lambda\in\fspec{\mathfrak B_{\mathcal A}}\,\big|\,\im\lambda=l\bigr\}$ and $m_\lambda$ is the algebraic multiplicity of a given $\lambda\in\Sigma$. This completes the result. 
\end{proof}

\subsection{Homogeneous operators with constant coefficients}

Let $\mathcal A$ be a $k\times k$ system of differential operators on $\R^n$ of order $(\mu,\nu)$. We say that $\mathcal A$ is a \emph{homogeneous constant coefficient operator} if, for each $i,j\in\{1,\dots,k\}$, the $ij$th entry of $\mathcal A$ is a constant coefficient scalar differential operator which is homogeneous of order $\mu_j-\nu_i$ (or equal to 0 if $\mu_j-\nu_i<0$).

\begin{rem}
\label{thisremnousedhere}
Suppose $\mathcal A$ is an elliptic homogeneous constant coefficient operator on $\R^n$ of order $(\mu,\nu)$. It follows immediately that $\mathcal A$ is uniformly elliptic on $\R^n$. On the other hand $A_{ij}(x,\dfo_x)$, the $ij$th entry of $\mathcal A$, is a constant coefficient scalar differential operator which is homogeneous of order $\mu_j-\nu_i$. Thus $A_{ij}(x,\dfo_x)$ is an admissible scalar operator with principal part $A_{ij}(x,\dfo_x)$ (see Definition \ref{earlydefofadmiss}). It follows that $\mathcal A$ is an admissible elliptic operator on $\R^n$ with principal part $\mathcal A_0=\mathcal A$. 
\end{rem}

\medskip

A straightforward application of Theorems \ref{stabofindex} and \ref{mressec5} allows us to generalise an index formula given in \cite{LM} to cover arbitrary admissible spaces. Before giving the result we need to introduce the following notation. For any $\beta\in\R$ let $P_n(\beta)$ denote the dimension of the set of polynomials in $n$ variables whose degree does not exceed $\beta$ (with $P_n(\beta)=0$ when $\beta<0$). Thus
\begin{equation}
\label{whthfnPnisdc1}
P_n(\beta)
\;=\;\begin{cases}
P_n(l)&\text{if $\beta\in[l,l+1)$ for some $l\in\NZ$,}\\
0&\text{if $\beta<0$.}
\end{cases}
\end{equation}
For any $\mu,\nu\in\NZ^k$ and $\beta\in\R$ we set 
\[
P_n^{\mu,\nu}\!(\beta)
\;=\;\sum_{i=1}^k\bigl(P_n(-\beta+\mu_i)-P_n(-\beta+\nu_i)
-P_n(\beta-\nu_i-n)+P_n(\beta-\mu_i-n)\bigr)\vtsp.
\]
It is clear from \eqref{whthfnPnisdc1} that the function $P_n^{\mu,\nu}$ is constant on the components of $\R\!\setminus\!\Z$.

\begin{thm}
\label{thethmihaflofse}
Suppose $\mathcal A$ is an admissible elliptic operator on $\R^n$ of order $(\mu,\nu)$ whose principal part is a homogeneous constant coefficient operator. Then $\Res{\mathcal A}=\bigl\{\beta\in\R\,\big|\,P_n^{\mu,\nu}\!(\beta^-)\ne P_n^{\mu,\nu}\!(\beta^+)\bigr\}\subseteq\Z$ and $\Index\mathcal A^{(E,\beta)}\negvtsp=\negvtsp P_n^{\mu,\nu}\!(\beta)$ for any admissible space $E$ and $\beta\negvtsp\in\negvtsp\R\negvtsp\setminus\negvtsp\Res{\mathcal A}$. In particular $\Res{\mathcal A}\cap(m,n)=\emptyset$ and $\Index\mathcal A^{(E,\beta)}=0$ for any $\beta\in(m,n)$, where $m=\max_i\mu_i$.
\end{thm}

\begin{proof}
By Theorems 3 and 4 in \cite{LM} and Remark \ref{actidentifofequivnorms} above we know that $\mathcal A^{(\Sob20\!,\beta)}$ is a Fredholm map with index $P_n^{\mu,\nu}\!(\beta)$ whenever $\beta\notin\Z$. On the other hand $\Res{\mathcal A}$ is a discrete subset of $\R$ (see Remark \ref{fspecthmrem}) whilst Theorems \ref{stabofindex} and \ref{mressec5} imply the function $\beta\mapsto\Index\mathcal A^{(E,\beta)}$ is independent of the admissible space $E$, is constant on the components of $\R\!\setminus\!\Res{\mathcal A}$, and satisfies
\[
\Index\mathcal A^{(E,\beta^-)}
\;\ne\;\Index\mathcal A^{(E,\beta^+)}
\]
for all $\beta\in\Res{\mathcal A}$. It follows that $\Res{\mathcal A}$ is precisely the set of points at which $P_n^{\mu,\nu}$ is discontinuous, whilst $\Index\mathcal A^{(E,\beta)}=P_n^{\mu,\nu}\!(\beta)$ for any admissible space $E$ and $\beta\in\R\!\setminus\!\Res{\mathcal A}$. 

If $\nu_i>m$ for some $i$ then the $i$th row of the matrix operator $\mathcal A$ would be 0. As this contradicts the assumption that $\mathcal A$ is elliptic we get $\max_i\nu_i\le m$. Coupled with \eqref{whthfnPnisdc1} and the fact that $\mu$ and $\nu$ are vectors of non-negative integers, it is clear that all the terms in the sum defining $P_n^{\mu,\nu}\!(\beta)$ are equal to 0 when $\beta\in(m,n)$. This completes the result.
\end{proof}

Theorem \ref{thethmihaflofse} can be made more explicit for several special classes of operators; these classes include Dirac type operators and the Laplacian.

\begin{thm}
\label{corofthethmihaflofse}
Suppose $\mathcal A$ is an admissible elliptic operator on $\R^n$ of order $(\mu,\nu)$ whose principal part is a homogeneous constant coefficient operator. Also suppose $n\ge2$, $\nu_1=\dots=\nu_k=0$ and $\mu_1=\dots=\mu_k=m$ for some $m\in\{1,2\}$. Then $\Res{\mathcal A}=\Z\!\setminus\!\{m+1,\dots,n-1\}$ (with $\Res{\mathcal A}=\Z$ if $n<m+2$). Furthermore, for any admissible space $E$ and $\beta\in\R\!\setminus\!\Res{\mathcal A}$, 
\begin{equation}
\label{explicitfor1}
\Index\mathcal A^{(E,\beta)}
\;=\;-\vtsp\frac k{(n-1)!}\,(l-1)\prod_{j=2}^{n-1}\vtsp\abs{l-j}
\end{equation}
if $m=1$ and 
\begin{equation}
\label{explicitfor2}
\Index\mathcal A^{(E,\beta)}
\;=\;-\vtsp\frac k{(n-1)!}\,(2l-n-1)\prod_{j=2}^{n-1}\vtsp\abs{l-j}
\end{equation}
if $m=2$, where $l\in\Z$ is chosen so that $\beta\in[l,l+1)$ and the product terms are defined to be equal to 1 if $n=2$. In particular, $\Index\mathcal A^{(E,\beta)}=0$ for any $\beta\in(m,n)$. 
\end{thm}

\begin{rem}
The special form of $\mu$ and $\nu$ means the domain and codomain of the map $\mathcal A^{(E,\beta)}$ are given as 
\[
\convsp\beta\mu E
\;=\;\prod_{i=1}^k\consp{\beta-m}{E^m}
\ \quad\text{and}\ \quad
\convsp\beta\nu E
\;=\;\prod_{i=1}^k\consp\beta E
\]
respectively.
\end{rem}

\begin{proof}[Proof of Theorem \ref{corofthethmihaflofse}]
We have $P_1(l)=l+1$ for any $l\in\NZ$ and $P_n(0)=1$ for any $n\in\N$. On the other hand, the set of polynomials in $(n+1)$ variables which are homogeneous of degree $l$ can be seen to have dimension $P_n(l)$. Hence $P_{n+1}(l)=P_{n+1}(l-1)+P_n(l)$ for all $n,l\in\N$. The resulting recurrence relations have the unique solution
\begin{equation}
\label{whthfnPnisdc2}
P_n(l)
\,=\,\frac{(l+n)!}{n!\, l!}
\end{equation}
for all $n\in\N$ and $l\in\NZ$. A straightforward calculation using \eqref{whthfnPnisdc1} and \eqref{whthfnPnisdc2} now shows that if $\beta\in(l,l+1)$ for some $l\in\Z$ then $P_n^{\mu,\nu}\!(\beta)$ is given explicitly by right hand side of \eqref{explicitfor1} or \eqref{explicitfor2} when $m=1$ or $2$ respectively. The result then follows from Theorem \ref{thethmihaflofse} and the fact that right hand sides of \eqref{explicitfor1} and \eqref{explicitfor2} are non-increasing functions of $l\in\Z$ which are strictly decreasing outside $\{m+1,\dots,n-1\}$.
\end{proof}

\section*{Acknowledgements}
The author wishes to thank D.~E.~Edmunds and V.~A.~Kozlov for several useful comments, especially those regarding existing results. This research was supported by EPSRC under grant GR/M21720.

\end{document}